\documentclass[journal]{IEEEtran}
\ifCLASSINFOpdf
\else
\fi
\usepackage{def}

\begin{document}
%
\title{Implicit Tracking-Based Distributed Constraint-Coupled Optimization}
%
%
%
%
%
\author{Jingwang Li and Housheng Su
    \thanks{The authors are with the School of Artificial Intelligence and Automation, and the Key Laboratory of Image Processing and Intelligent Control of Education Ministry of China, Huazhong University of Science and Technology, Wuhan 430074, China (email: jingwangli@outlook.com, houshengsu@gmail.com).}
    \thanks{Digital Object Identifier \href{https://doi.org/10.1109/TCNS.2022.3203486}{10.1109/TCNS.2022.3203486}}
}

\maketitle


\begin{abstract}
    A class of distributed optimization problem with a globally coupled equality constraint and local constrained sets is studied in this paper. For its special case where local constrained sets are absent, an augmented primal-dual gradient dynamics is proposed and analyzed, but it cannot be implemented distributedly since the violation of the coupled constraint needs to be used. Benefiting from the brand-new comprehending of a classical distributed unconstrained optimization algorithm, the novel implicit tracking approach is proposed to track the violation distributedly, which leads to the birth of the \underline{i}mplicit tracking-based \underline{d}istribut\underline{e}d \underline{a}ugmented primal-dual gradient dynamics (IDEA). A projected variant of IDEA, i.e., Proj-IDEA, is further designed to deal with the general case where local constrained sets exist. With the aid of the Lyapunov stability theory, the convergences of IDEA and Proj-IDEA over undigraphs and digraphs are analyzed respectively. As far as we know, Proj-IDEA is the first constant step-size distributed algorithm which can solve the studied problem without the need of the strict convexity of local cost functions. Besides, if local cost functions are strongly convex and smooth, IDEA can achieve exponential convergence with a weaker condition about the coupled constraint. Finally, numerical experiments are taken to corroborate our theoretical results.
\end{abstract}


%
\IEEEpeerreviewmaketitle

\section{Introduction}
\label{intro}
Due to its promising application prospects in large-scale machine learning, distributed control, decentralized estimation, smart grid, and many other fields \cite{nedic2018network}, distributed optimization has become one of the most popular topics in the control community. Consider a multi-agents system consisting of $n$ agents, where the communication topology among agents is modeled by a connected graph $\mG$ (may be directed),
the objective is to solve the following constrained optimization problem:
\begin{equation} \label{cenp2} \tag{P1}
    \begin{aligned}
        \min_{x_i \in \mR^{d_i}} \  & \sum_{i=1}^n f_{i}(x_{i})                    \\
        \text{s.t.} \               & \sum_{i=1}^{n}A_ix_i = b,                    \\
                                    & x_{i} \in \mathcal{X}_{i}, \ i=1, \ldots, n,
    \end{aligned}
\end{equation}
where $f_i:\mR^{d_i} \rightarrow \mR$ and $\mathcal{X}_{i} \subseteq \mR^{d_i}$ are the local cost function and the local constrained set of agent $i$ respectively, and $A_i \in \mR^{p \times d_i}$ is the local constraint matrix of agent $i$, which specifies the globally coupled equality constraint together with $b \in \mR^{p}$. It is assumed that \cref{cenp2} has at least a finite optimal solution. Many optimization problems emerging in engineering and management can be formulated as \cref{cenp2}, such as resource allocation, economic dispatch, network utility maximization, and so on. We aim to design distributed algorithms to solve \cref{cenp2}, which means that: first, agent $i$ can only communicate with its neighbors over $\mG$; second, $b$ is known by all agents, but $f_i$, $A_i$, and $\mathcal{X}_{i}$ are all private information of agent $i$ and cannot be exchanged with its neighbors. For the sake of privacy preservation, transmitting $\nabla f_i$ to its neighbors is also forbidden, since it is possible to reconstruct $f_i$ from $\nabla f_i$.

The classical distributed unconstrained optimization problem is
\begin{equation} \label{uc2} \tag{P2}
    \min_{x \in \mR^m} \ \frac{1}{n}\sum_{i=1}^{n}h_i(x),
\end{equation}
where $h_i:\mR^{m} \rightarrow \mR$ is the local cost function. The most significant difference between \cref{cenp2} and \cref{uc2} is that the former has an equality constraint which couples all agent's private information and decision variables, which is exactly the reason that it is called a distributed constraint-coupled optimization problem. Though \cref{uc2} is equivalent to
\begin{equation} \label{ud22} \tag{P3}
    \begin{aligned}
        \min_{x_i \in \mR^m} \  & \frac{1}{n}\sum_{i=1}^{n}h_i(x_i) \\
        \text{s.t.} \           & x_i = x_j, \ i, j = 1, \cdots, n,
    \end{aligned}
\end{equation}
the equality constraint of \cref{ud22} is naturally distributed, which is essentially different from that of \cref{cenp2}. Specifically speaking, to meet the former, each agent only needs to keep its own decision variable the same with its neighbors', under the assumption that $\mG$ is connected. Benefiting from that, it is convenient to apply primal (usually consensus-based) methods \cite{nedic2009distributed,nedic2017achieving,qu2017harnessing,pu2020push} and dual (or primal-dual) methods \cite{wang2010control,gharesifard2013distributed,kia2015distributed,shi2014linear,shi2015extra,yuan2018exact} to derive distributed algorithms for solving \cref{ud22}. However, both the above two approaches fail to work for \cref{cenp2}. The reason is obvious for the former, while for the latter, it is that the resulting algorithms are not distributed since the dual update needs the violation of the coupled constraint.

Nevertheless, it is easy to verify that the dual of \cref{cenp2} has the same form with \cref{uc2}, thus we can employ existing distributed unconstrained optimization algorithms to solve the dual of \cref{cenp2}, then the solution of \cref{cenp2} can be easily obtained from that of its dual if the strong duality holds. In fact, the above idea has been used in many earlier works \cite{chang2014multi,yi2016initialization,zhu2019distributed,alghunaim2019proximal,zhang2020distributed}. \cref{cenp2} and \cref{uc2} are both considered in \cite{chang2014multi}, C-ADMM, a distributed version of Alternating Direction Method of Multipliers (ADMM), is first proposed to solve \cref{uc2}, then DC-ADMM is further derived to solve \cref{cenp2} by applying C-ADMM to its dual. Therefore, DC-ADMM can be seen as a dual variant of C-ADMM.
Similarly, the algorithms proposed in \cite{yi2016initialization}, \cite{zhu2019distributed}, \cite{alghunaim2019proximal}, and \cite{zhang2020distributed} correspond to those proposed in \cite{wang2010control}, \cite{kia2015distributed}, \cite{yuan2018exact}, and \cite{pu2020push} respectively.
However, the convergences of all algorithms mentioned above rely on the strict or strong convexity of the local cost function $f_i$, which limits their application scenarios.
Though the algorithms proposed in \cite{falsone2017dual} and \cite{notarnicola2019constraint} can guarantee convergence without the strict convexity of $f_i$, both of them suffer from a slow convergence rate due to the use of diminishing step-sizes.

Can we design a constant step-size distributed algorithm for \cref{cenp2}, which can converge when $f_i$ is only convex? To meet this goal, it is natural to consider the augmented version of \cref{cenp2}, but the decomposability of the dual problem will be destroyed \cite{boyd2011distributed}. Specifically, the dual of the augmented version of \cref{cenp2} no longer has the same form with \cref{uc2}, which makes it impossible to apply the above approach. An alternative approach is employing ADMM to solve \cref{cenp2} directly, but the primal and dual updates of ADMM both need the violation of the coupled constraint, which prevents it from being implemented distributedly. Nevertheless, it is feasible to use dynamic average consensus (DAC) algorithms to track the violation distributedly \cite{zhu2010discrete,kia2015distributed}, this is how Tracking-ADMM works, which is another distributed version of ADMM \cite{falsone2020tracking}. Beforehand, the idea that tracking the violation has been employed to design distributed versions of existing centralized algorithms in \cite{cherukuri2016initialization} and \cite{kia2017distributed}. As far as we know, though the algorithms proposed in \cite{falsone2020tracking,su2021distributed,cherukuri2016initialization,kia2017distributed,nedic2018improved} can converge without the need of the strict convexity of $f_i$, all of whom can only deal with special cases of \cref{cenp2}: $A_i=1$\cite{cherukuri2016initialization}; $A_i=I$\cite{nedic2018improved}; $A_i\in \mR^{p\times 1}$, $\mathcal{X}_i=\mR$\cite{kia2017distributed}; $\mathcal{X}_i$ is bounded \cite{falsone2020tracking,su2021distributed}. Particularly, in the algorithms proposed in \cite{cherukuri2016initialization} and \cite{nedic2018improved}, each agent needs to transmit the gradient information to its neighbors, which is unfavorable for privacy preservation and not allowed in our scenario.

In this work, we first consider the special case of \cref{cenp2}: $\mathcal{X}_{i} = \mR^{d_i}$, and propose an augmented primal-dual gradient dynamics (APGD) to solve it, but APGD cannot be implemented distributedly since it needs to use the violation of the coupled constraint. Benefiting from the brand-new comprehending of a classical distributed unconstrained optimization algorithm \cite{kia2015distributed}, we propose the novel implicit tracking approach, which is essentially different from the tracking approach adapted in \cite{falsone2020tracking,cherukuri2016initialization,kia2017distributed}, then design a distributed version of APGD, i.e., IDEA. Correspondingly, we call the latter the explicit tracking approach. Besides, a projected variant of IDEA, i.e., Proj-IDEA, is further designed to handle the general case of \cref{cenp2}.

Our major contributions are summarized as follows.
\begin{enumerate}[(i)]
    \item Under undigraphs, the convergence of Proj-IDEA can be guaranteed when $f_i$ is only convex. As far as we know, Proj-IDEA is the first constant step-size distributed algorithm which can solve \cref{cenp2} without the need of the strict or strong convexity of $f_i$. As stated above, existing algorithms \cite{falsone2020tracking,su2021distributed,cherukuri2016initialization,kia2017distributed,nedic2018improved} can only deal with some special cases of \cref{cenp2}. Among these works, the set-up considered in \cite{falsone2020tracking} and \cite{su2021distributed} is the closest one with \cref{cenp2}, where $\mathcal{X}_i$ is bounded. However, the boundedness of $\mathcal{X}_i$ cannot be satisfied by many practical optimization problems, which limits their application scope significantly. Besides, there is a subproblem that needs to be solved exactly at each iteration of Tracking-ADMM, as is well known, obtaining the exact solution of an optimization problem is usually difficult and costly \cite{devolder2014first}. By contrast, $\mathcal{X}_i$ of \cref{cenp2} can be unbounded and no subproblem needs to be solved in Proj-IDEA.
    \item When local constrained sets are absent, there exists several exponentially convergent algorithms for \cref{cenp2}, all of whom require that $f_i$ is strongly convex and smooth, and further impose some different conditions on $A_i$, such as $A_i=1$\cite{kia2017distributed}; $A_i=I$\cite{yi2016initialization,zhu2019distributed,zhang2020distributed,nedic2018improved}; $A_i$ has full row rank\cite{chang2014multi,alghunaim2019proximal}. With the same condition about $f_i$, IDEA can achieve exponential convergence (also called linear convergence) if $A = [A_1, \cdots, A_n]$ has full row rank, which is weaker than existing ones.
    \item Under digraphs, the convergence of Proj-IDEA can be guaranteed when $f_i$ is strongly convex. Though the algorithms proposed in \cite{zhang2020distributed,zhu2019distributed,kia2017distributed} can also converge over digraphs, the problems they considered are all special cases of \cref{cenp2}. In addition, IDEA can converge exponentially over digraphs if we further assume that $A_i = I$. More than that, the exponential convergence of IDEA can be extended to time-varying graphs, no matter for undirected ones or directed ones.
    \item The convergence analysis is based on the Lyapunov approach and LaSalle's invariance principle. Benefiting from the deep understanding of the relation between APGD and IDEA, we are able to design nice Lyapunov functions and depict the largest invariant sets artfully, which is the key in proving the convergences of IDEA and Proj-IDEA.
    \item Thanks to the implicit tracking approach, the number of state variables that need to be exchanged of IDEA is only half that of the \underline{e}xplicit tracking-based \underline{d}istribut\underline{e}d \underline{a}ugmented primal-dual gradient dynamics (EDEA).
          Even though, it is shown that IDEA usually has a faster convergence rate in numerical experiments, which means that the total communication cost of IDEA is usually no more than half that of EDEA for a given accuracy. Besides, the implicit tracking mechanism can offer us a brand-new viewpoint to understand existing distributed unconstrained optimization algorithms. Inspired by the implicit tracking mechanism, we develop a unified algorithmic framework--unified gradient tracking \cite{li2022gradient}, which can unify most existing first-order distributed unconstrained optimization algorithms.
\end{enumerate}

The rest of this paper is organized as follows. In \cref{Pre}, notations and some necessary preliminaries are introduced. In \cref{design}, IDEA and Proj-IDEA are designed. In \cref{convergence}, the convergences of IDEA and Proj-IDEA are analyzed. In \cref{simu}, numerical experiments are taken to corroborate our theoretical results and evaluate the performance of IDEA and Proj-IDEA. Finally, the conclusion is given in \cref{conclusion}.

\section{Preliminaries} \label{Pre}
\subsection{Notations}
The vector of $n$ zeros (ones) and the $n\times n$ identity matrix are denoted by $\0_n$ ($\1_n$) and $I_n$ respectively. Notice that the dimension of a vector or matrix might not be explicitly given, if it could be deduced from the context. For $x \in \mR^{m}$, $\|x\|$ denotes its Euclidean norm. For $A \in \mR^{m \times n}$, $\us(A)$ and $\os(A)$ denote its smallest nonzero and largest singular values respectively, $\|A\|$ denotes its spectral norm, $\mN(A)$ and $\mC(A)$ denote its null and column spaces respectively. For $A \in \mR^{n\times n}$, $A \succ 0$ ($A \succeq 0$) denotes that $A$ is positive definite (semi-definite), $\ue(A)$, $\eta_{2}(A)$, and $\ove(A)$ denote its smallest nonzero, second smallest, and largest eigenvalues respectively. $\text{diag}\{\cdot\}$ denotes a diagonal or block diagonal matrix, which depends on its components specifically. Given y $\in \mR^m$ and $\mathcal{X} \subseteq \mR^m$, the projection operator is defined as $\Px(y) = \arg{\min_{x \in \mathcal{X}}{\|y-x\|}}$.Finally, iff means if and only if.

\subsection{Graph Theory}
A digraph composing of $n$ agents is described by $\mG=(\mathcal{V}, \mathcal{E})$, where $\mathcal{V} = \{1, \cdots, n\}$ and $\mathcal{E} \subset \mathcal{V} \times \mathcal{V}$ are the vertex and edge sets respectively. The edge $(i, j)$ denotes that agent $i$ can receive information sent by agent $j$, naturally, agent $i$ is an out-neighbor of agent $j$ and agent $j$ is an in-neighbor of agent $i$. $\mathcal{A} = [a_{ij}] \in \mR^{n \times n}$ is the weighted adjacency matrix, where $a_{ij} > 0$ if $(i, j) \in \mathcal{E}$ and $a_{ij} = 0$ otherwise. In particular, $(i, i) \in \mathcal{E}$, $\forall i \in \mathcal{V}$. The in-degree and out-degree of agent $i$ are defined as $d_{i}^{in} = \sum_{j=1}^n a_{ji}$ and $d_{i}^{out} = \sum_{j=1}^n a_{ij}$ respectively. $\mG$ is weight-balanced iff $d_{i}^{in} = d_{i}^{out}, \ \forall i \in \mathcal{V}$, and $\mG$ is undirected iff $a_{ij}=a_{ji}, \ \forall i \in \mathcal{V}$, hence undigraphs are naturally weight-balanced.
A directed path is a sequence of edges which connect two agents. $\mG$ is strongly connected if there exists a directed path for any pair of agents of it. The Laplacian of $\mG$ is defined by $L = D^{out} - \mathcal{A}$, where $D^{out} = \text{diag}\{d_1^{out}, \cdots, d_n^{out}\}$, and $0$ is a simple eigenvalue of $L$ if $\mG$ is strongly connected. Define $\hat{L} = \frac{L+L\T }{2}$, then the following statements are equivalent \cite{bullo2009distributed}: (i) $\mG$ is weight-balanced; (ii) $\1_n\T L = \0_n\T$; (iii) $\hat{L} \succeq 0$. Let $m \in \mR^n$ and $M \in \mR^{n\times (n-1)}$ satisfy
\eqe{ \label{mM}
    m = \frac{1}{\sqrt{n}}\1_n, \ m\T M = \0_n\T, \\
    M\T M = I_{n-1}, \ MM\T  = I_n - mm\T,
}
if $\mG$ is strongly connected and weight-balanced, then $0$ is a simple eigenvalue of $\hat{L}$ and the corresponding eigenvector is $a\1_n$, where $a \in \mR$, and it holds that \cite{kia2015distributed}
\begin{equation} \label{sec}
    0 \prec \eta_2(\hat{L}) I \preceq M\T \hat{L}M \preceq \ove(\hat{L}) I.
\end{equation}

\subsection{Convex Analysis}
\begin{lemma} \cite{nesterov2018lectures} \label{strPro}
    Given a convex set $\mathcal{X} \subseteq \mR^{d}$ and a differentiable function $f:\mR^d \rightarrow \mR$, if $f$ is $l$-smooth on $\mathcal{X}$ with $l>0$, then
    $$\|\nabla f(x) - \nabla f(y)\|^2 \leq l(\nabla f(x) - \nabla f(y))\T (x-y), \ \forall x, y \in \mathcal{X}.$$
    If $f$ is convex on $\mathcal{X}$, then
    $$(\nabla f(x) - \nabla f(y))\T (x-y) \geq 0, \ \forall x, y \in \mathcal{X}.$$
    If $f$ is $\mu$-strongly convex on $\mathcal{X}$ with $\mu>0$, then
    $$(\nabla f(x) - \nabla f(y))\T (x-y) \geq \mu\|x-y\|^2, \ \forall x, y \in \mathcal{X}.$$
\end{lemma}

\begin{lemma} \cite{hiriart2012fundamentals} \label{proj}
    Given a closed and convex set $\mathcal{X} \subseteq \mR^d$, for any $x, y \in \mR^d$, it holds that
    \eqe{
        (\Px(y)-y)\T (\Px(y) - \Px(x)) \leq 0, \\
        (x-y)\T (\Px(x)-\Px(y)) \geq \|\Px(x)-\Px(y)\|^2.
        \nonumber
    }
\end{lemma}

\begin{lemma} \cite{liu2013one} \label{prd}
    Given a closed and convex set $\mathcal{X} \subseteq \mR^d$, define $\psi:\mR^d \rightarrow \mR$ as
    $$\psi(x) = \|x-x^*\|^2 - \|x-\Px(x)\|^2,$$
    where $x^* \in \mathcal{X}$. For any $x \in \mR^d$, we have
    \begin{enumerate}[(i)]
        \item $\psi(x) \geq \|\Px(x)-x^*\|^2$;
        \item $\psi$ is differentiable and $\nabla \psi(x) = 2(\Px(x)-x^*)$.
    \end{enumerate}
\end{lemma}

\section{Algorithm Design} \label{design}
\subsection{Centralized Augmented Primal-dual Gradient Dynamics}
\cref{cenp2} can be reformulated as
\begin{equation} \label{cenp} \tag{P4}
    \begin{aligned}
        \min_{\mx \in \mR^{d}} \  & f(\mx) = \sum_{i=1}^n f_{i}(x_{i}) \\
        \text{s.t.} \             & A\mx = b,                          \\
                                  & \mx \in \mathcal{X},
    \end{aligned}
\end{equation}
where $\mx = [x_1\T, \cdots, x_n\T]\T\in \mR^{d}$, $d = \sum_{i=1}^nd_i$, $A = [A_1, \cdots, A_n] \in \mR^{p \times d}$, and $\mathcal{X}=\prod_{i=1}^{N} \mathcal{X}_{i} \in \mR^{d}$ is the Cartesian product of $\mathcal{X}_{1}, \cdots$, $\mathcal{X}_{n-1}$ and $\mathcal{X}_{n}$. The following assumption holds throughout this paper.

\begin{assumption} \label{dif}
    $\mathcal{X}_i$ is closed and convex, $f_{i}$ is convex on $\mathcal{X}_i$ and differentiable on an open set containing $\mathcal{X}_i$, and $\nabla f_i$ is locally Lipschitz on an open set containing $\mathcal{X}_i$, $\forall i \in \mathcal{V}$. The Slater’s condition is satisfied by \cref{cenp}, i.e., there exists a relative interior point of $\mathcal{X}$ such the equality constraint holds.
\end{assumption}

Define the augmented Lagrangian of \cref{cenp} as
\begin{equation} \label{Lag}
    \begin{aligned}
        \cLa(\mx, \lambda) = \underbrace{f(\mx) + \lambda\T (A\mx - b)}_{\cL(\mx, \lambda)} + \frac{1}{2\alpha}\|A\mx-b\|^2,
    \end{aligned}
    \nonumber
\end{equation}
where $\lambda \in \mR^{p}$ is the Lagrange multiplier, $\alpha > 0$ is a constant, and $\cL(\mx, \lambda)$ is the standard (or unaugmented) Lagrangian. Since \cref{cenp} is convex and Slater's condition holds, the strong duality holds, which implies that the existence of the saddle point of $\cLa$ can be guaranteed. Notice that if $(\mx^*, \lambda^*)$ is a saddle point of $\cLa$, then $\mx^*$ and $\lambda^*$ are the optimal solutions of \cref{cenp} and its dual respectively \cite{bertsekas1997nonlinear}. Therefore, we can obtain the optimal solution of \cref{cenp} by solving its corresponding saddle point problem:
\eqe{ \label{sadP}
    \min_{\mx \in \mathcal{X}}\max_{\lambda\in \mathbb{R}^p}\cLa(\mx, \lambda).
    \nonumber
}

As is well known, primal-dual gradient dynamics is a popular and efficient approach to solve saddle point problems. The seminal works of primal-dual gradient dynamics can date back to 1950s \cite{kose1956solutions}, and recently, it has drawn many research interests \cite{feijer2010stability,wang2010control,gharesifard2013distributed,cherukuri2017saddle,qu2018exponential,cortes2019distributed}.
To solve the special case of \cref{cenp}: $\mathcal{X} = \mR^{d}$, we propose the following primal-dual gradient dynamics:
\eqe{ \label{sadDy}
    \dot{\mx} &= -\alpha\nabla_\mx \cLa(\mx, \lambda), \\
    \dot{\lambda} &= \nabla_\mathbf{\lambda} \cLa(\mx, \lambda),
    \nonumber
}
which can be unfolded as
\begin{equation} \label{cena}
    \begin{aligned}
        \dot{\mx}     & = -\alpha(\nabla f(\mx) + A\T \lambda) - A\T (A\mx-b), \\
        \dot{\lambda} & = A\mx - b,
    \end{aligned}
    \nonumber
\end{equation}
this is APGD mentioned before.

If $(\mx^*, \lambda^*)$ is an equilibrium point of APGD, it holds that
\eqe{ \label{kkt}
    \0 &= \nabla f(\mx^*) + A\T \lambda^*, \\
    \0 &= A\mx^* - b,
}
obviously $(\mx^*, \lambda^*)$ is a saddle point of $\cLa$. Let $(\mx(t), \lambda(t))$ be the trajectory of APGD with the initial point $(\mx(0), \lambda(0)) \in \mR^d \times \mR^{p}$, then we have the following two propositions.

\begin{proposition} \label{cenCon}
    Suppose \cref{dif} holds, $(\mx(t), \lambda(t))$ converges to a saddle point of $\cLa$.
\end{proposition}

\begin{proposition} \label{cenStr}
    Suppose \cref{dif} holds, $\lambda(0)=\0$, and $f_i$ is $\mu_i$-strongly convex, $\forall i \in \mathcal{V}$, $(\mx(t), \lambda(t))$ converges exponentially to a saddle point of $\cLa$.
\end{proposition}

Though $A$ of \cref{cenp} has a specific structure, \cref{cenCon,cenStr} are applicable for any $A$ since none of its structure information is used in the proof.

\begin{remark}
    The key role in APGD that allows it to converge without the strict convexity of $f_i$, is $A\T (A\mx-b)$, which corresponds to the augmented term of $\cLa$. From the view of control, $A\T (A\mx-b)= A\T \dot{\lambda}$ can be seen as a ``derivative feedback'' technique. Without it, APGD will no longer possess the above property. Nevertheless, it is not indispensable for the exponential convergence. In fact, APGD can still achieve exponential convergence without it, which can be easily proved based on the proof of \cref{cenStr}.
\end{remark}

\begin{remark}
    It is worth mentioning that \cref{cenCon} can also be proved by applying Corollary 4.5 of \cite{cherukuri2017saddle}, where the asymptotic stability of general primal-dual gradient dynamics is analyzed. Nevertheless, we offer a different proof idea, which is quite beneficial to the convergence analysis of IDEA. Beforehand, the exponential convergence of primal-dual gradient dynamics has been studied in \cite{qu2018exponential} and \cite{cortes2019distributed}. However, \cite{qu2018exponential} only considers the standard Lagrangian, which can be seen as a special case of $\cLa$: $\alpha \rightarrow \infty$. In \cite{cortes2019distributed}, $\alpha$ of $\cLa$ is limited in $(0, 1)$, which is a necessary condition for the exponential convergence. Besides, both \cite{qu2018exponential} and \cite{cortes2019distributed} require that $A$ has full row rank, but we do not impose any condition on $A$.
\end{remark}

\begin{remark} \label{rin}
    APGD is essentially a kind of inexact augmented Lagrangian method \cite{nedelcu2014computational}, which is given as
    \eqe{
        \mx\+ &\approx \arg \min_{\mx \in \mR^{d}}\cLa(\mx, \lambda^k), \\
        \lambda\+ &= \lambda^k + \frac{1}{\alpha}(A\mx\+-b).
        \nonumber
    }
    To see that, consider the discretized form of APGD:
    \eqe{ \label{dce}
        \mx\+ &= \mx^k -\nabla f(\mx^k) - A\T \lambda - \frac{1}{\alpha}A\T (A\mx^k-b), \\
        \lambda\+ &= \lambda^k + \frac{1}{\alpha}(A\mx\+-b),
        \nonumber
    }
    where we take $\frac{1}{\alpha}$ as the discretized step-size and use $\mx\+$ instead of $\mx^k$ in the dual update. $\mx\+$ can be seen as a rough approximation of the optimal solution of $\min_{\mx \in \mR^{d}}\cLa(\mx, \lambda^k)$, where only one step gradient descent of $\cLa(\mx, \lambda^k)$ is taken.
    Furthermore, note that $\cLa$ can be regarded as the standard Lagrangian of
    \eqe{ \label{ecenp}
        \min_{\mx \in \mR^{d}} \ &f(\mx) + \frac{1}{2\alpha}\|A\mx-b\|^2 \\
        \text{s.t.} \ &A\mx = b,
    }
    which is equivalent to \cref{cenp}. The dual function of \cref{ecenp} is
    \eqe{
        g(\lambda) = \inf_{\mx \in \mR^{d}}\cLa(\mx, \lambda),
        \nonumber
    }
    define $\mx^*(\lambda) = \arg \min_{\mx \in \mR^{d}}\cLa(\mx, \lambda)$, $\mx^*(\lambda)$ is unique \cite{bertsekas1997nonlinear} if $\alpha$ is small enough, then $g$ is differentiable and
    \eqe{
        \nabla g(\lambda) = A\mx^*(\lambda) - b.
        \nonumber
    }
    Therefore, $A\mx-b$ of APGD can be seen as the inexact dual gradient.
\end{remark}

\subsection{Implicit Tracking Approach}
Though APGD can solve \cref{cenp}, it cannot be implemented in a distributed manner since its primal and dual updates both need the global information $A\mx-b= \sum_{i=1}^{n}A_ix_i-b$. One solution that obtain distributed versions of APGD is employing DAC algorithms to track $A\mx-b$ distributedly, an example is EDEA:
\begin{subequations} \label{adac2}
    \begin{align}
        \dot{\mx} & = -\alpha(\nabla f(\mx) + \mA\T \ml) - \mA\T \mr, \label{adac2a} \\
        \dot{\ml} & = \mr-\mL\ml, \label{adac2b}                                     \\
        \dot{\mr} & = - \gamma(\mr-(\mA\mx-\mb)) -\mz - \beta\mL\mr, \label{adac2c}  \\
        \dot{\mz} & = \gamma\beta\mL\mr, \label{adac2d}
    \end{align}
\end{subequations}
where $\mL = L \otimes I_p$, $\mA = \text{diag}\{A_1, \cdots, A_n\} \in \mR^{np \times d}$, $\mb = [b_1\T, \cdots, b_n\T]\T\in \mR^{np}$ (where $b_i$ satisfies $\sum_{i=1}^{n}b_i=b$), $\ml = [\lambda_1\T, \cdots, \lambda_n\T]\T\in \mR^{np}$, $\mr = [r_1\T, \cdots, r_n\T]\T\in \mR^{np}$, and $\mz = [z_1\T, \cdots, z_n\T]\T\in \mR^{np}$. The key role in EDEA is \cref{adac2c,adac2d}, which is the continuous-time DAC algorithm proposed in \cite{kia2015dynamic}, theoretically $r_i$ will track to $\sum_{i=1}^{n}A_ix_i-b$ if $\gamma$ and $\beta$ are chosen appropriately, i.e.,
$$\lim_{t \rightarrow \infty}r_i(t) = \frac{1}{n}\sum_{i=1}^{n}A_ix_i(t)-b, \ i \in \mathcal{V}.$$
Therefore, EDEA will act like APGD as time goes on. However, there are two state variables that need to be exchanged in EDEA. Is is possible design another distributed version of APGD, with fewer state variables that need to be exchanged? To achieve that goal, we can resort to the novel implicit tracking approach.

To illustrate the implicit tracking approach, consider the unconstrained set-up \cref{uc2} and assume that it has at least an optimal solution $x^*$. As mentioned before, \cref{uc2} is equivalent to \cref{ud22}.
If each agent takes the following dynamics:
\eqe{ \label{gf}
    \dot{x}_i(t) = -\frac{1}{n}\sum_{j=1}^{n} \nabla h_j(x_j(t)) - \lt(x_i(t)-\frac{1}{n}\sum_{j=1}^{n}x_j(t)\rt),
}
it is easy to verify that $x_i(t)$ will eventually converge to $x^*$, but \cref{gf} is obviously not distributed. Similarly, we can employ the DAC algorithm to track the global gradient $\frac{1}{n}\sum_{j=1}^{n} \nabla h_j(x_j(t))$, which leads to the following distributed algorithm:
\eqe{ \label{adac}
    \dot{\mx} &= -\mr- \mL\mx, \\
    \dot{\mr} &= - \gamma(\mr-\nabla h(\mx)) -\mz - \beta\mL\mr, \\
    \dot{\mz} &= \gamma\beta\mL\mr,
}
where $\mx = [x_1\T, \cdots, x_n\T]\T\in \mR^{mn}$, $\mr = [r_1\T, \cdots, r_n\T]\T\in \mR^{mn}$, and $\mz = [z_1\T, \cdots, z_n\T]\T\in \mR^{mn}$. It is worth mentioning that \cref{adac} is not a new algorithm essentially, since its discrete-time counterpart has been studied in \cite{nedic2017achieving} and \cite{qu2017harnessing}, where the discrete-time DAC algorithm \cite{zhu2010discrete} is used.

Except \cref{adac}, there exists another distributed version of \cref{gf}. Notice that \cref{gf} can be rewritten as
\eqe{ \label{gf2}
\dot{x}_i(t) = \frac{1}{n}\sum_{j=1}^{n}(x_j(t)-\nabla h_j(x_j(t))) - x_i(t),
\nonumber
}
it is reasonable to track $\frac{1}{n}\sum_{j=1}^{n}(x_j(t)-\nabla h_j(x_j(t)))$ directly, which leads to the following distributed algorithm:
\eqe{ \label{ap4}
    \dot{\mx} &= - \gamma(\mx-(\mx-\nabla h(\mx))) -\mz - \beta\mL\mx, \\
    \dot{\mz} &= \gamma\beta\mL\mx,
    \nonumber
}
i.e.,
\eqe{ \label{ap}
    \dot{\mx} &= -\gamma\nabla h(\mx) -\mz - \beta\mL\mx, \\
    \dot{\mz} &= \gamma\beta\mL\mx,
}
which is originally proposed in \cite{kia2015distributed}, but with a different design idea. The dynamics of agent $i$ is given as
\eqe{ \label{api}
\dot{x}_i(t) &= - \gamma\nabla h_i(x_i(t)) -z_i(t) - \beta\sum_{j=1}^na_{ij}(x_i(t)-x_j(t)), \\
\dot{z}_i(t) &= \gamma\beta\sum_{j=1}^na_{ij}(x_i(t)-x_j(t)).
}
Along \cref{api}, theoretically we have
$$\lim_{t \rightarrow \infty}x_i(t) = \frac{1}{n}\sum_{j=1}^{n}(x_j(t)-\nabla h_j(x_j(t))),$$
if $\gamma$ and $\beta$ are chosen appropriately. Therefore, \cref{api} will act like \cref{gf} as time goes on, which implies that
$$\lim_{t \rightarrow \infty}\gamma\nabla h_i(x_i(t))+z_i(t) = \frac{1}{n}\sum_{j=1}^{n} \nabla h_j(x_j(t)).$$
That is, $\gamma\nabla h_i(x_i(t))+z_i(t)$ plays a role in tracking the global gradient $\frac{1}{n}\sum_{j=1}^{n} \nabla h_j(x_j(t))$ and $z_i(t)$ can be seen as agent $i$'s estimation of the gap between its local gradient and the global gradient. Different from \cref{adac}, there is no explicit state variable to track the global gradient in \cref{ap}, instead, the global gradient is tracked in an implicit way. Therefore, we call the tracking approach used in \cref{ap} the implicit tracking approach (mechanism). In the following, we will show how to use the implicit tracking approach to design a more efficient distributed version of APGD.

\begin{remark}
    The implicit tracking mechanism dose not simply exist in \cref{ap}, but many existing distributed unconstrained optimization algorithms. Inspired by the implicit tracking mechanism, we develop a unified algorithmic framework--unified gradient tracking \cite{li2022gradient}, which can unify most existing first-order distributed unconstrained optimization algorithms.
\end{remark}

\begin{remark} \label{rem4}
    Though \cref{ap} is not a new algorithm, the above derivation of \cref{ap} and the resulting implicit tracking approach are novel. In \cite{kia2015distributed}, the design of \cref{ap} is based on the feedback control idea. Specifically speaking, it is observed that
    \eqe{
        \dot{\mx} = -\gamma\nabla h(\mx) - \beta\mL\mx
        \nonumber
    }
    cannot converge to the optimal solution since local gradients are generally different. Therefore, the integral feedback term $\mz$ is constructed to correct the error among agents caused by local gradients. In fact, we can offer the third explanation of \cref{ap} when $\mG$ is undirected and connected (which implies that $L \succeq 0$). Let $L = \sqrt{L}\sqrt{L}$, where $\sqrt{L}=PD^{\frac{1}{2}}P\T$, $P$ is an orthogonal matrix, and $D$ is the diagonal matrix consisting of the eigenvalues of $L$. Like $L$, $\sqrt{L}$ is also symmetric. Besides, $L$ and $\sqrt{L}$ possess the following two properties: (i) $\1\T L = \0\T$ is equivalent to $\1\T \sqrt{L} = \0\T$; (ii) for any $z = [z_1, \cdots, z_n]\T\in \mR^n$, $\sqrt{L} z = \0$ is equivalent to $Lz = \0$, which implies that $z_1 = \cdots = z_n$.
    Define $\sL = \sqrt{L} \otimes I_m$, then \cref{uc2} is equivalent to
    \eqe{ \label{ud}
        \min_{\mx \in \mR^{mn}} \ &h(\mx) \\
        \text{s.t.} \ &\sL\mx = \0,
    }
    whose augmented Lagrangian can be defined as
    \eqe{ \label{Lc}
        \cL_{c}(\mx, \my) = h(\mx) + \my\T \sL\mx + \frac{\beta}{2\gamma}\mx\mL\mx.
        \nonumber
    }
    Then we can employ APGD (let $A = \sL$ and $b = \0$):
    \begin{subequations} \label{ap3}
        \begin{align}
            \dot{\mx} & = - \gamma(\nabla h(\mx) +\sL\my) - \beta\mL\mx, \label{ua} \\
            \dot{\my} & = \beta\sL\mx, \label{ub}
        \end{align}
    \end{subequations}
    to solve \cref{ud}, letting $\mz = \sL\my$ gives that
    \eqe{
        \dot{\mx} &= - \gamma(\nabla h(\mx) +\mz) - \beta\mL\mx, \\
        \dot{\mz} &= \beta\mL\mx,
        \nonumber
    }
    which is equivalent to \cref{ap}. Therefore, \cref{ap} can be seen as a special case of APGD, when $\mG$ is undirected and connected. It is worth noting that this perspective can help us study the convergence of \cref{ap}. Concretely speaking, due to the equivalence of \cref{ap} and \cref{ap3}, \cref{cenCon,cenStr} are also applicable for \cref{ap}. In \cite{kia2015distributed}, it is proved that \cref{ap} can converge if $h_i$ is convex and $h$ is strictly convex, but according to \cref{cenCon}, we know that the strict convexity of $h$ is not necessary.
\end{remark}

\subsection{Implicit Tracking-based Distributed Augmented Primal-dual Gradient Dynamics}
The dual function of \cref{cenp} is given as
\eqe{
    G(\lambda) &= \inf_{\mx \in \mathcal{X}}\cL(\mx, \lambda) \\
    &= \inf_{\mx \in \mathcal{X}}f(\mx) + \lambda\T (A\mx - b) \\
    &= \sum_{i=1}^n\inf_{x_i \in \mathcal{X}_i}f_i(x_i) + \lambda\T (A_ix_i - b_i) \\
    &= \sum_{i=1}^ng_i(\lambda),
    \nonumber
}
where $b_i$ satisfies $\sum_{i=1}^nb_i = b$, then the dual problem of of \cref{cenp} can be written as
\eqe{
    \max_{\lambda \in \mR^{p}} \frac{1}{n}\sum_{i=1}^ng_i(\lambda),
    \nonumber
}
which is equivalent to
\eqe{ \label{dD}
    \max_{\ml \in \mR^{np}} \ &g(\ml) = \frac{1}{n}\sum_{i=1}^ng_i(\lambda_i), \\
    \text{s.t.} \ &\lambda_i = \lambda_j, \ i,j \in \mathcal{V},
}
where $\ml = [\lambda_1\T, \cdots, \lambda_n\T]\T$.
If $f_i$ is strictly convex, $x_i^*(\lambda_i)$ is unique, then $g_i$ is differentiable and
\eqe{
    \nabla g_i(\lambda_i) = A_ix_i^*(\lambda_i)-b_i,
    \nonumber
}
where
$$x_i^*(\lambda_i) = \arg \min_{x_i \in \mathcal{X}_i}\Big\{f_i(x_i) + \lambda_i\T (A_ix_i - b_i)\Big\}.$$
Employing \cref{ap} to solve \cref{dD} gives that
\eqe{ \label{apd}
    \dot{\ml} &= \gamma(\mA\mx^*(\ml)-\mb) -\mz - \beta\mL\ml, \\
    \dot{\mz} &= \gamma\beta\mL\ml,
}
where $\mx^*(\ml) = [x_1^*(\lambda_1)\T, \cdots, x_n^*(\lambda_n)\T]\T$. Ignoring local constrained sets for the moment and replacing the exact dual gradient with the inexact one, we can further obtain the following inexact counterpart of \cref{apd}:
\eqe{ \label{pid2}
    \dot{\mx} &= -\alpha (\nabla f(\mx) + \mA\T \ml), \\
    \dot{\ml} &= \gamma(\mA\mx-\mb)-\mz - \beta \mL\ml, \\
    \dot{\mz} &= \gamma\beta \mL\ml.
    \nonumber
}
Recall the implicit tracking mechanism discussed before, then $\gamma(\mA\mx-\mb)-\mz$ can serve as a distributed approximation of the inexact global dual gradient $A\mx-b$, which leads to another distributed version of APGD:
\eqe{ \label{pid}
    \dot{\mx} &= -\alpha (\nabla f(\mx) + \mA\T \ml) -\mA\T (\gamma(\mA\mx-\mb)-\mz), \\
    \dot{\ml} &= \gamma(\mA\mx-\mb)-\mz - \beta \mL\ml, \\
    \dot{\mz} &= \gamma\beta \mL\ml.
    \nonumber
}
Based on the ideology that using as few parameters as possible on the premise of guaranteeing the convergence, we further derive IDEA:
\eqe{ \label{sys0}
    \dot{\mx} &= -\alpha (\nabla f(\mx) + \mA\T \ml) -\mA\T (\mA\mx-\mb-\mz), \\
    \dot{\ml} &= \mA\mx-\mb-\mz - \beta \mL\ml, \\
    \dot{\mz} &= \alpha\beta \mL\ml.
    \nonumber
}

To deal with the general case of \cref{cenp}, we further develop a projected variant of IDEA, which is called Proj-IDEA:
\eqe{ \label{pidea}
    \dot{\mw} &= -\alpha (\mw-\mx+\nabla f(\mx) + \mA\T \ml) -\mA\T (\mA\mx-\mb-\mz), \\
    \dot{\ml} &= \mA\mx-\mb-\mz - \beta \mL\ml, \\
    \dot{\mz} &= \alpha\beta \mL\ml, \\
    \mx &= \Px(\mw).
    \nonumber
}
For agent $i$, Proj-IDEA is implemented as
\eqe{
\dot{w}_i(t) =& -\alpha \lt(w_i(t)-x_i(t)+\nabla f_i(x_i(t)) + A_i\T \lambda_i(t)\rt) \\
&-A_i\T m_i(t), \\
\dot{\lambda}_i(t) =& m_i(t) - \beta \sum_{j=1}^na_{ij}(\lambda_i(t)-\lambda_j(t)), \\
\dot{z}_i(t) =& \alpha\beta \sum_{j=1}^na_{ij}(\lambda_i(t)-\lambda_j(t)), \\
x_i(t) =& \Pxi(w_i(t)), \\
m_i(t) =& A_ix_i(t)-b_i - z_i(t).
\nonumber
}

\begin{remark}
    As two distributed versions of APGD, $\mA\T (\mA\mx-\mb-\mz)$ of IDEA and $\mA\T \mr$ of EDEA correspond to $A\T (A\mx-b)$ of APGD, which are the keys that allow IDEA and EDEA to converge when $f_i$ is not strictly convex. Though IDEA has only one state variable that needs to be exchanged, it usually shows a better convergence rate in numerical experiments, which means that the total communication cost of IDEA is usually no more than half that of EDEA for a given accuracy, this is exactly the reason we choose IDEA instead of EDEA. Without $\mA\T (\mA\mx-\mb-\mz)$, IDEA will become
    \eqe{ \label{idea2}
        \dot{\mx} &= -\alpha (\nabla f(\mx) + \mA\T \ml), \\
        \dot{\ml} &= \mA\mx-\mb-\mz - \beta \mL\ml, \\
        \dot{\mz} &= \alpha\beta \mL\ml.
        \nonumber
    }
    we call it Unaugmented IDEA. It is worth mentioning that a special case of Unaugmented IDEA has been studied in \cite{zhu2019distributed} ($A_i=I$, $\alpha = 1$), and before that, a similar algorithm has been proposed in \cite{yi2016initialization}:
    \eqe{ \label{yi}
        \dot{\mx} &= -\nabla f(\mx) - \ml, \\
        \dot{\ml} &= \mx-\mb-\mL\mz - \mL\ml, \\
        \dot{\mz} &= \mL\ml.
    }
    The major limitation of Unaugmented IDEA and \cref{yi} is that they cannot converge for non-strictly convex cases. More than that, for strongly convex cases, the convergence rate of IDEA is usually slightly faster than Unaugmented IDEA in numerical experiments.
\end{remark}

\section{Convergence Analysis} \label{convergence}
In this section, the convergences of IDEA and Proj-IDEA over undigraphs and digraphs are analyzed respectively.
For the sake of readability, all proofs are placed in the appendix.
\subsection{Convergence of IDEA}
We first analyze the convergence of IDEA over undigraphs, then extend the results to digraphs.
To begin the convergence analysis, it is necessary to figure out the relation between the equilibrium point of IDEA and the saddle point of $\cLa$, which is depicted by the following lemma.
\begin{lemma} \label{optCon}
    Suppose \cref{dif} holds, $\mG$ is a strongly connected and weight-balanced digraph, and $\sum_{i=1}^nz_i(0) =\0$, then $(\mx^*, \lambda^*)$ is a saddle point of $\cLa$ iff there exists $\mz^* \in \mR^{np}$ such that $(\mx^*, \ml^*, \mz^*)$ is an equilibrium point of IDEA, where $\ml^* = \1_n \otimes \lambda^*$.
\end{lemma}

\begin{assumption} \label{ug}
    $\mG$ is undirected and connected.
\end{assumption}

As mentioned before, the major advantage of IDEA is that its convergence can be guaranteed when $f_i$ is only convex, which benefits from the distributed augmented term $\mA\T (\mA\mx-\mb-\mz)$.
As a distributed version of APGD, this property of IDEA is in line with \cref{cenCon}, which is presented in the following theorem.

\begin{theorem} \label{the1}
    Suppose \cref{dif,ug} hold and $\mathcal{X}_{i} = \mR^{d_i}$, $\forall i \in \mathcal{V}$. Then, given $\alpha, \ \beta > 0$, for any $(\mx(0), \ml(0), \mz(0)) \in \mR^d \times \mR^{np} \times \mR^{np}$ satisfies $\sum_{i=1}^nz_i(0) =\0$, $(\mx(t), \ml(t),\mz(t))$ driven by IDEA converges to $(\mx^*, \ml^*, \mz^*)$, where $\mx^*$ is an optimal solution of \cref{cenp}.
\end{theorem}

The following theorem shows that IDEA can achieve exponential convergence when $f_i$ is strongly convex and smooth and $A$ has full row rank.
\begin{theorem} \label{the2}
    Suppose \cref{dif,ug} hold, $A$ has full row rank, $f_i$ is $\mu_i$-strongly convex and $l_i$-smooth, $\mathcal{X}_{i} = \mR^{d_i}$, $\forall i \in \mathcal{V}$. Define $\mu = \min_{i \in \mathcal{V}}{\mu_i}$, $l = \max_{i \in \mathcal{V}}{l_i}$, and $\varphi>0$, let $\beta > 0$, $\varphi$ and $\alpha$ satisfy
    \eqe{ \label{th2step}
        \varphi &> \max\lt\{\frac{l}{2}-1, 2\os^2(\mA)\rt\}, \\
        \alpha &> \max\lt\{1, \frac{\frac{\varphi\os^2(\mA)}{\mu}+\frac{l}{2}}{\varphi+1-\frac{l}{2}}\rt\}.
    }
    Then, for any $(\mx(0), \ml(0), \mz(0)) \in \mR^d \times \mR^{np} \times \mR^{np}$ satisfies $\sum_{i=1}^nz_i(0) =\0$, $(\mx(t), \ml(t),\mz(t))$ driven by IDEA converges exponentially to $(\mx^*, \ml^*, \mz^*)$, where $\mx^*$ is the unique optimal solution of \cref{cenp}.
\end{theorem}

\begin{remark}
    An important property used in the proof of \cref{cenStr} is that $z\T AA\T z \geq \us^2(A)\|z\|^2$ for $z \in \mC(A)$. For APGD, if we let $\lambda(0)=\0$, then $\lambda(t)$ will always stay in $\mC(A)$. As a consequence, the exponential convergence of APGD can be guaranteed even if $A$ does not have full row rank. However, for IDEA, we cannot make $\ml(t)$ always stay in $\mC(\mA)$ or $\mC(\mL)$, hence the above property cannot be applied, then some extra conditions have to be imposed to $\mA$ (or $A$). As shown in \cref{the2}, the extra condition is that $A$ has full row rank, which is weaker than existing results.
\end{remark}

\begin{remark}
    Recently we noticed the independent work \cite{alghunaim2021dual}, where the dual consensus proximal algorithm (DCPA) is proposed, which can achieve exponential convergence under the same condition with IDEA. However, there are several significant differences between IDEA and DCPA: (i) DCPA is essentially a distributed version of the unaugmented primal-dual gradient algorithm, hence its convergence relies on the strong convexity of $f_i$, but IDEA, as a distributed augmented primal-dual gradient dynamics, does not suffer from that issue; (ii) As a continuous-time algorithm, the convergence analysis of IDEA is based on the Lyapunov stability theory, which is of independent interest.
\end{remark}

\begin{remark}
    Note that $\os(\mA) = \max_{i \in \mathcal{V}}\os(A_i)$, it is not difficult to estimate the lower bound of $\mu$ and the upper bounds of $l$ and $\os(\mA)$ distributedly if $\mu_i$, $l_i$, and $\os(A_i)$ are known by agent $i$, and $\sum_{i=1}^nz_i(0) =\0$ can be satisfied by simply setting $z_i(0) =\0$. Therefore, IDEA can be initialized distributedly.
\end{remark}

Since the Lyapunov function defined in the proof of \cref{the2} does not contain any information about $\mG$, it is straightforward to extend \cref{the2} to time-varying graphs, as the following corollary shows.
\begin{corollary} \label{cor1}
    Suppose \cref{dif} holds, the time-varying $\mG$ is always connected and its adjacency matrix is piecewise constant and uniformly bounded, and other conditions are the same with \cref{the2}.
    Then, for any $(\mx(0), \ml(0), \mz(0)) \in \mR^d \times \mR^{np} \times \mR^{np}$ satisfies $\sum_{i=1}^nz_i(0) =\0$, $(\mx(t), \ml(t),\mz(t))$ driven by IDEA converges exponentially to $(\mx^*, \ml^*, \mz^*)$, where $\mx^*$ is the unique optimal solution of \cref{cenp}.
\end{corollary}

The proof of \cref{cor1} is almost the same with \cref{the2}, except that $\eta_2(L)$ needs to be replaced by $\min_{s \in \mathcal{S}}\eta_2(L_{s})$, where $\mathcal{S}$ is the index set of all possible structures of $\mG$, hence we omit it.

\begin{assumption} \label{dg}
    $\mG$ is directed, strongly connected and weight-balanced.
\end{assumption}

For undigraphs, IDEA can converge even if $f_i$ is only convex, but for digraphs, we need further to assume that $f_i$ is strongly convex, as shown in the following theorem.

\begin{theorem} \label{the3}
    Suppose \cref{dif,dg} hold, $f_i$ is $\mu_i$-strongly convex, $\mathcal{X}_{i} = \mR^{d_i}$, $\forall i \in \mathcal{V}$. Define $\mu = \min_{i \in \mathcal{V}}{\mu_i}$ and $\varphi > 0$, let $\alpha$ and $\beta$ satisfy
    \eqe{ \label{115}
        \alpha &\geq \frac{(\varphi^2+3\varphi+3)\os^2(\mA)}{(\varphi+1)\mu}, \\
        \beta &\geq \frac{(\varphi+1)^2\alpha}{\varphi\eta_2(\hat{L})}.
    }
    Then, for any $(\mx(0), \ml(0), \mz(0)) \in \mR^d \times \mR^{np} \times \mR^{np}$ satisfies $\sum_{i=1}^nz_i(0) =\0$, $(\mx(t), \ml(t),\mz(t))$ driven by IDEA converges to $(\mx^*, \ml^*, \mz^*)$, where $\mx^*$ is the unique optimal solution of \cref{cenp}.
\end{theorem}

Though $\eta_2(\hat{L})$ is global information, it can be calculated distributedly \cite{charalambous2015distributed}. Therefore, IDEA can still be initialized distributedly for digraphs.

The following theorem explores the exponential convergence of IDEA over digraphs, where we assume that $A_i = I_p$.
\begin{theorem} \label{the4}
    Suppose \cref{dif,dg} hold, $A_i = I_p$, $f_i$ is $\mu_i$-strongly convex and $l_i$-smooth, $\mathcal{X}_{i} = \mR^{d_i}$, $\forall i \in \mathcal{V}$. Define  $\mu = \min_{i \in \mathcal{V}}{\mu_i}$, $l = \max_{i \in \mathcal{V}}{l_i}$, and $\varphi > 0$, let $\alpha$ and $\beta$ satisfy
    \eqe{ \label{116}
        \alpha &\geq \max\bigg\{\frac{1}{2}, \frac{(\varphi^2+3\varphi+3)+l^2+\frac{3}{2} - \mu}{(\varphi+1)\mu}\bigg\}, \\
        \beta &\geq \frac{2(\varphi+1)^2\alpha^2+1}{2\varphi\alpha\eta_2(\hat{L})}.
    }
    Then, for any $(\mx(0), \ml(0), \mz(0)) \in \mR^d \times \mR^{np} \times \mR^{np}$ satisfies $\sum_{i=1}^nz_i(0) =\0$, $(\mx(t), \ml(t),\mz(t))$ driven by IDEA converges exponentially to $(\mx^*, \ml^*, \mz^*)$, where $\mx^*$ is the unique optimal solution of \cref{cenp}.
\end{theorem}

The Lyapunov function defined in the proof of \cref{the4} does not contain any information about $\mG$, hence we can also extend \cref{the4} to time-varying graphs, which leads to the following corollary.

\begin{figure*}[tb]
    \begin{center}
        \subfigure[]{
            \begin{minipage}[b]{0.48\textwidth}
                \includegraphics[scale=0.2]{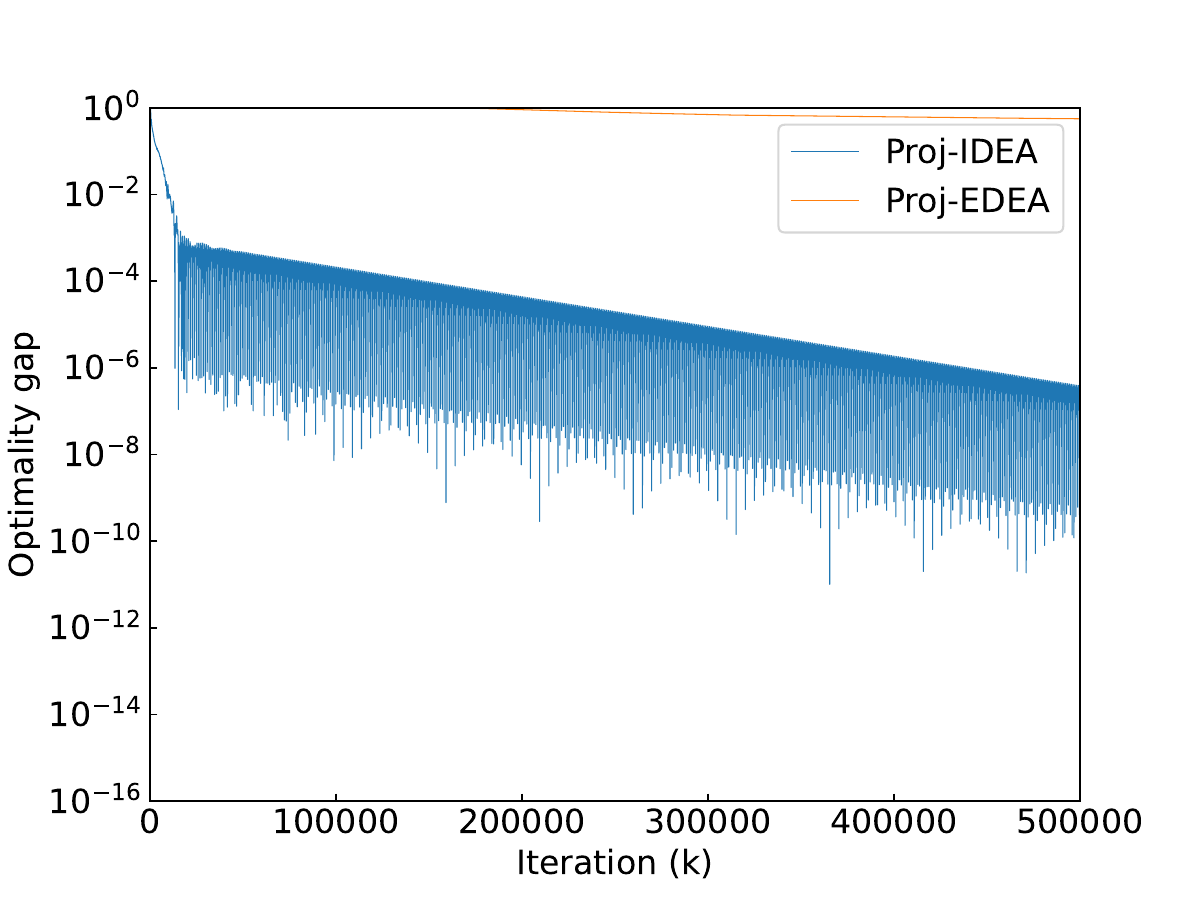}
                \includegraphics[scale=0.2]{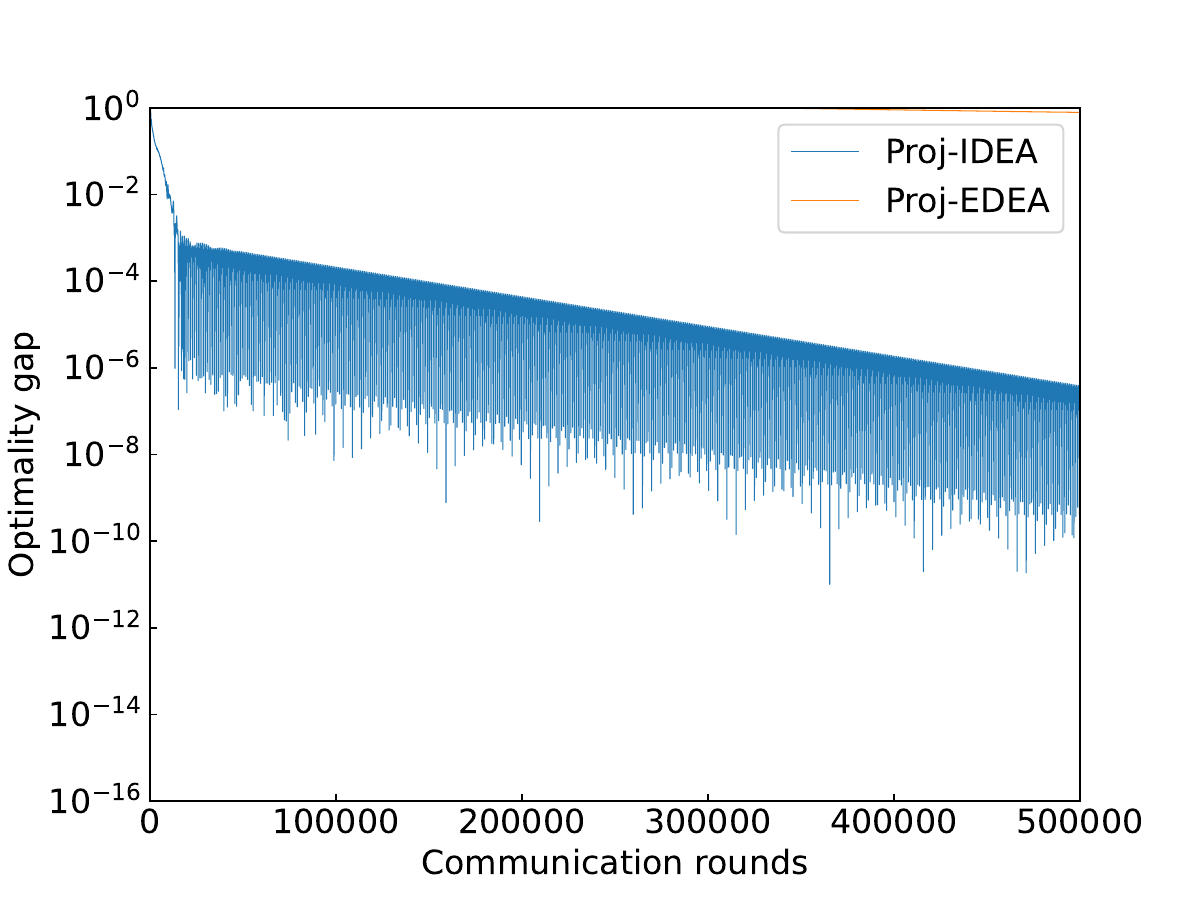}
            \end{minipage}
        }
        \subfigure[]{
            \begin{minipage}[b]{0.48\textwidth}
                \includegraphics[scale=0.2]{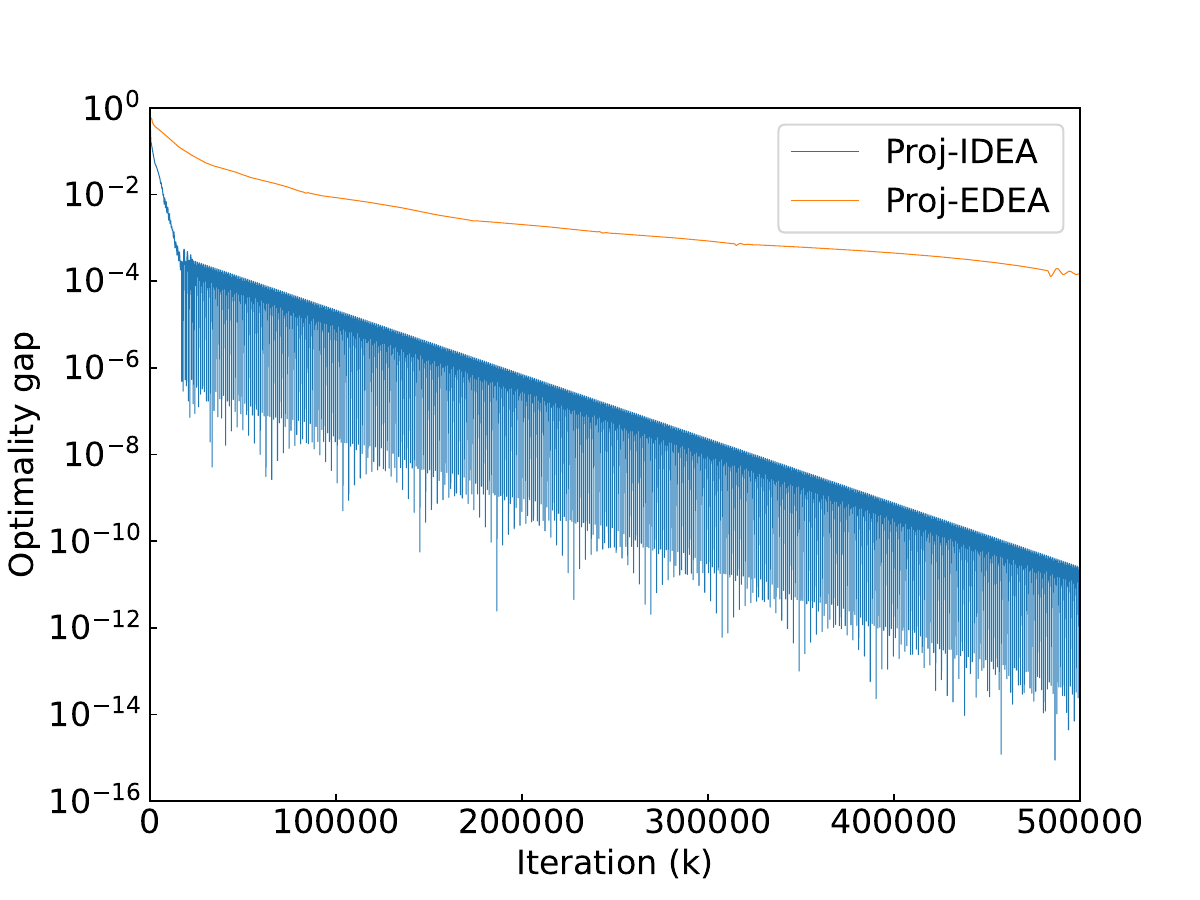}
                \includegraphics[scale=0.2]{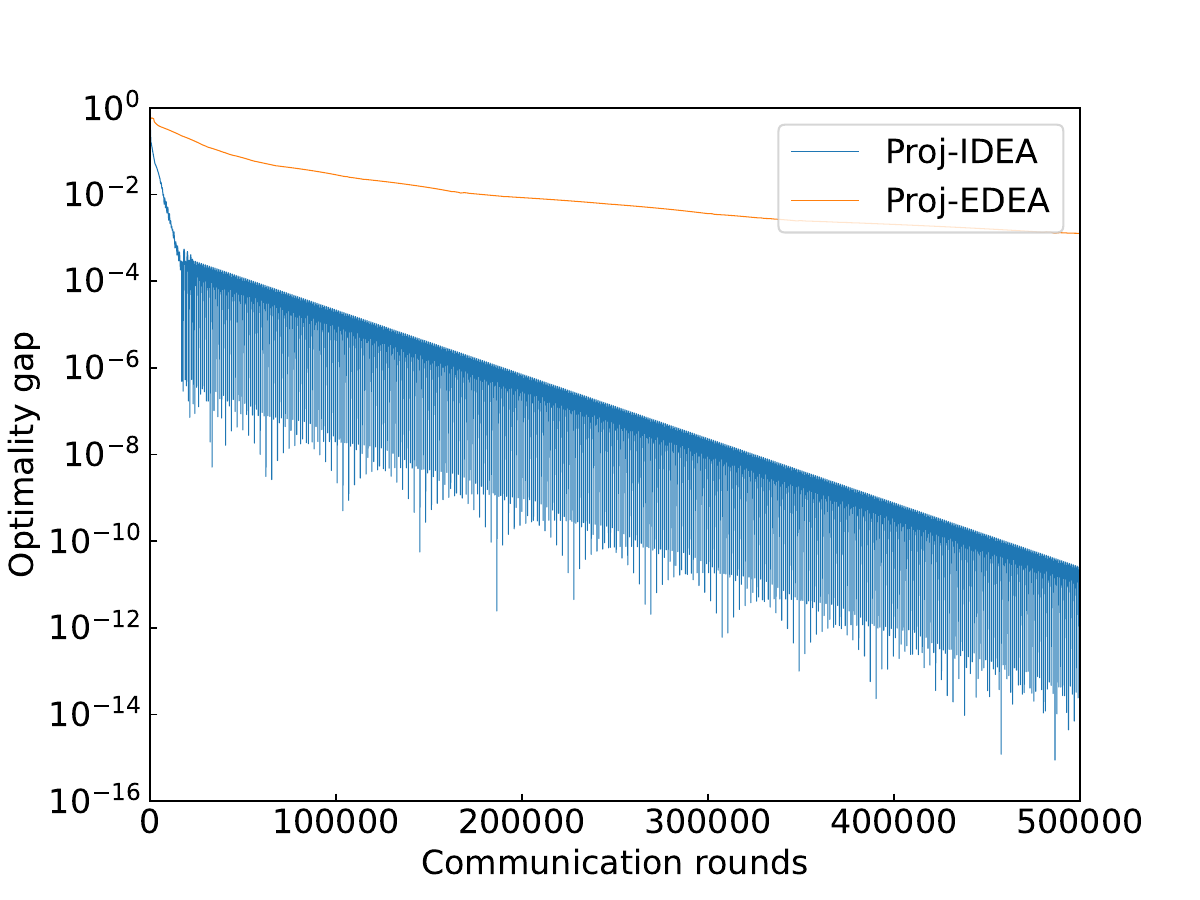}
            \end{minipage}
        }
        \subfigure[]{
            \begin{minipage}[b]{0.48\textwidth}
                \includegraphics[scale=0.2]{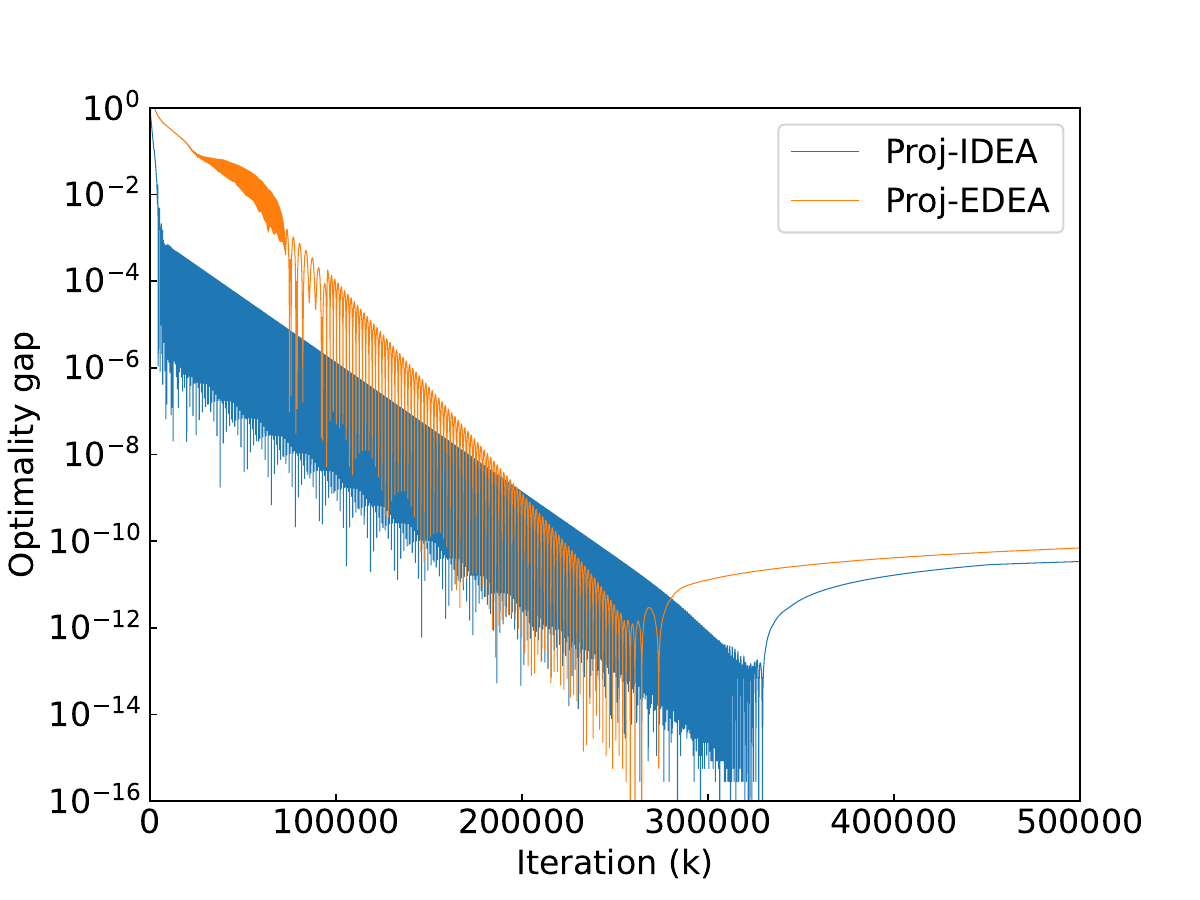}
                \includegraphics[scale=0.2]{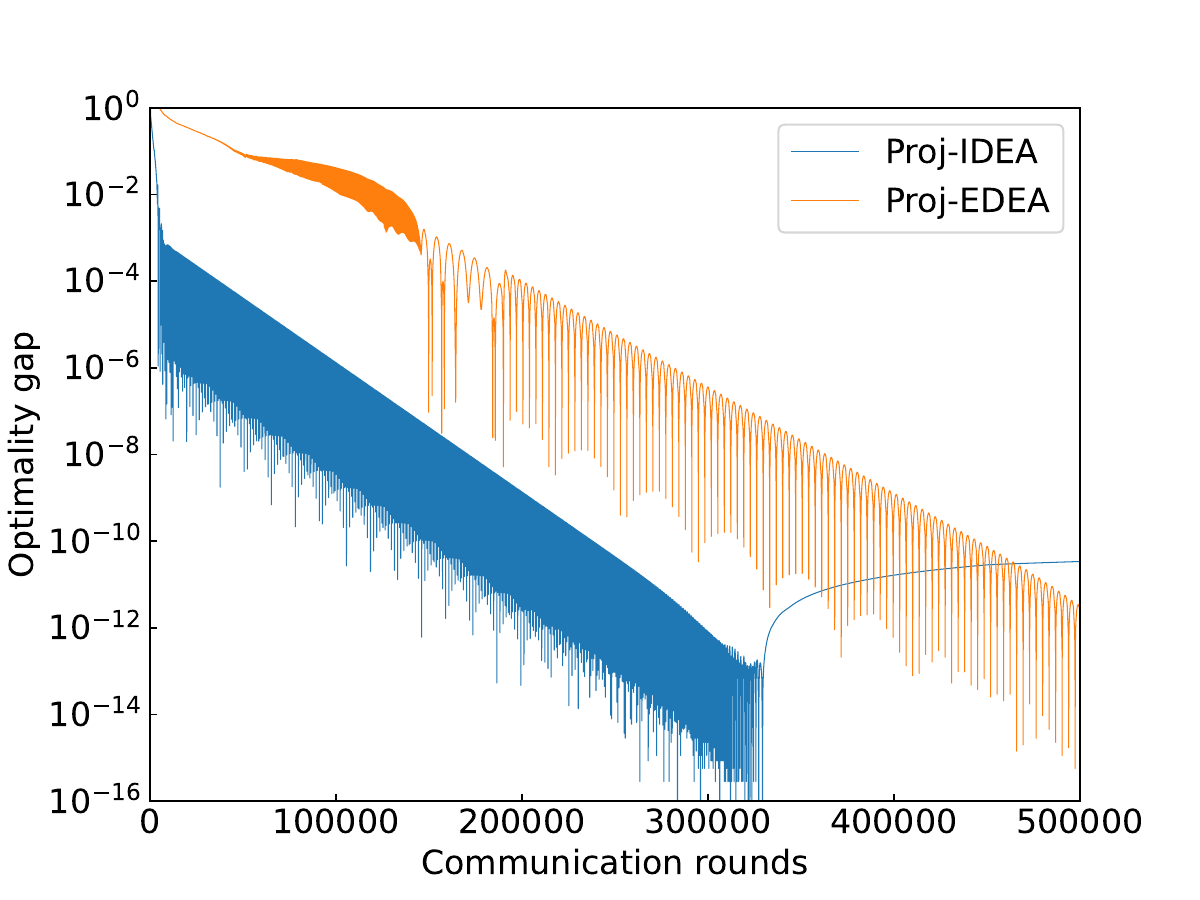}
            \end{minipage}
        }
        \subfigure[]{
            \begin{minipage}[b]{0.48\textwidth}
                \includegraphics[scale=0.2]{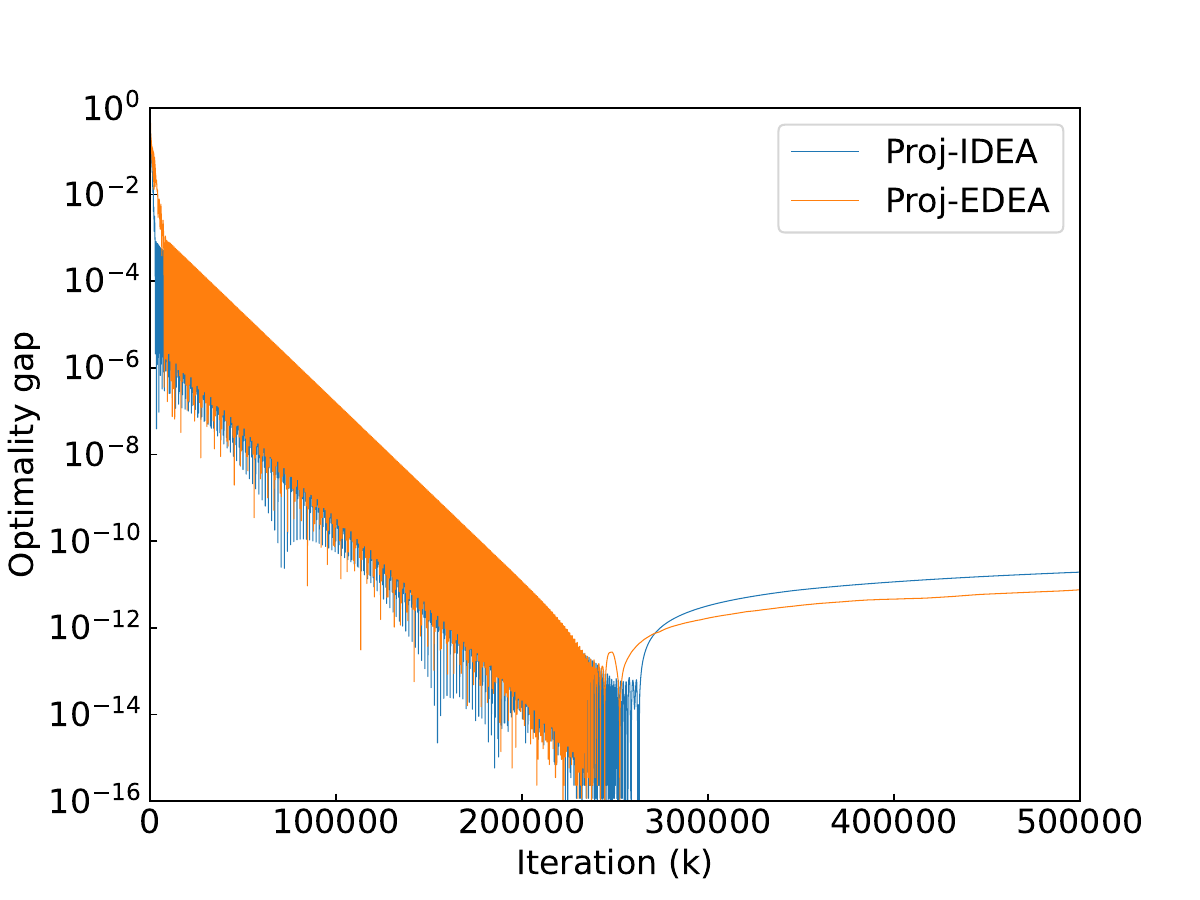}
                \includegraphics[scale=0.2]{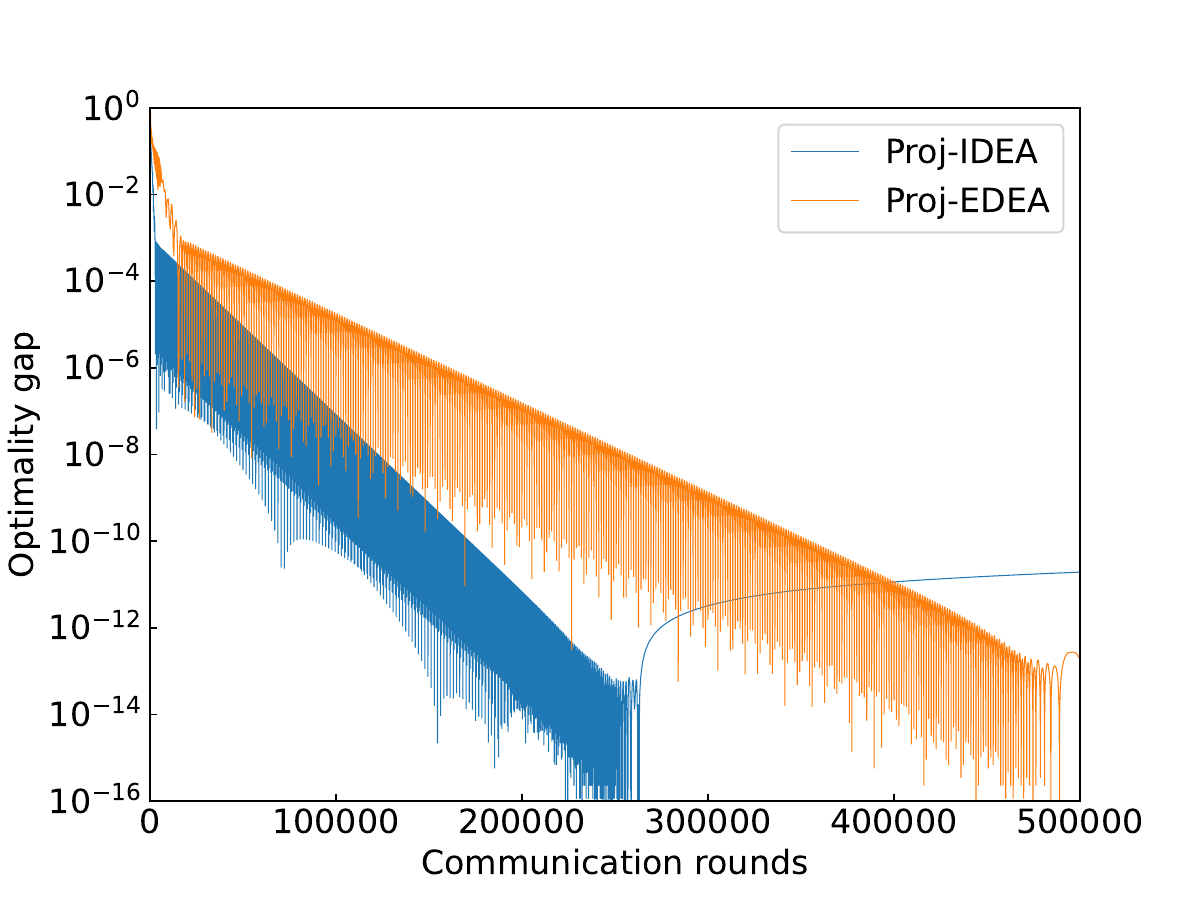}
            \end{minipage}
        }
        \caption {Experiment results of Case 1. (a) Circle graph, $\eta_2(L) \approx 0.02$. (b) Random graph with $p=0.05$, $\eta_2(L) \approx 0.21$. (c) Random graph with $p=0.1$, $\eta_2(L) \approx 0.54$. (d) Random graph with $p=0.3$, $\eta_2(L) \approx 6.34$.}
        \label{fig1}
    \end{center}
\end{figure*}

\begin{figure*}[tb]
    \begin{center}
        \subfigure[]{
            \begin{minipage}[b]{0.48\textwidth}
                \includegraphics[scale=0.2]{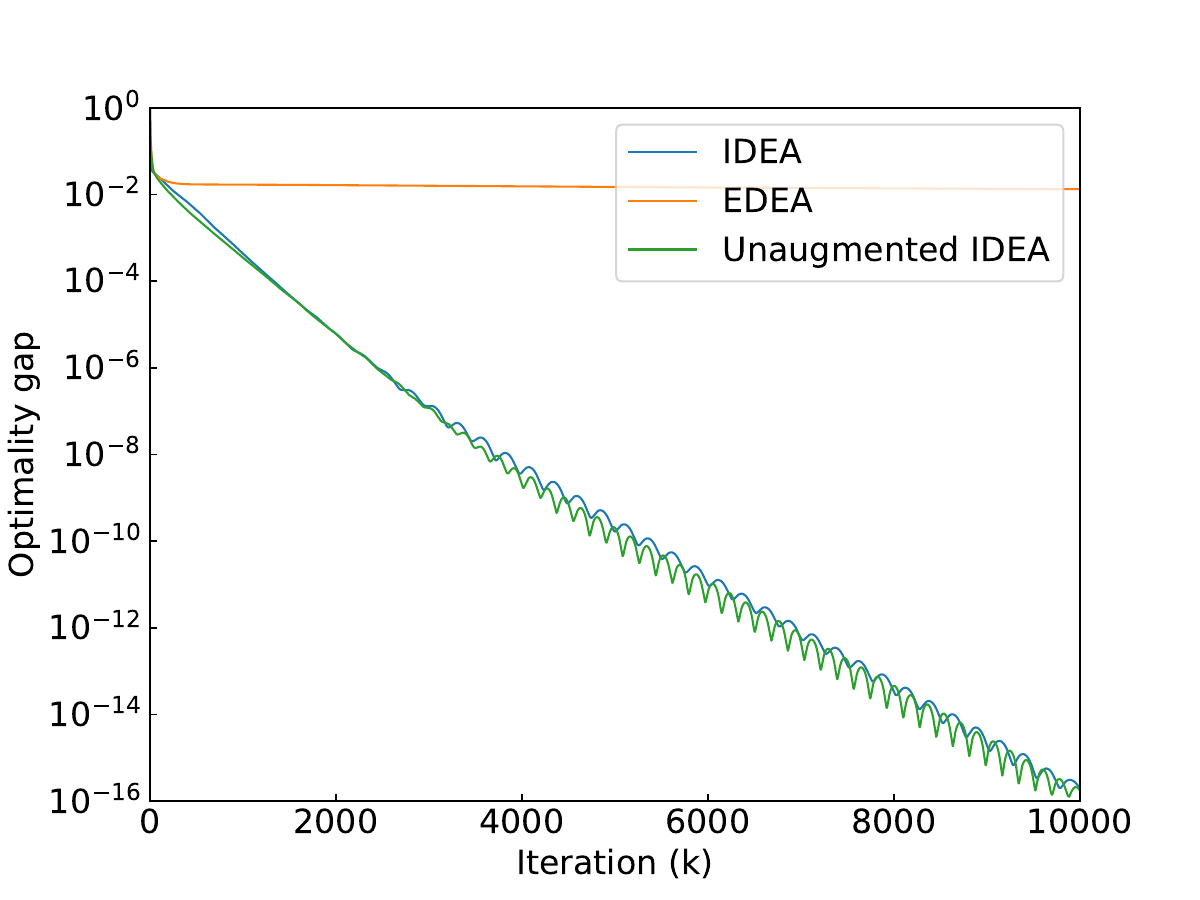}
                \includegraphics[scale=0.2]{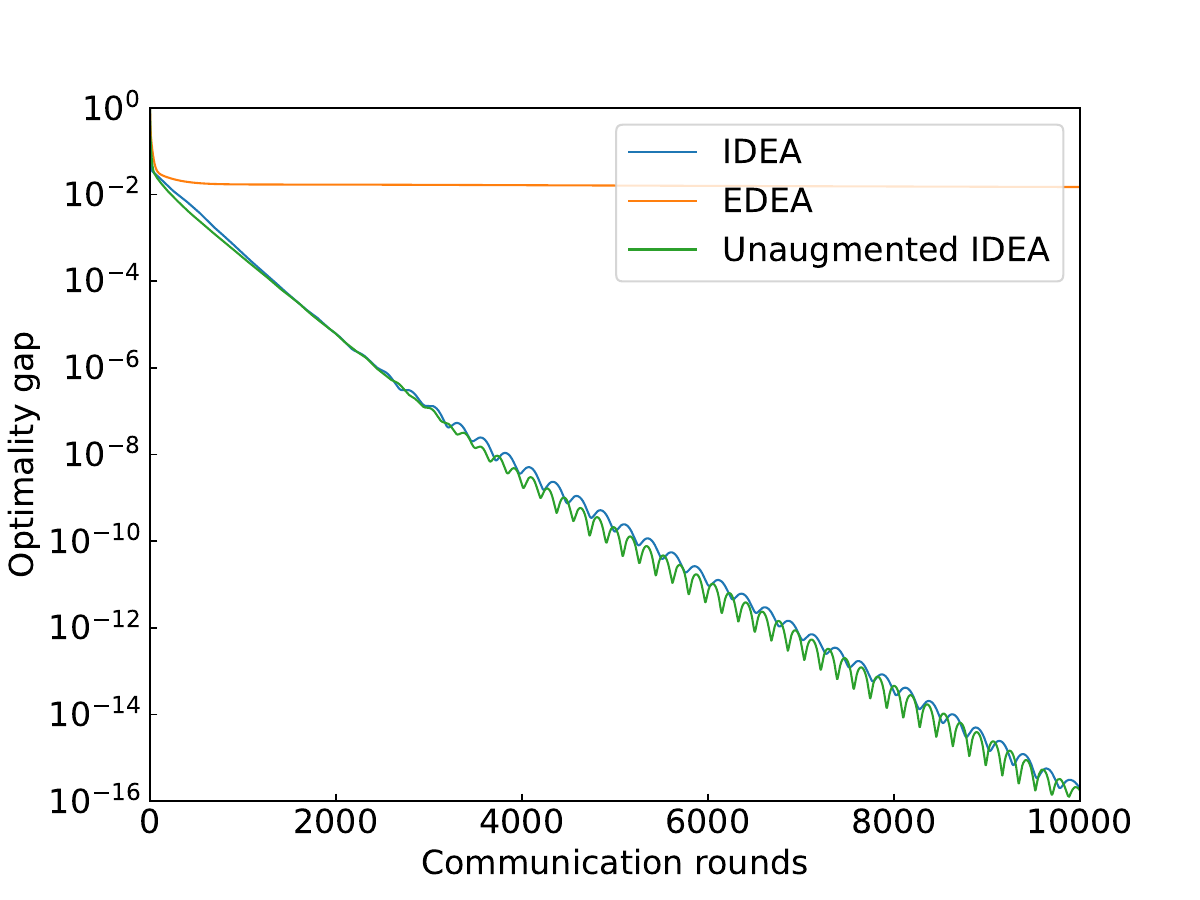}
            \end{minipage}
        }
        \subfigure[]{
            \begin{minipage}[b]{0.48\textwidth}
                \includegraphics[scale=0.2]{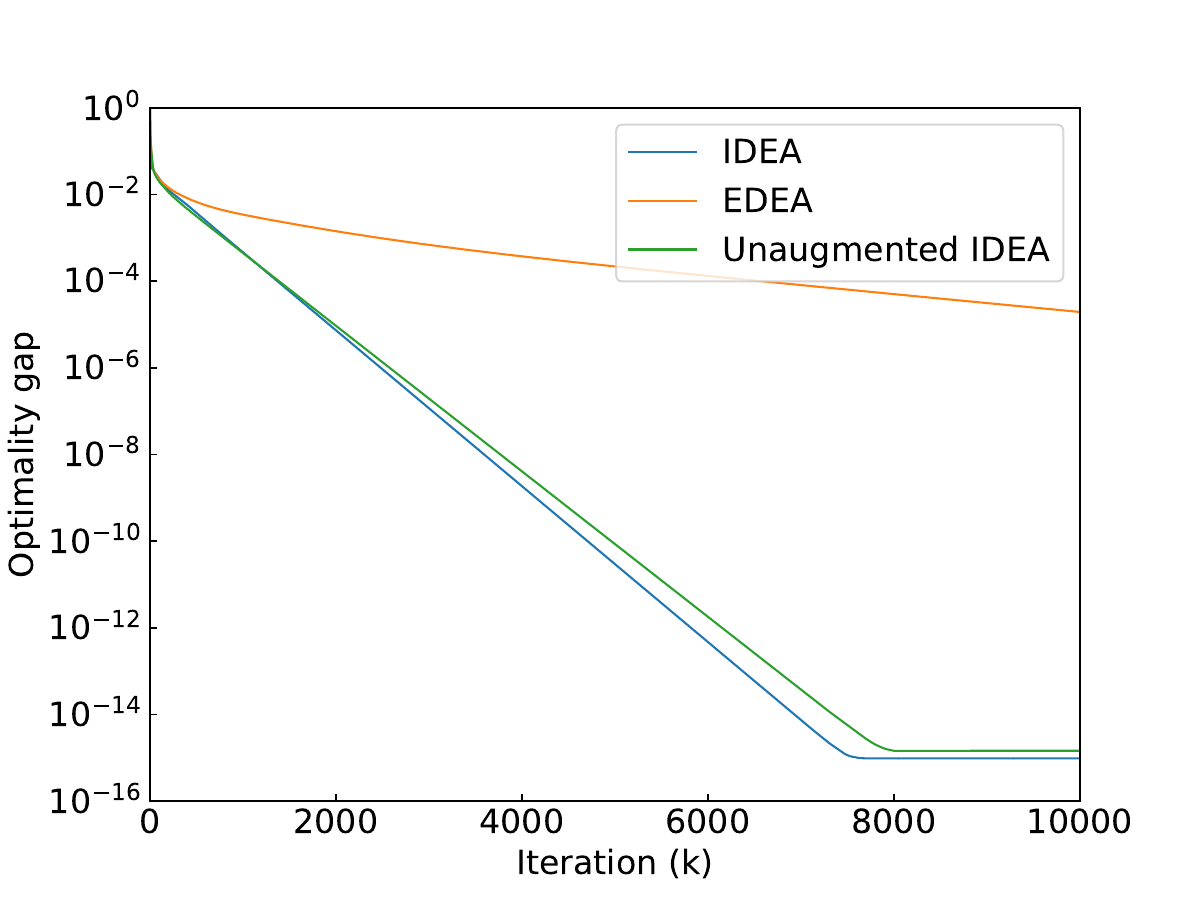}
                \includegraphics[scale=0.2]{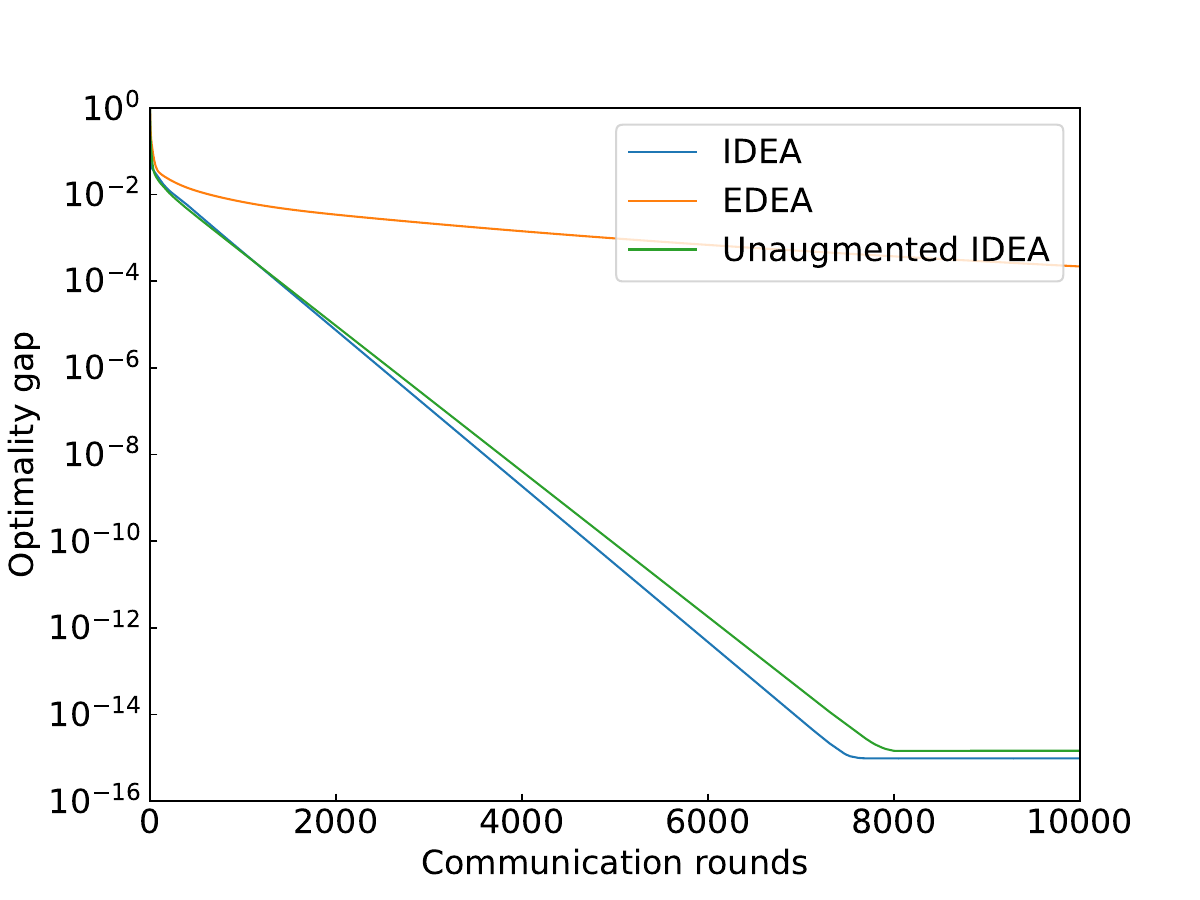}
            \end{minipage}
        }
        \subfigure[]{
            \begin{minipage}[b]{0.48\textwidth}
                \includegraphics[scale=0.2]{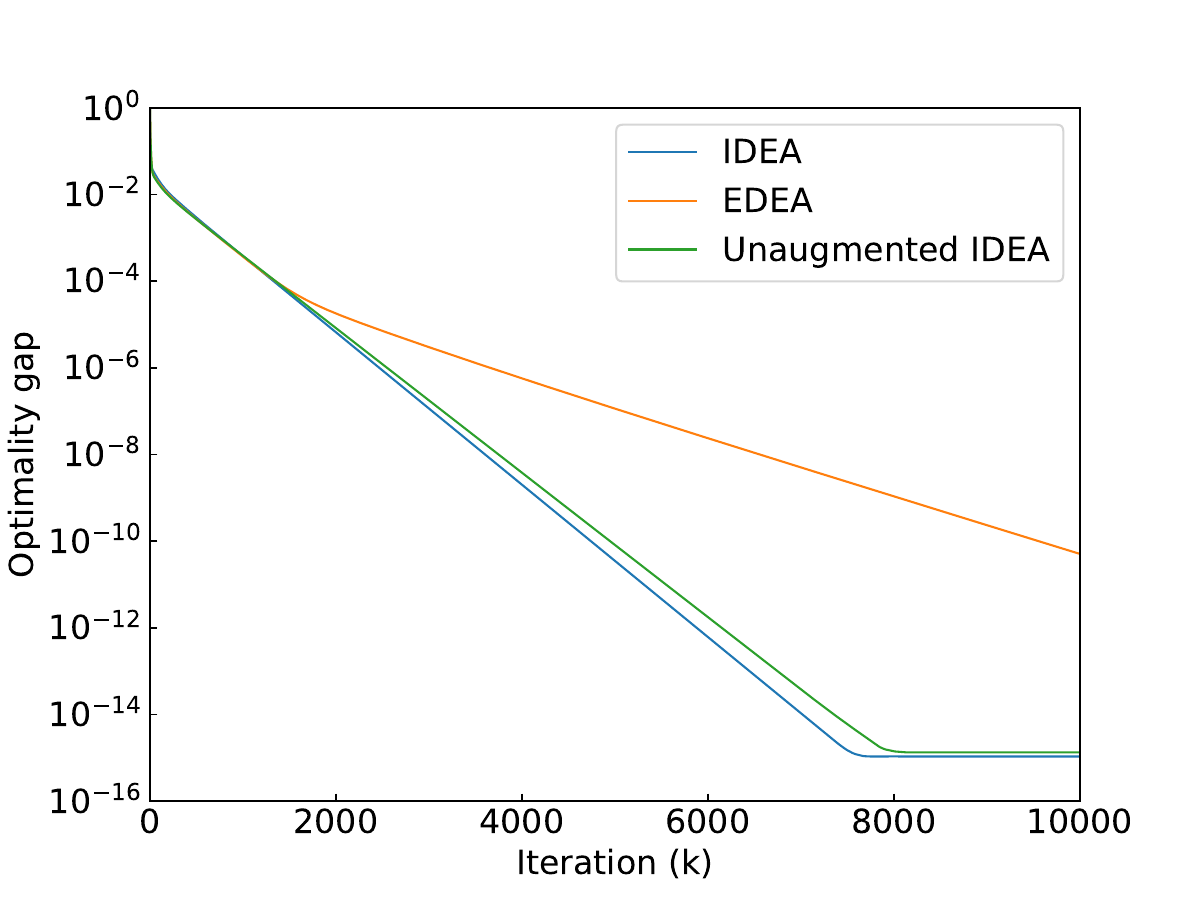}
                \includegraphics[scale=0.2]{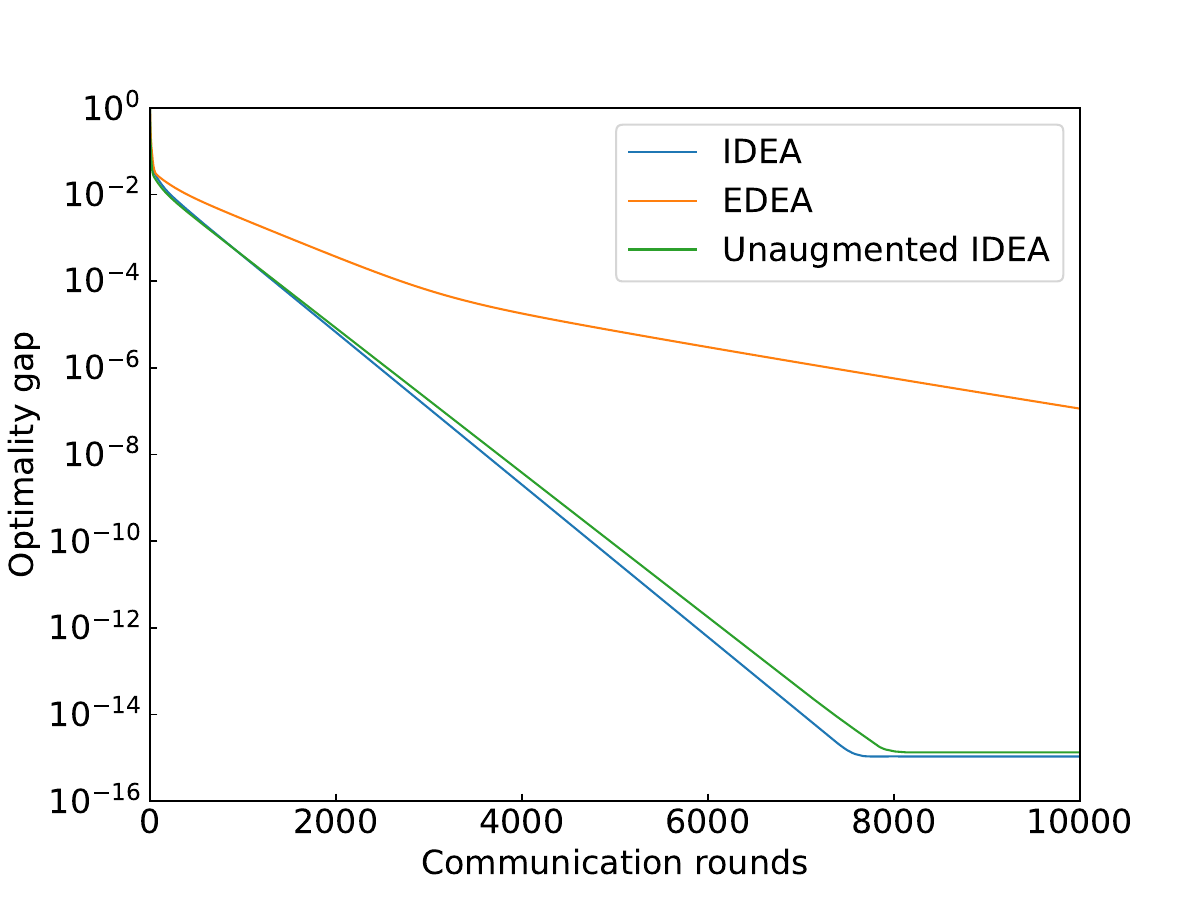}
            \end{minipage}
        }
        \subfigure[]{
            \begin{minipage}[b]{0.48\textwidth}
                \includegraphics[scale=0.2]{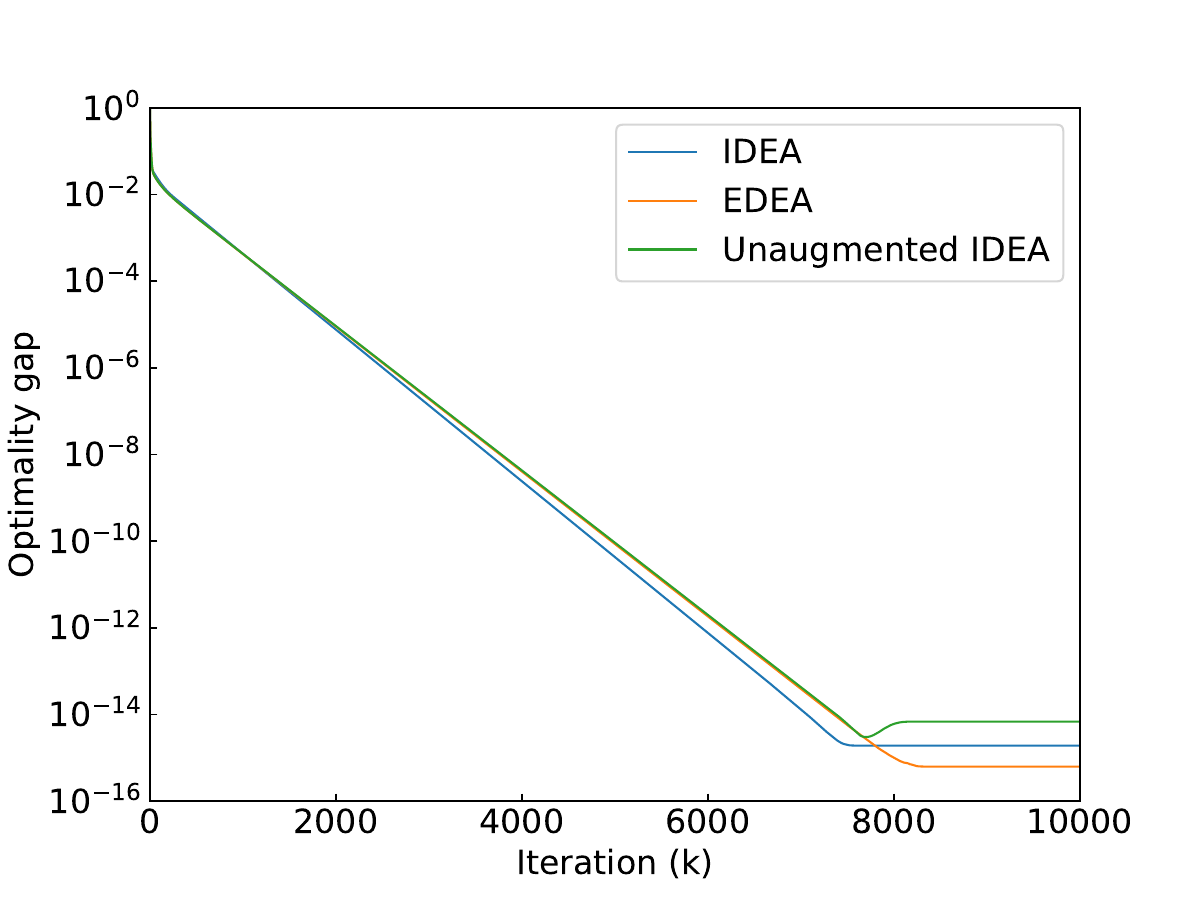}
                \includegraphics[scale=0.2]{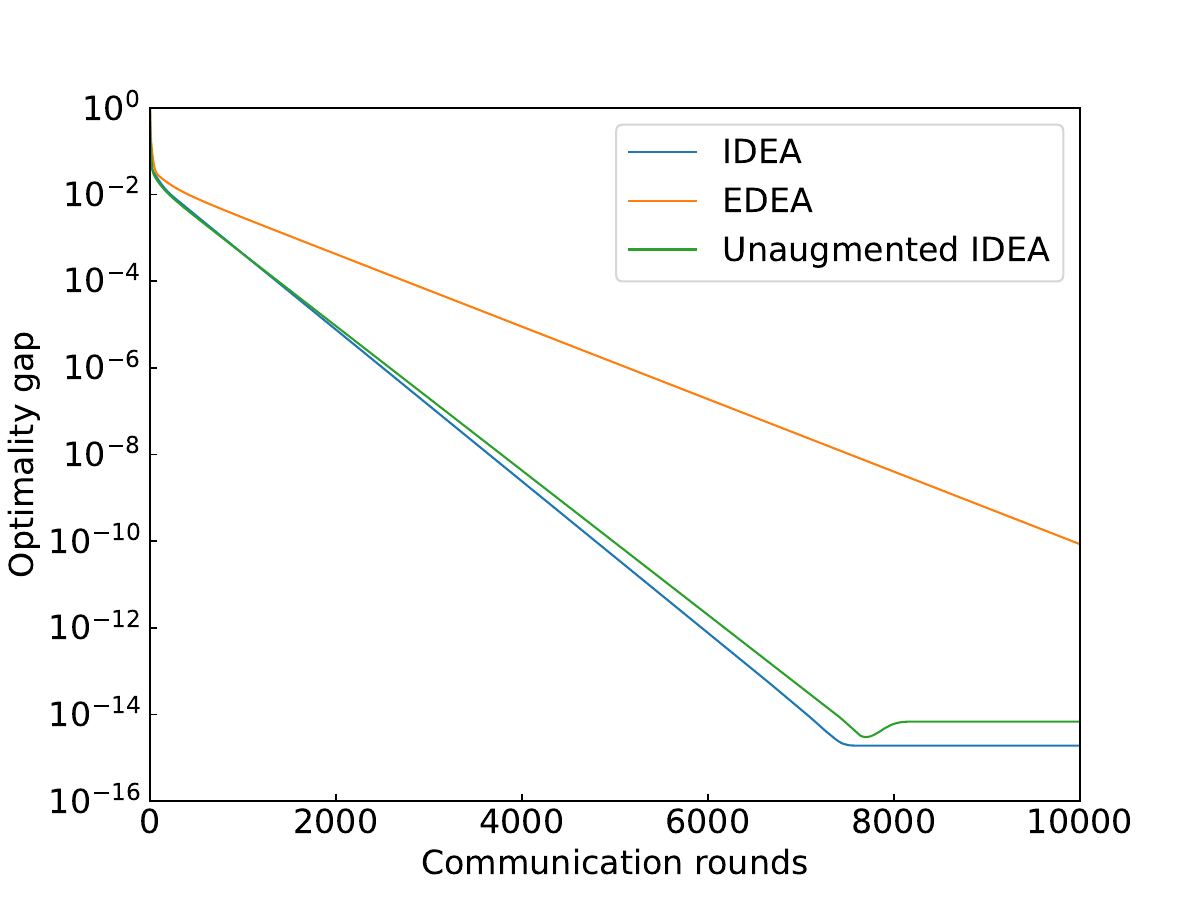}
            \end{minipage}
        }
        \caption {Experiment results of Case 2. (a) Circle graph, $\eta_2(L) \approx 0.02$. (b) Random graph with $p=0.05$, $\eta_2(L) \approx 0.21$. (c) Random graph with $p=0.1$, $\eta_2(L) \approx 0.54$. (d) Random graph with $p=0.3$, $\eta_2(L) \approx 6.34$.}
        \label{fig2}
    \end{center}
\end{figure*}

\begin{figure*}[tb]
    \begin{center}
        \subfigure[]{
            \begin{minipage}[b]{0.48\textwidth}
                \includegraphics[scale=0.2]{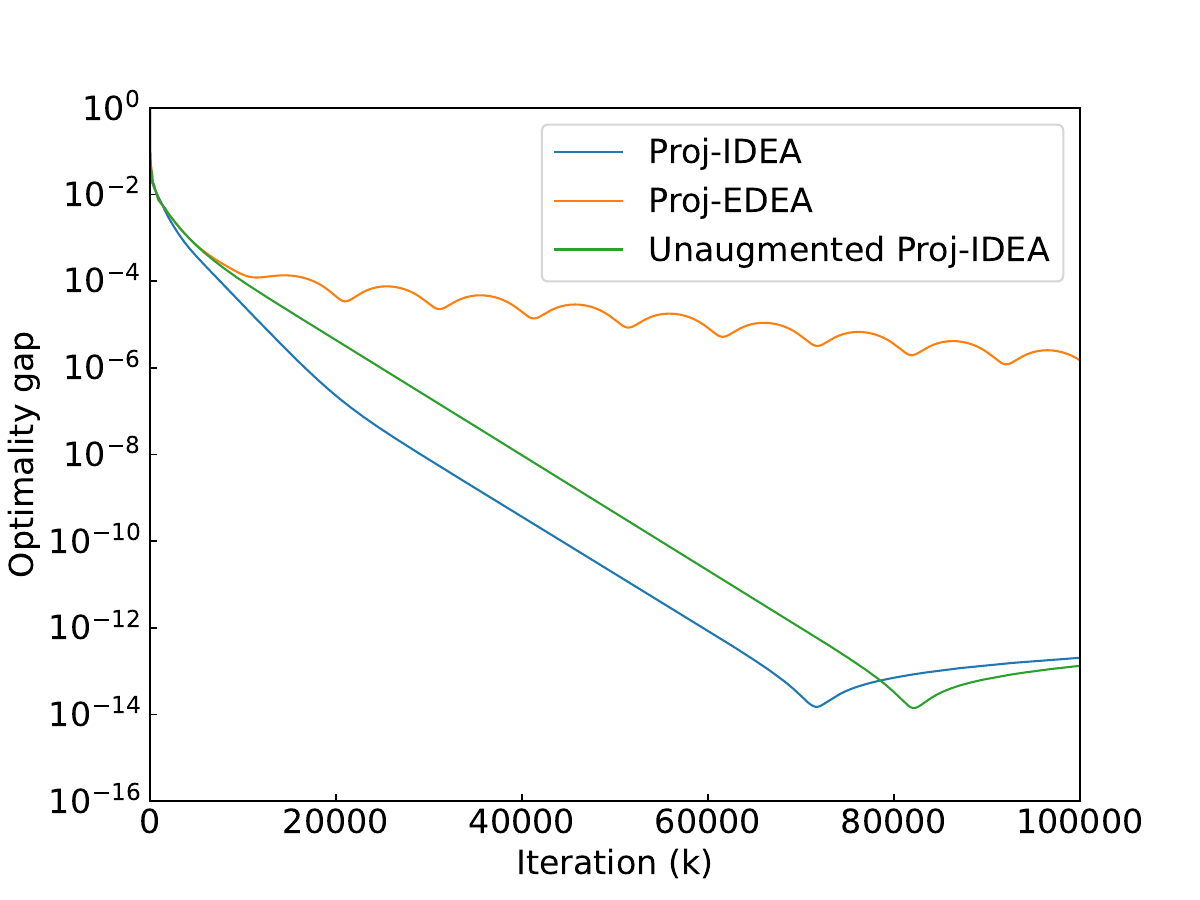}
                \includegraphics[scale=0.2]{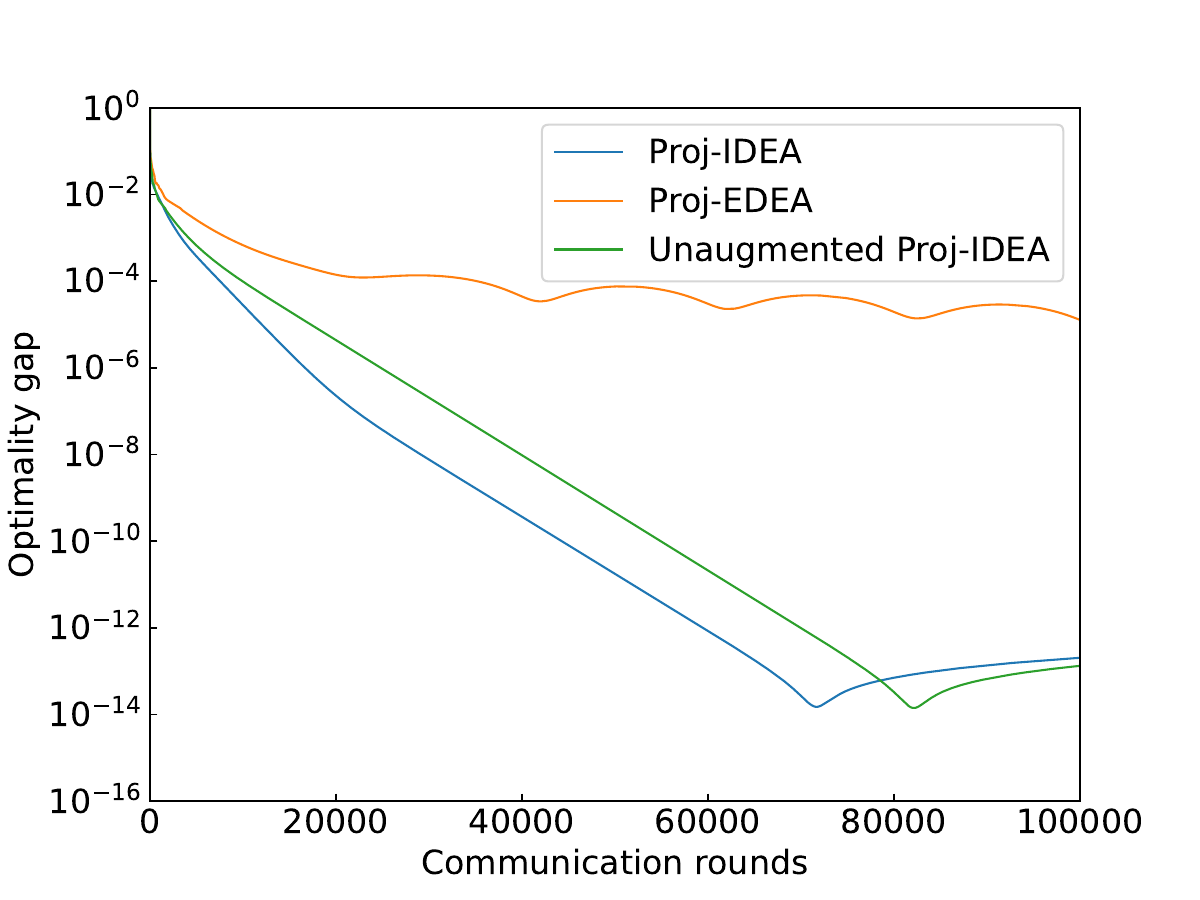}
            \end{minipage}
        }
        \subfigure[]{
            \begin{minipage}[b]{0.48\textwidth}
                \includegraphics[scale=0.2]{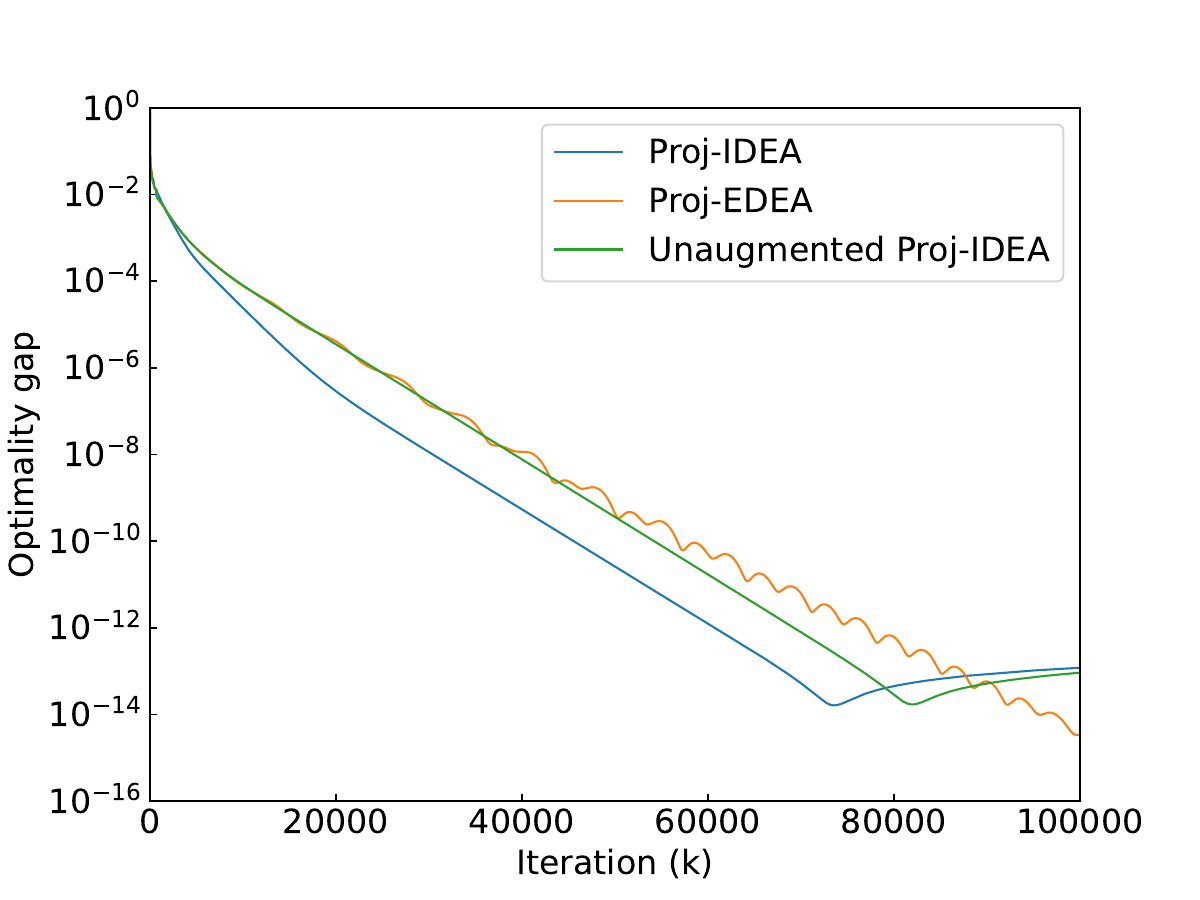}
                \includegraphics[scale=0.2]{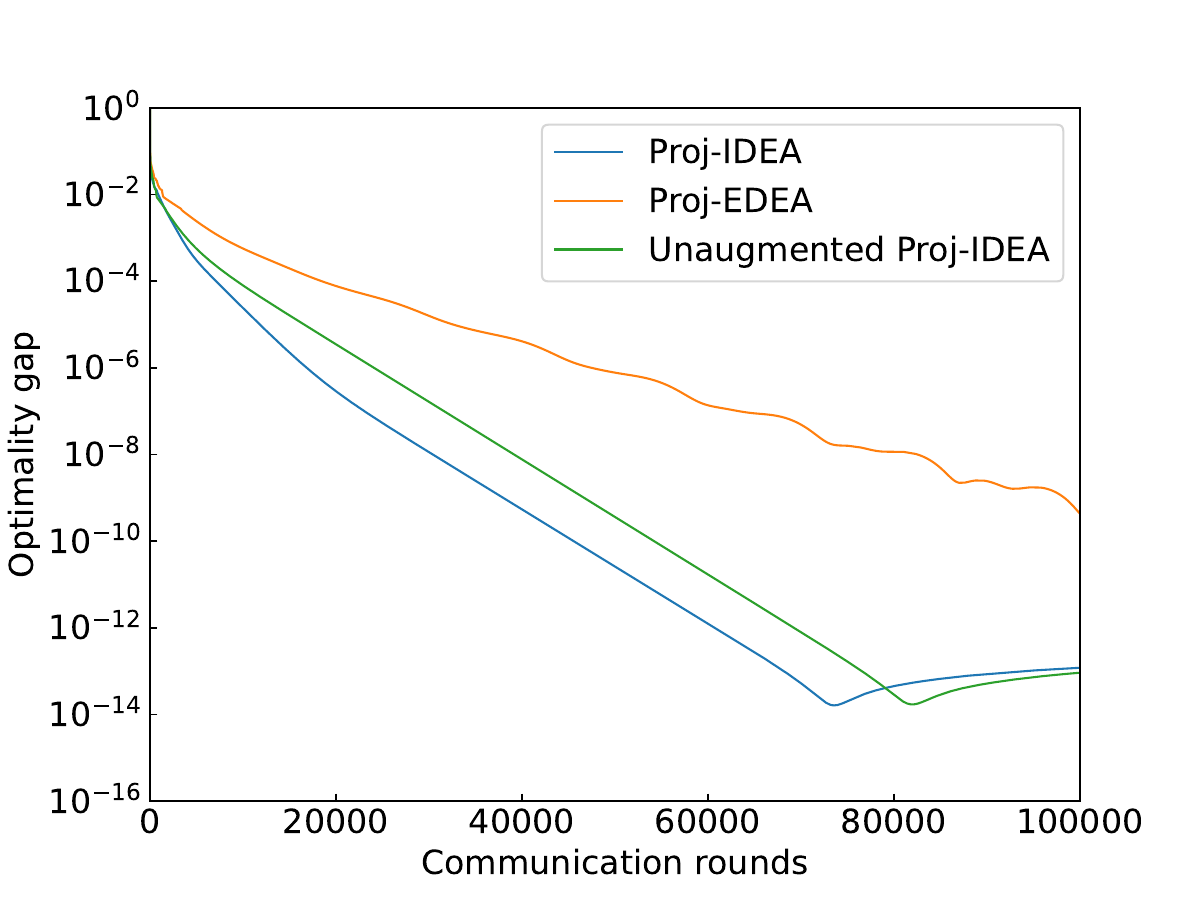}
            \end{minipage}
        }
        \subfigure[]{
            \begin{minipage}[b]{0.48\textwidth}
                \includegraphics[scale=0.2]{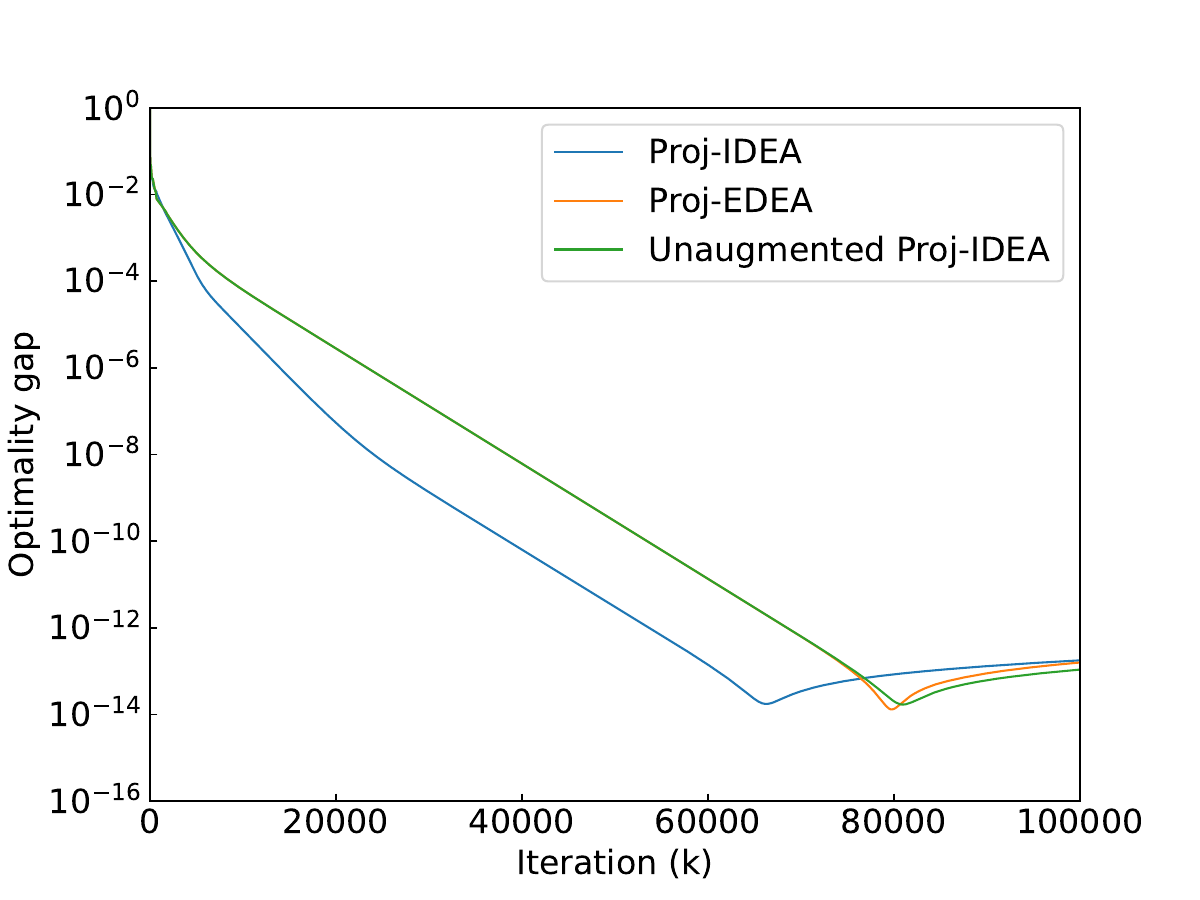}
                \includegraphics[scale=0.2]{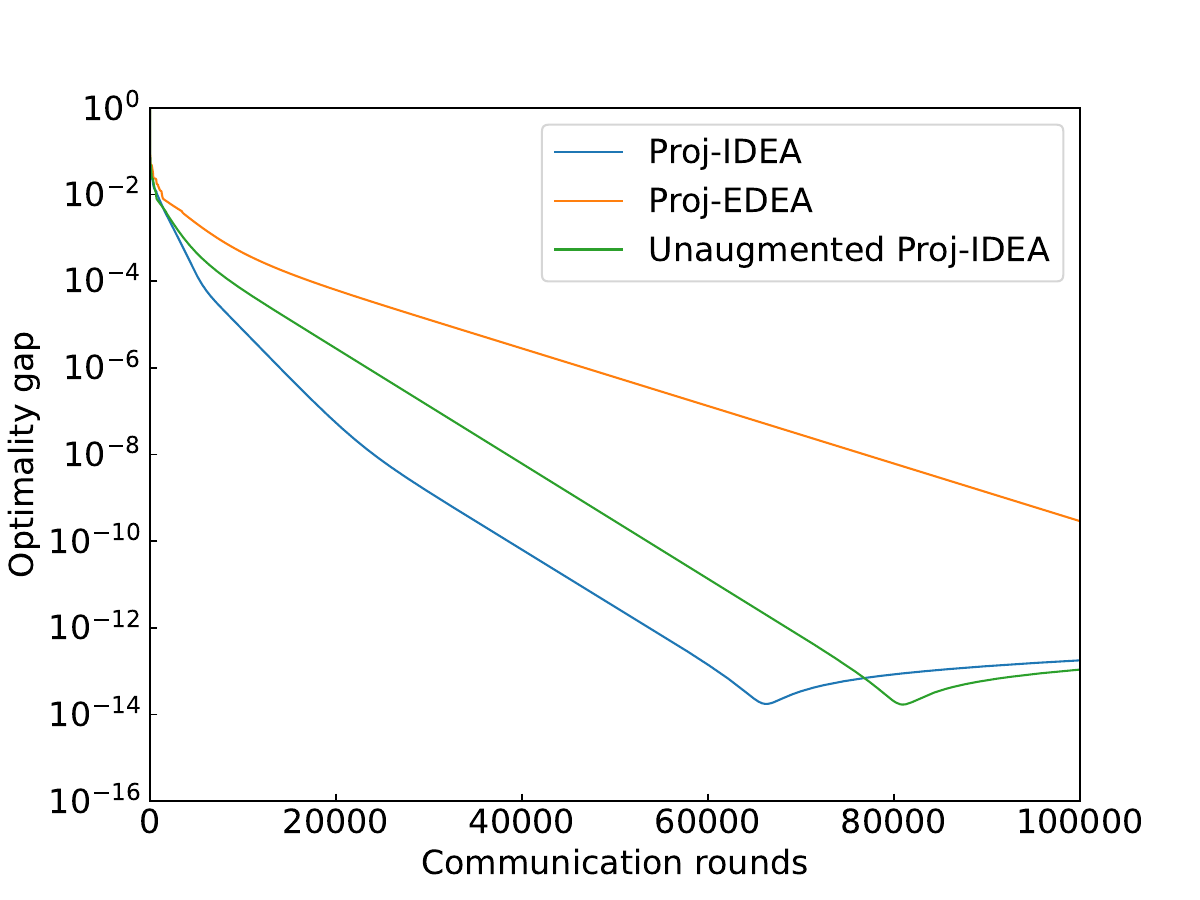}
            \end{minipage}
        }
        \subfigure[]{
            \begin{minipage}[b]{0.48\textwidth}
                \includegraphics[scale=0.2]{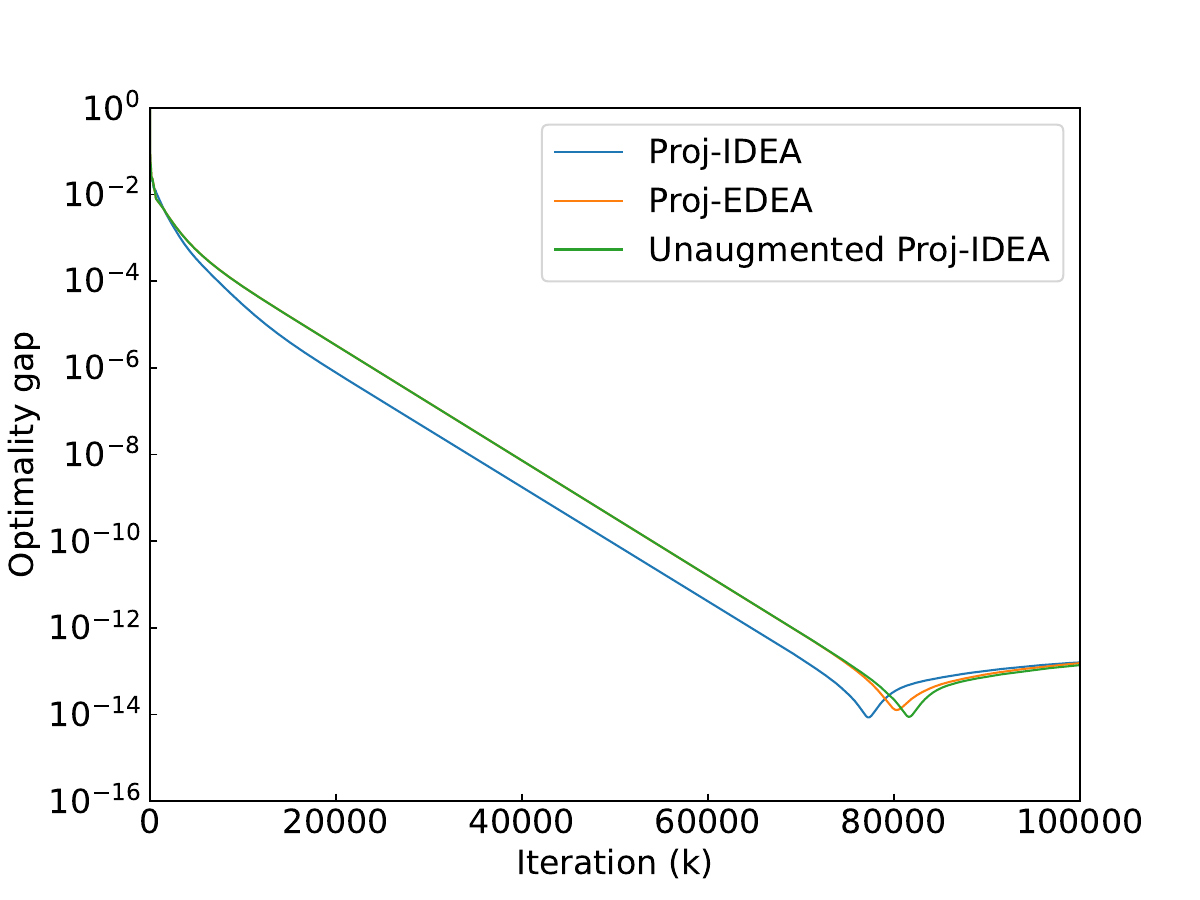}
                \includegraphics[scale=0.2]{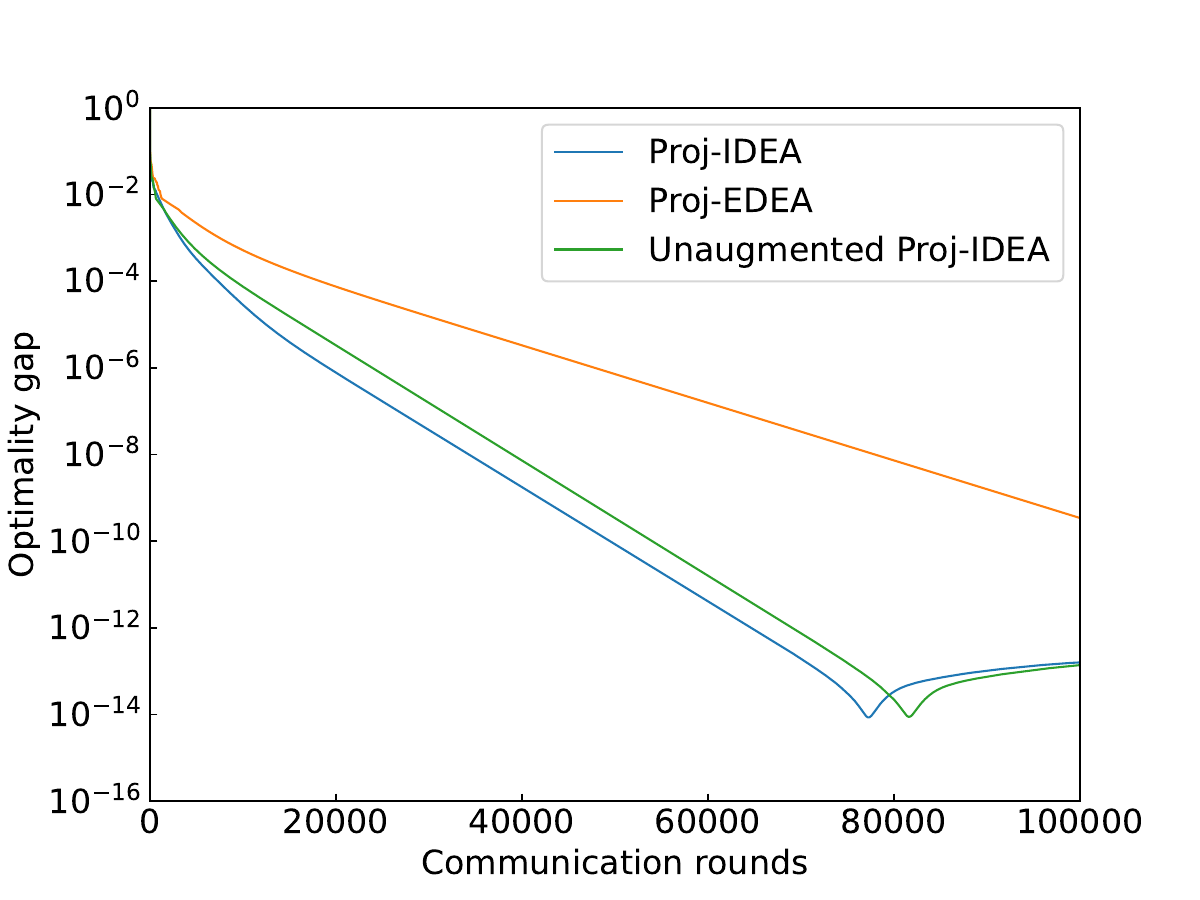}
            \end{minipage}
        }
        \caption {Experiment results of Case 3. (a) Directed circle graph. (b) Directed exponential graph with $e=2$. (c) Directed exponential graph with $e=4$. (d) Directed exponential graph with $e=6$.}
        \label{fig3}
    \end{center}
\end{figure*}

\begin{figure*}[tb]
    \begin{center}
        \subfigure[]{
            \begin{minipage}[b]{0.48\textwidth}
                \includegraphics[scale=0.2]{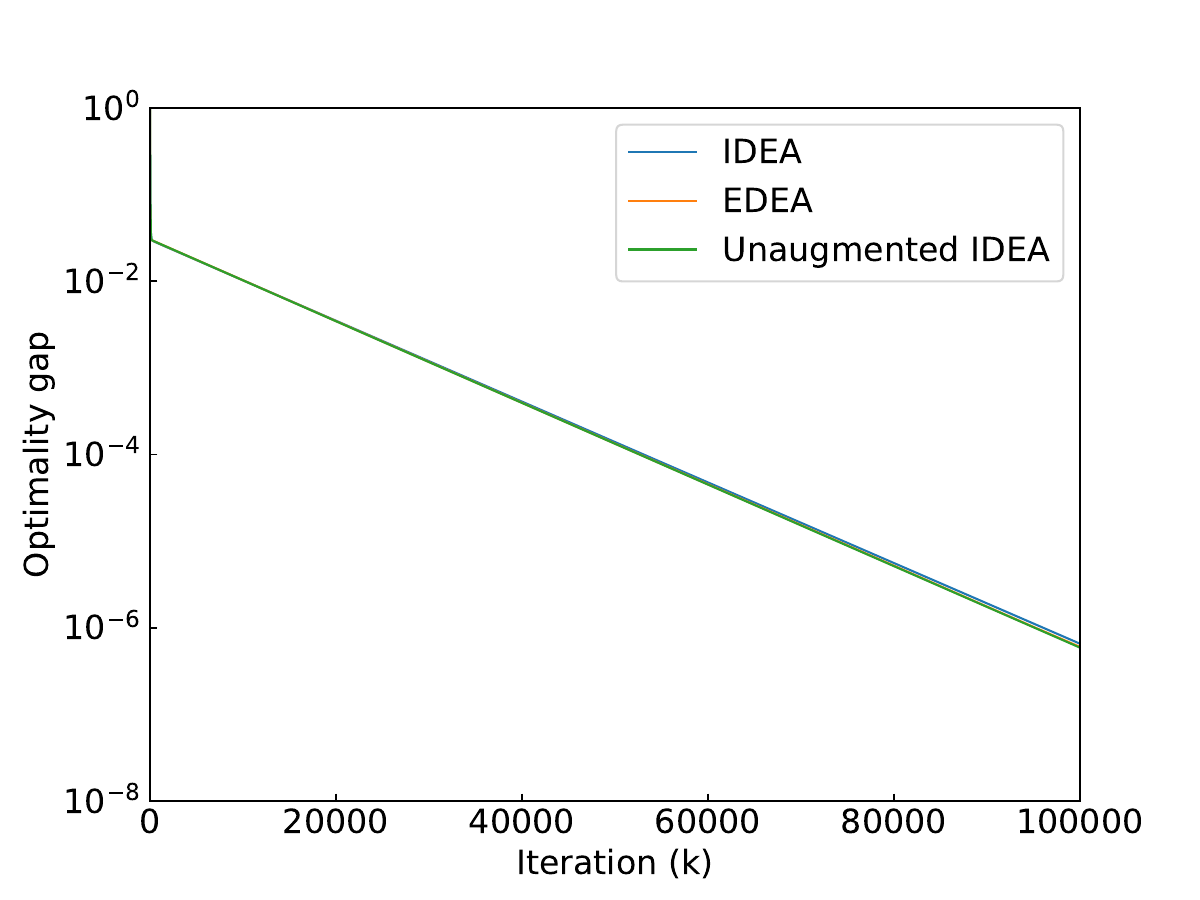}
                \includegraphics[scale=0.2]{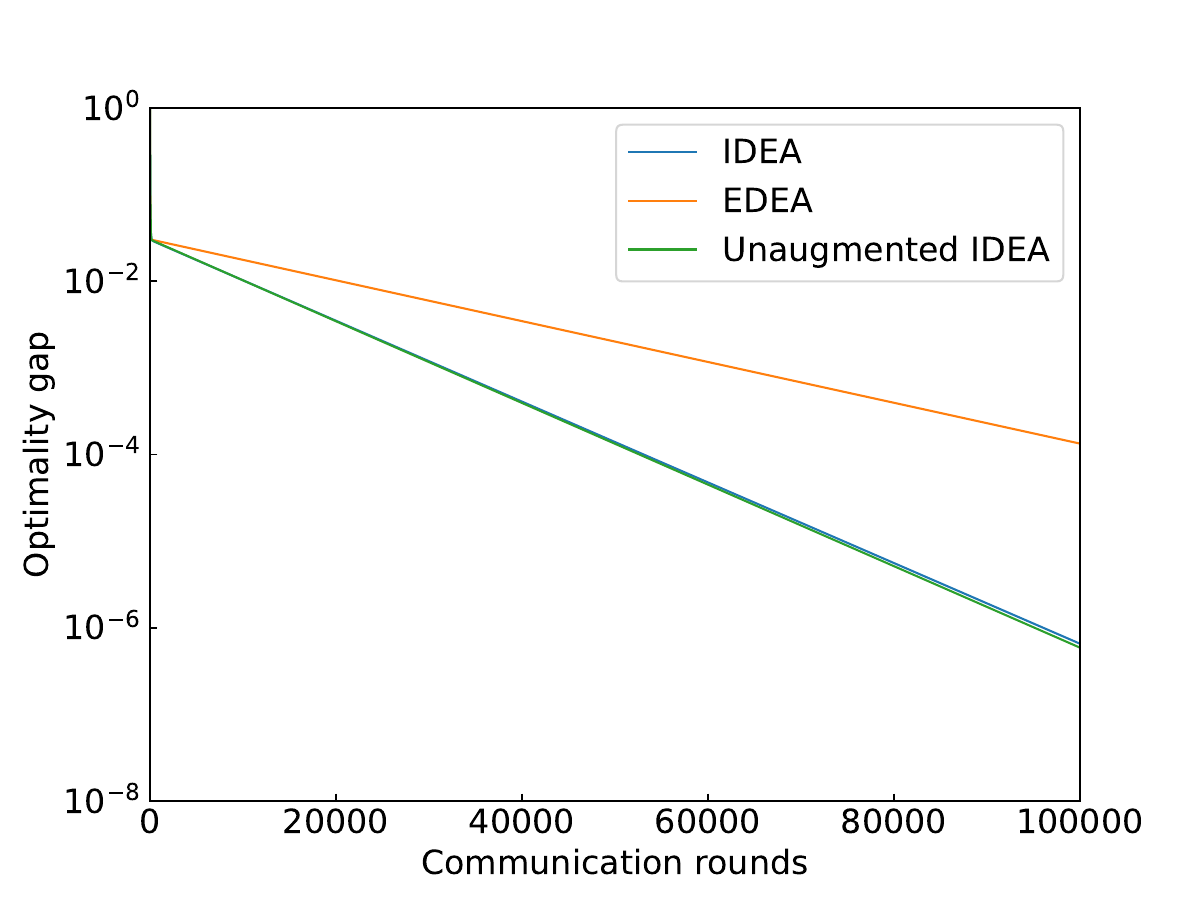}
            \end{minipage}
        }
        \subfigure[]{
            \begin{minipage}[b]{0.48\textwidth}
                \includegraphics[scale=0.2]{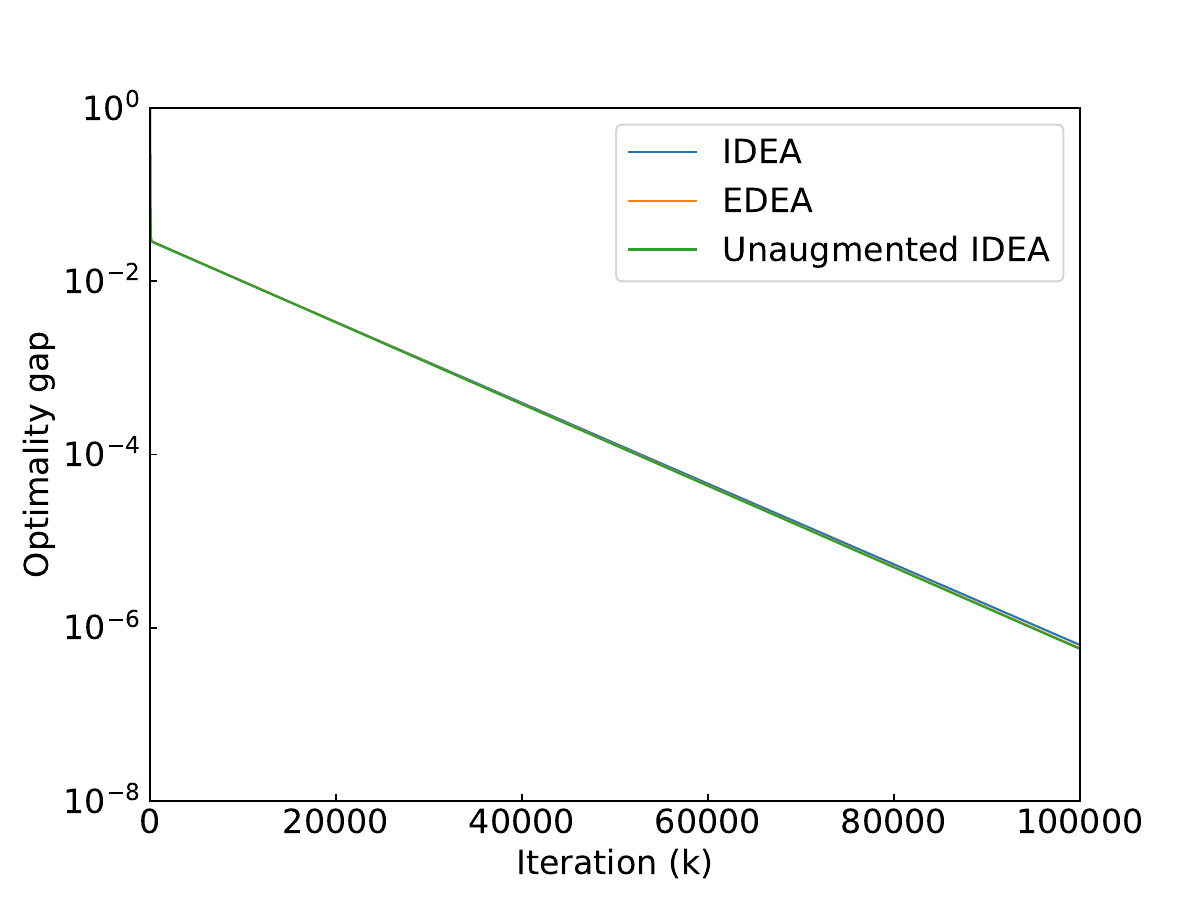}
                \includegraphics[scale=0.2]{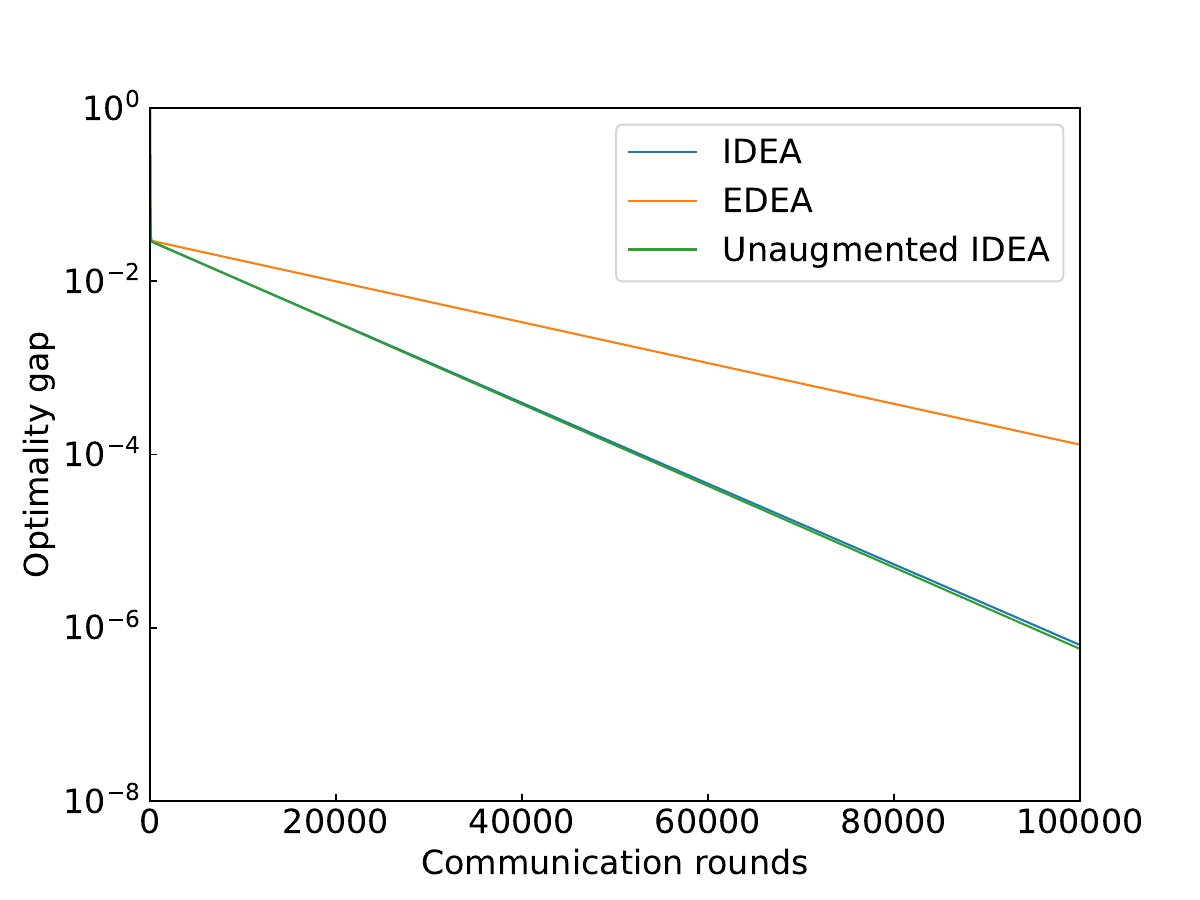}
            \end{minipage}
        }
        \subfigure[]{
            \begin{minipage}[b]{0.48\textwidth}
                \includegraphics[scale=0.2]{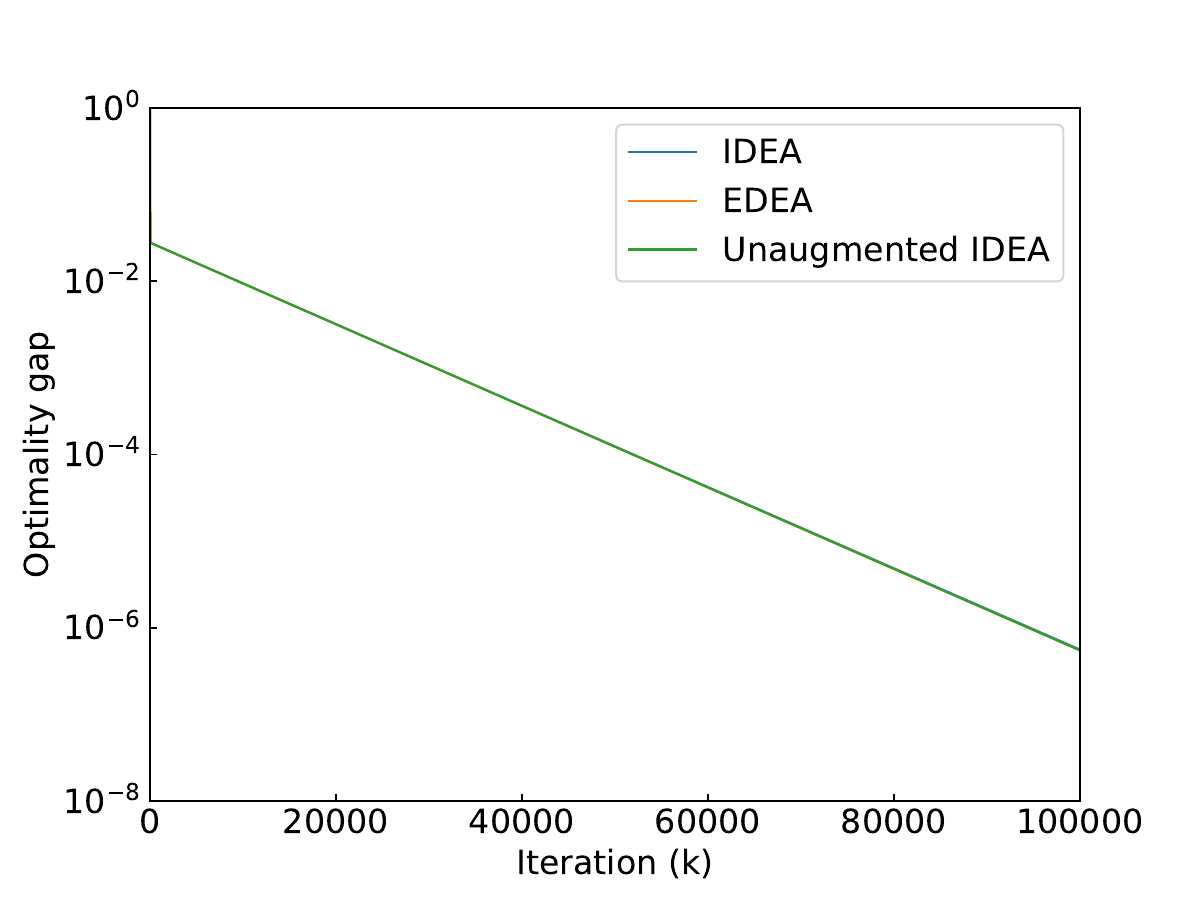}
                \includegraphics[scale=0.2]{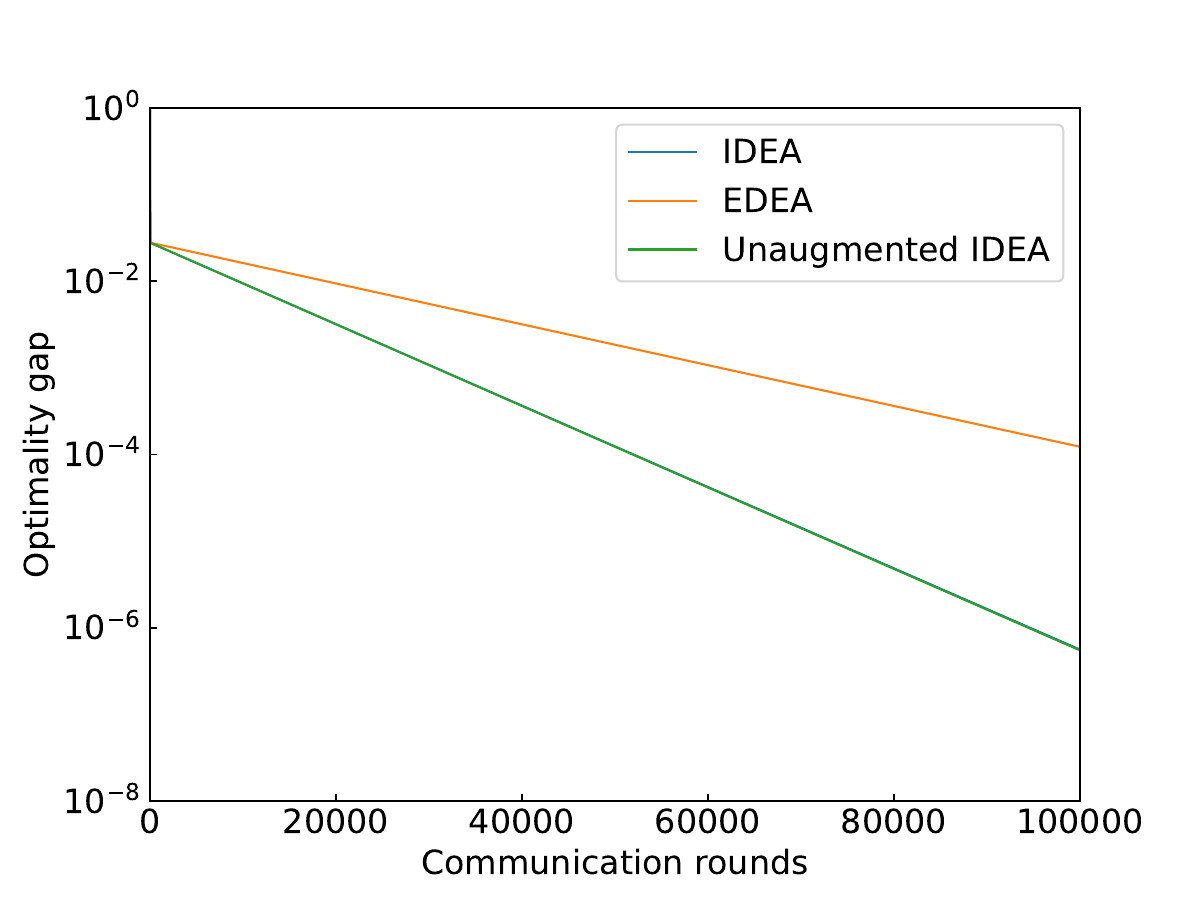}
            \end{minipage}
        }
        \subfigure[]{
            \begin{minipage}[b]{0.48\textwidth}
                \includegraphics[scale=0.2]{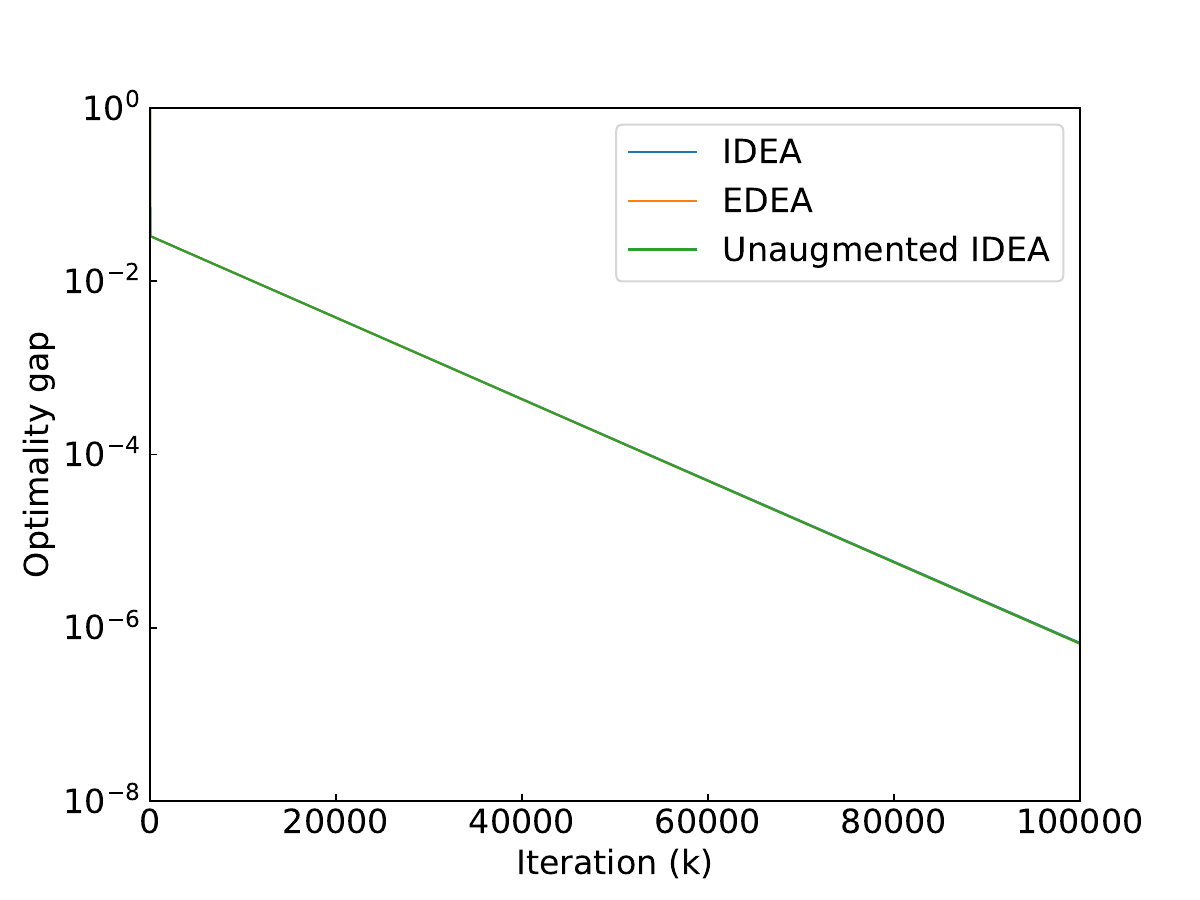}
                \includegraphics[scale=0.2]{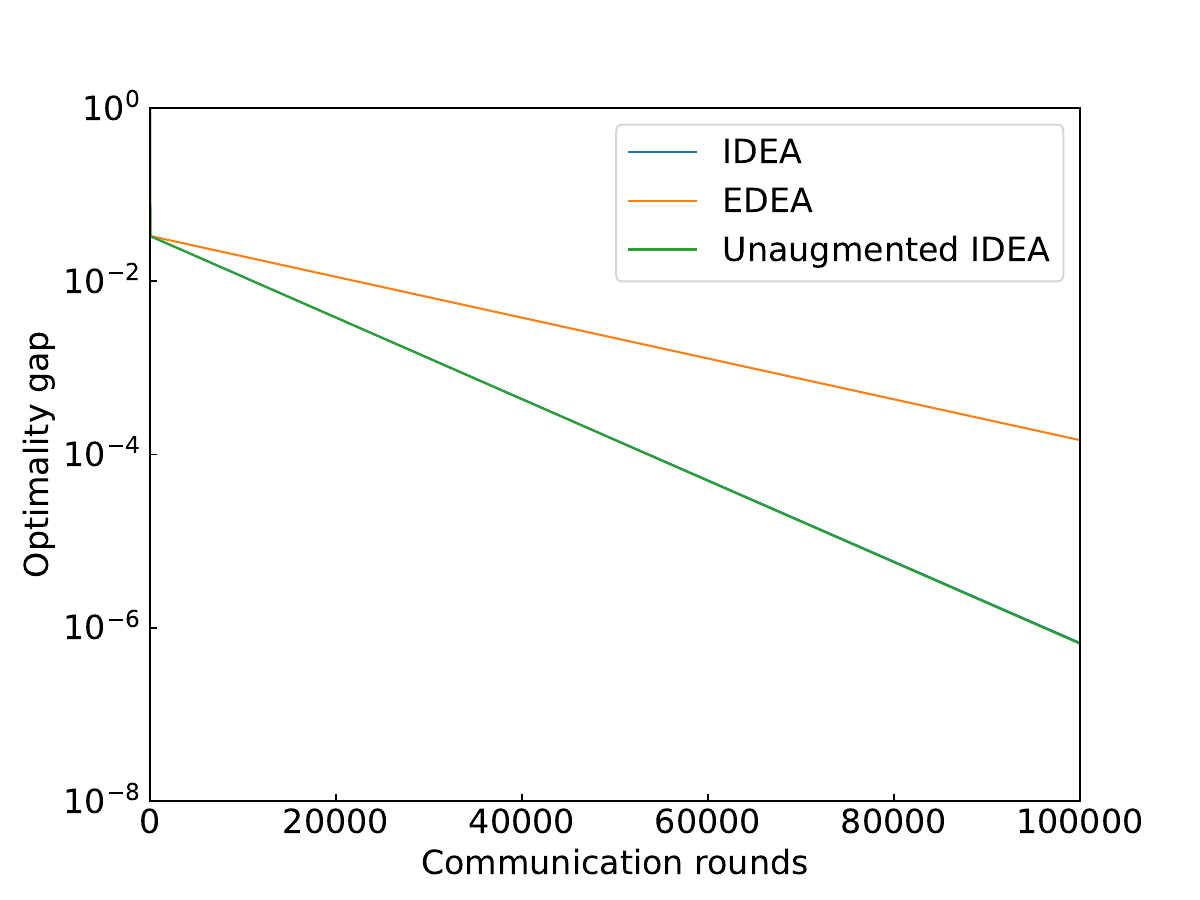}
            \end{minipage}
        }
        \caption {Experiment results of Case 4. (a) Directed circle graph. (b) Directed exponential graph with $e=2$. (c) Directed exponential graph with $e=4$. (d) Directed exponential graph with $e=6$.}
        \label{fig4}
    \end{center}
\end{figure*}
\begin{corollary} \label{cor2}
    Suppose \cref{dif} holds, the time-varying $\mG$ is always strongly connected and weight-balanced and its adjacency matrix is piecewise constant and uniformly bounded, and other conditions are the same with \cref{the4}, except that $\beta$ satisfies
    \eqe{
        \beta \geq \frac{2(\varphi+1)^2\alpha^2+1}{2\varphi\alpha\min_{s \in \mathcal{S}}\eta_2(\hat{L}_{s})},
        \nonumber
    }
    where $\mathcal{S}$ is the index set of all possible structures of $\mG$. Then, for any $(\mx(0), \ml(0), \mz(0)) \in \mR^d \times \mR^{np} \times \mR^{np}$ satisfies $\sum_{i=1}^nz_i(0) =\0$, $(\mx(t), \ml(t),\mz(t))$ driven by IDEA converges exponentially to $(\mx^*, \ml^*, \mz^*)$, where $\mx^*$ is the unique optimal solution of \cref{cenp}.
\end{corollary}

Due to the similar reason with \cref{cor1}, we omit the proof of \cref{cor2}.

\subsection{Convergence of Proj-IDEA} \label{con-p-idea}
The following lemma describes the relation between the equilibrium point of Proj-IDEA and the saddle point of $\cLa$.
\begin{lemma} \label{optCon2}
    Suppose \cref{dif} holds, $\mG$ is a strongly connected and weight-balanced digraph, and $\sum_{i=1}^nz_i(0) =\0$, then $(\mx^*, \lambda^*)$ is a saddle point of $\cLa$ iff there exists $\mz^* \in \mR^{np}$ such that $(\mw^*, \ml^*, \mz^*)$ is an equilibrium point of Proj-IDEA, where $\mw^*= \mx^*-\nabla f(\mx^*) - A\T \lambda^*$ and $\ml^* = \1_n \otimes \lambda^*$.
\end{lemma}

As a projected variant of IDEA, Proj-IDEA does inherit some properties from IDEA, for example, Proj-IDEA can also converge over undigraphs when $f_i$ is only convex, as shown in the following theorem.
\begin{theorem} \label{the5}
    Suppose \cref{dif,ug} hold. Then, given $\alpha, \ \beta > 0$, for any $(\mw(0), \ml(0), \mz(0)) \in \mR^d \times \mR^{np} \times \mR^{np}$ satisfies $\sum_{i=1}^nz_i(0) =\0$, $\mx(t)$ driven by Proj-IDEA converges to an optimal solution of \cref{cenp}.
\end{theorem}

Similar to IDEA, we have to introduce the strong convexity of $f_i$ to guarantee the convergence of Proj-IDEA over digraphs.
\begin{theorem} \label{the6}
    Suppose \cref{dif,dg} hold, $f_i$ is $\mu_i$-strongly convex on $\mathcal{X}_i$, $\forall i \in \mathcal{V}$. Define $\mu = \min_{i \in \mathcal{V}}{\mu_i}$ and $\varphi > 0$, let $\alpha$ and $\beta$ satisfy
    \eqe{ \label{stepthe6}
        \alpha &\geq \frac{(\varphi^2+3\varphi+3)\os^2(\mA)}{(\varphi+1)\mu}, \\
        \beta &\geq \frac{(\varphi+1)^2\alpha}{\varphi\eta_2(\hat{L})}.
    }
    Then, for any $(\mw(0), \ml(0), \mz(0)) \in \mR^d \times \mR^{np} \times \mR^{np}$ satisfies $\sum_{i=1}^nz_i(0) =\0$, $\mx(t)$ driven by Proj-IDEA converges to the unique optimal solution of \cref{cenp}.
\end{theorem}

\begin{remark}
    Notice that IDEA and Proj-IDEA use the gradient of $f_i$, thereby requiring that $f_i$ is differentiable. Nevertheless, it is worth emphasizing that IDEA and Proj-IDEA can still work even if $f_i$ is not differentiable, where the gradient can be replaced with the generalized gradient. Specifically, on the basis of our proof ideas, it is trivial to extend \cref{the1,the3,the5,the6} to the nondifferentiable cases by feat of nonsmooth analysis techniques \cite{clarke1990optimization,cortes2008discontinuous}.
\end{remark}

\section{Numerical Experiments}
\label{simu}
\newcolumntype{Y}{>{\centering\arraybackslash}X}
\strutlongstacks{T}

In this section, some numerical experiments are taken to validate our theoretical results and evaluate the performance of IDEA and Proj-IDEA. In particular, we compare the convergence rates of IDEA and EDEA in all numerical experiments, including their projected variants, Proj-IDEA and Proj-EDEA, where Proj-EDEA is given as
\eqe{
    \dot{\mw} &= -\alpha (\mw-\mx+\nabla f(\mx) + \mA\T \ml) - \mA\T \mr, \\
    \dot{\ml} &= \mr-\mL\ml, \\
    \dot{\mr} &= - \gamma(\mr-(\mA\mx-\mb)) -\mz - \beta\mL\mr, \\
    \dot{\mz} &= \gamma\beta\mL\mr, \\
    \mx &= \Px(\mw).
    \nonumber
}
Besides, for strongly convex cases, Unaugmented IDEA or Unaugmented Proj-IDEA is also added to the comparison.

We consider four cases, which correspond to \cref{the5} (\cref{the1}), \cref{the2}, \cref{the6} (\cref{the3}), and \cref{the4} respectively, where the first case is borrowed from \cite{falsone2020tracking} (the difference is that $\mathcal{X}_i$ may be unbounded in our case).

Case 1: undigraphs, $n = 50$, $p = 10$, $d_i = 2$, $f_i = c_i\T x_i$, $\mathcal{X}_i = \{x_i \in \mR^2 \ | \ l_{ij} \leq x_{ij} \leq r_{ij}, \ j = 1,2\}$, where $l_{ij}$ may be $-\infty$ and $r_{ij}$ may be $+\infty$, $i=1, \cdots, n, \ j = 1,2$.

Case 2: undigraphs, $n = 50$, $p = 10$, $d_i = 2$, $f_i = x_i\T a_ix_i + c_i\T x_i$, $\mathcal{X}_i = \mR^2$, $a_i \in \mR^{2\times2}$ is a positive definite diagonal matrix, $i=1, \cdots, n$, $A$ is guaranteed to have full row rank.

Case 3: digraphs, $n = 20$, $p = 4$, $d_i = 2$, $f_i = x_i\T a_ix_i + c_i\T x_i$, $\mathcal{X}_i = \{x_i \in \mR^2 \ | \ l_{ij} \leq x_{ij} \leq r_{ij}, \ j = 1,2\}$, where $l_{ij}$ may be $-\infty$ and $r_{ij}$ may be $+\infty$, $a_i \in \mR^{2\times2}$ is a positive definite diagonal matrix, $i=1, \cdots, n, \ j = 1,2$.

Case 4: Except that $A_i = I_2$ and $\mathcal{X}_i = \mR^2$, $i=1, \cdots, n$, the remainder is the same with Case 3.

We construct four undigraphs: cycle graph, Erdos–Renyi random graphs with the connectivity probability $p = 0.05, 0.1, 0.3$. The algebraic connectivities $\eta_2(L)$ of these undigraphs are incremental. There are also four digraphs: directed cycle graph, directed exponential graphs with $e = 2, 4, 6$. The directed exponential graph with $n$ nodes is constructed as: node $i$ can only send information to node $(i+2^j)\mod n$, $j=0, 1, \cdots, e-1$. Obviously, the numbers of these digraphs' edges are incremental.

Since IDEA, EDEA, and their variants are all continuous-time algorithms, we use their Euler discretizations in practical implementation. For example, the Euler discretization of IDEA is given as
\eqe{
    \mx(k+1) =& \mx(k) - \alpha \delta(\nabla f(\mx(k)) + \mA\T \ml(k)) \\
    &-\delta\mA\T (\mA\mx(k)-\mb-\mz(k)), \\
    \ml(k+1) =& \ml(k) + \delta(\mA\mx(k)-\mb)-\delta\mz(k) - \beta \delta\mL\ml(k), \\
    \mz(k+1) =& \mz(k) + \alpha\beta \delta\mL\ml(k),
    \nonumber
}
where $\delta$ is the discretized step-size. The discretized step-sizes of the four cases are chosen as $0.01$, $0.005$, $0.001$, and $0.001$ respectively. For each algorithm, we try different combinations of the parameters to obtain the fastest convergence rate. The experiment results of the four cases are illustrated in \cref{fig1,fig2,fig3,fig4}. The relative optimality gap is defined as $\frac{|f(\mx(k))-f^*|}{|f(\mx(0))-f^*|}$ and $\frac{\|\mx(k)-\mx^*\|}{\|\mx(0)-\mx^*\|}$ for convex and strongly convex cases respectively.

As shown in \cref{fig1,fig2,fig3,fig4}, IDEA usually has a faster convergence rate than EDEA, especially when the graph is poorly connected, though the communication cost of IDEA is only half that of EDEA at each iteration. Besides, the convergence rate of IDEA is slightly better than Unaugmented IDEA for most cases.

\section{Conclusion}
\label{conclusion}
This paper studys a class of distributed constraint-coupled optimization problem, i.e., \cref{cenp2}. APGD is first proposed to solve the special case of \cref{cenp2} where local constrained sets are absent, but it cannot be implemented distributedly. Benefiting from the brand-new comprehending of a classical distributed unconstrained optimization algorithm, we propose the novel implicit tracking approach and design a distributed version of APGD, i.e., IDEA. To deal with the general case of \cref{cenp2}, a projected variant of IDEA, i.e., Proj-IDEA, is further designed. With the help of the Lyapunov stability theory, the convergences of IDEA and Proj-IDEA over undigraphs and digraphs are analyzed respectively. Finally, numerical experiments are taken to corroborate our theoretical results. Compared with the results obtained under undigraphs, the conditions to guarantee the convergences of IDEA and Proj-IDEA over digrpahs is stronger, hence future work will focus on obtaining weaker convergence conditions under digraphs.

\bibliographystyle{ieeetr}
\bibliography{../public/bib}

\begin{IEEEbiography}
    [{\includegraphics[width=1in,height=1.25in,clip,keepaspectratio]{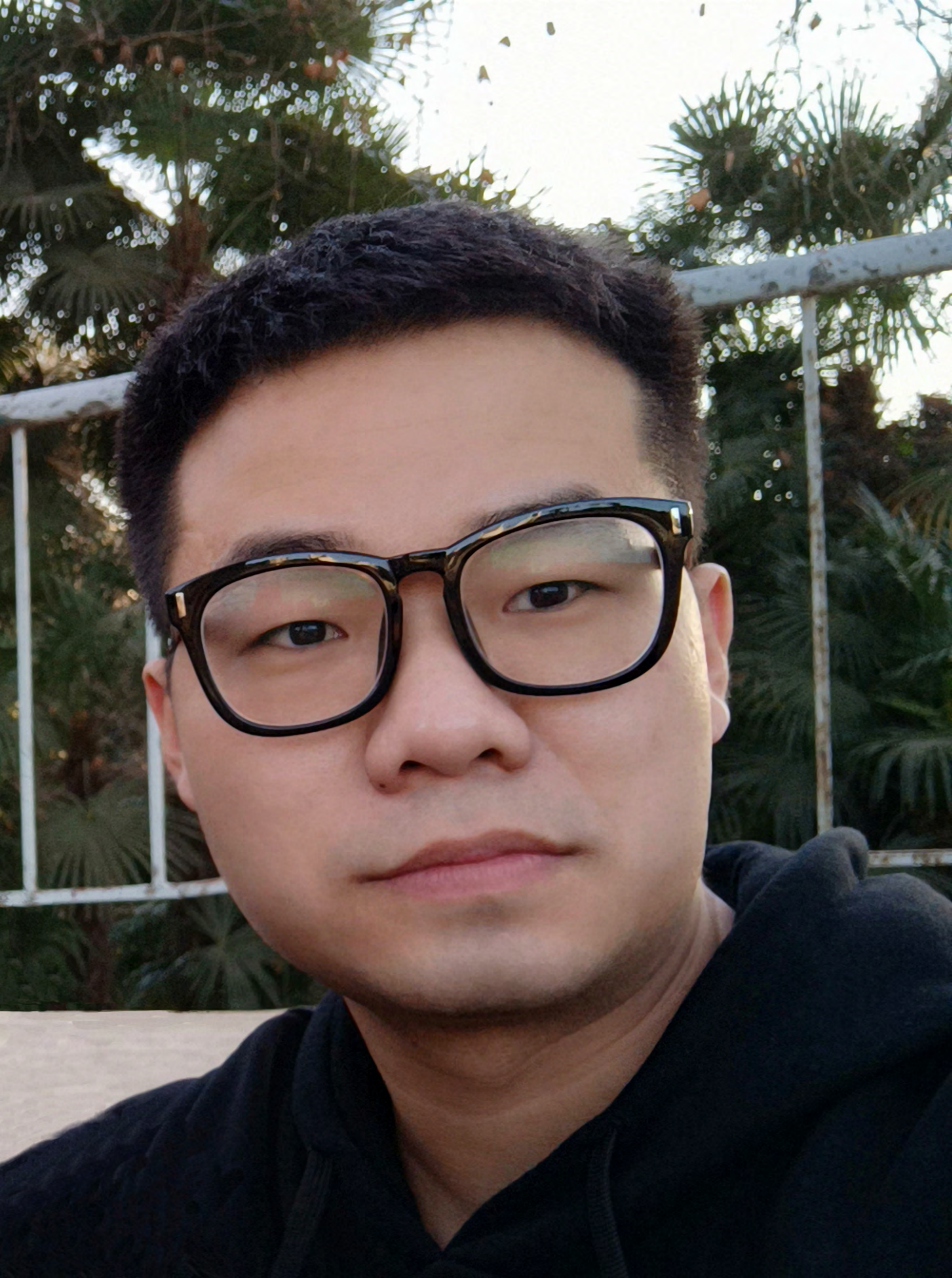}}]
    {Jingwang Li}
    received the B.S. degree in engineering management from Huazhong Agricultural University, Wuhan, China, in 2019 and the M.S. degree in control science and engineering from Huazhong University of Science and Technology, Wuhan, China, in 2022. His research interests include distributed/decentralized optimization and learning.
\end{IEEEbiography}

\begin{IEEEbiography}
    [{\includegraphics[width=1in,height=1.25in,clip,keepaspectratio]{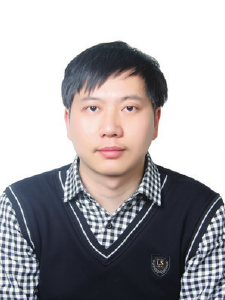}}]
    {Housheng Su}
    received his B.S. degree in automatic control and his M.S. degree in control theory and control engineering from Wuhan University of Technology, Wuhan, China, in 2002 and 2005, respectively, and his Ph.D. degree in control theory and control engineering from Shanghai Jiao Tong University, Shanghai, China, in 2008. From December 2008 to January 2010, he was a Postdoctoral researcher with the Department of Electronic Engineering, City University of Hong Kong, Hong Kong. Since November 2014, he has been a full professor with the School of Artificial Intelligence and Automation, Huazhong University of Science and Technology, Wuhan, China. He is an Associate Editor of IET Control Theory and Applications. His research interests lie in the areas of multi-agent coordination control theory and its applications to autonomous robotics and mobile sensor networks.
\end{IEEEbiography}

\begin{appendix} \label{app}

    \textit{Proof of \cref{cenCon}.} Let $(\mx^*, \lambda^*)$ be an equilibrium point of APGD, then APGD can be rewritten as
    \eqe{ \label{cenaeq}
        \dot{\mx} &= -\alpha(\nabla f(\mx)-\nabla f(\mx^*) + A\T (\lambda-\lambda^*)) - A\T A(\mx-\mx^*), \\
        \dot{\lambda} &= A(\mx - \mx^*).
    }
    Consider the following positive definite and radially unbounded candidate Lyapunov function:
    \eqe{
        V(\mx, \lambda) = \frac{1}{2}(\mx-\mx^*)\T (\mx-\mx^*) + \frac{\alpha}{2}(\lambda-\lambda^*)\T (\lambda-\lambda^*),
        \nonumber
    }
    whose Lie derivative along \cref{cenaeq} is
    \eqe{
        \dot{V} &= -\alpha(\mx-\mx^*)\T (\nabla f(\mx)-\nabla f(\mx^*)) -\|A(\mx-\mx^*)\|^2 \\
        &\leq 0,
        \nonumber
    }
    where the inequality holds since $f$ is convex. Therefore, $(\mx(t), \lambda(t))$ is bounded.

    Define $\mathcal{M} = \lt\{(\mx, \lambda) \in \mR^d \times \mR^{p} \ | \ \dot{V}(\mx, \lambda) = 0\rt\}$, let the largest invariant set of $\mathcal{M}$ and the equilibrium point set of APGD be $\mathcal{I}$ and $\mathcal{E}$ respectively, obviously $\mathcal{E}\subset \mathcal{I}$, next we will prove $\mathcal{I}\subset \mathcal{E}$. For any $(\hat{\mx}, \hat{\lambda}) \in \mathcal{I} \in \mathcal{M}$, we have
    \eqe{ \label{ic}
        A(\hat{\mx}-\mx^*) &= \0, \\
        (\hat{\mx}-\hat{\mx}^*)\T (\nabla f(\mx)-\nabla f(\mx^*)) &= 0,
    }
    combining \cref{ic} and \cref{kkt} gives that
    \eqe{
        A\hat{\mx}-b &= \0, \\
        (\hat{\mx}-\mx^*)\T (\nabla f(\hat{\mx})+A\T \lambda^*) &= 0.
        \nonumber
    }
    The convexity of $f$ implies that
    $$f(\hat{\mx})\geq f(\mx^*)+(\hat{\mx}-\mx^*)\T \nabla f(\hat{\mx}),$$
    notice that $(\mx^*, \lambda^*)$ is also a saddle point of $\cLa$, hence $\cLa(\hat{\mx}, \lambda^*) \geq \cLa(\mx^*, \lambda^*)$, it follows that
    \eqe{
        0 &= (\hat{\mx}-\mx^*)\T (\nabla f(\hat{\mx})+A\T \lambda^*) \\
        &= (\hat{\mx}-\mx^*)\T \nabla f(\hat{\mx}) + \lambda^{*\top}A(\hat{\mx}-\mx^*) \\
        &\geq f(\hat{\mx}) - f(\mx^*) + \lambda^{*\top}A(\hat{\mx}-\mx^*) \\
        &= \cLa(\hat{\mx}, \lambda^*) - \cLa(\mx^*, \lambda^*) \\
        &\geq 0,
        \nonumber
    }
    recall that $A(\hat{\mx}-\mx^*)=0$, then we have
    $$0 \leq f(\hat{\mx}) - f(\mx^*) \leq 0,$$
    which means that $f(\hat{\mx}) = f(\mx^*)$, hence $\hat{\mx}$ is an optimal solution of \cref{cenp}.
    Let $(\mx(t), \lambda(t))$ be the trajectory of APGD initiated from $(\hat{\mx}, \hat{\lambda})$, since $(\hat{\mx}, \hat{\lambda})\in \mathcal{I}$, we have $(\mx(t), \lambda(t)) \in \mathcal{I}$ for any $t \geq 0$, which means that $\mx(t)$ is always an optimal solution of \cref{cenp}. Since the strong duality holds, there must exist $\hat{\lambda}^*(t) \in \mR^p$ such that $(\mx(t), \hat{\lambda}^*(t))$ is a saddle point of $\cLa$ for any $t \geq 0$ \cite[Theorem 3.34]{ruszczynski2011nonlinear}. As a saddle point of $\cLa$, $(\mx(t), \hat{\lambda}^*(t))$ must satisfy $\nabla f(\mx(t)) + A\T \hat{\lambda}^*(t) = \0$, besides, $(\mx(t), \lambda(t)) \in \mathcal{M}$ gives that $A\mx(t) - b = \0$, hence $(\mx(t), \lambda(t))$ must be governed by the following system:
    \eqe{
        \dot{\mx}(t) &= -\alpha A\T (\hat{\lambda}-\hat{\lambda}^*(t)), \\
        \dot{\lambda}(t) &= \0,
        \nonumber
    }
    where $\nabla f(\mx(t))$ is replaced by $-A\T \hat{\lambda}^*(t)$. Furthermore, $A\mx(t) - b = \0$ implies that $A\dot{\mx}(t) = \0$, then
    $$(\hat{\lambda}-\hat{\lambda}^*(t))\T A\dot{\mx}(t) = -\alpha \|A\T (\hat{\lambda}-\hat{\lambda}^*(t))\|^2 = 0,$$
    hence $\dot{\mx}(t) = -\alpha A\T (\hat{\lambda}-\hat{\lambda}^*(t)) = \0$, which implies that $(\hat{\mx}, \hat{\lambda})$ is an equilibrium point of APGD, then we can conclude that $\mathcal{I}\subset \mathcal{E}$ since $(\hat{\mx}, \hat{\lambda})$ can be any point of $\mathcal{I}$. Therefore, the largest invariant set of $\mathcal{M}$ is the equilibrium point set of APGD.

    Let $(\mx(t), \lambda(t))$ be the trajectory of APGD initiated from any point of $\mR^d \times \mR^{p}$, according to LaSalle's invariance principle, $(\mx(t), \lambda(t))$ converges to the equilibrium point set of APGD. Since $(\mx^*, \lambda^*)$ of $V$ can be any equilibrium point of APGD, all equilibrium point of APGD are Lyaponov stable, it follows that APGD is semistable \cite[Theorem 4.20]{haddad2011nonlinear}, which implies that $(\mx(t), \lambda(t))$ converges to an equilibrium point of APGD, then the conclusion follows from the equivalence of the equilibrium point of APGD and the saddle point of $\cLa$. $\hfill\blacksquare$

    \textit{Proof of \cref{cenStr}.} Let $\mu = \min_{i \in \mathcal{V}}{\mu_i}$, obviously $f$ is $\mu$-strongly convex, which guarantees that \cref{cenp} has a unique optimal solution $\mx^*$. According to \cref{kkt}, however, the optimal solution of its dual problem, $\lambda^*$, is not unique unless $A$ has full row rank. Even though, we can easily verify that $\lambda_c^*$, the projection of $\lambda^*$ onto $\mC(A)$, is unique \cite[Lemma 2]{alghunaim2020linear}, and $(\mx^*, \lambda_c^*)$ is a saddle point of $\cLa$. Notice that $\lambda(0)=\0$, obviously it holds that $\lambda(t)\in \mC(A), \ \forall t\geq 0$. In this case, i.e., $(\mx(t), \lambda(t)) \in \mR^d \times \mC(A)$, $(\mx^*, \lambda_c^*)$ is the unique saddle point of $\cLa$. Therefore, to obtain the desired conclusion, we only need to prove that $(\mx(t), \lambda(t))$ converges exponentially to $(\mx^*, \lambda_c^*)$.

    Let $(\mx(t), \lambda(t))$ be the trajectory of APGD initiated from any point of $\mR^d \times \mR^{p}$, based on the proof of \cref{cenCon}, we know that $(\mx(t), \lambda(t))$ is bounded. Besides, recall that $\nabla f_i$ is locally Lipschitz, then we have
    \eqe{ \label{lS}
        \|\nabla f(\mx(t))-\nabla f(\mx^*)\| \leq l\|\mx-\mx^*\|, \ \forall t \geq 0,
    }
    where $l>0$ is a constant.
    Define
    \eqe{
        \mX &= \mx - \mx^*, \\
        \Lambda &= \lambda - \lambda_c^*,
        \nonumber
    }
    then APGD can be rewritten as
    \eqe{ \label{cena2}
        \dot{\mX} &= -\alpha(g(\mX, \mx^*) + A\T \Lambda) - A\T A\mX, \\
        \dot{\Lambda} &= A\mX,
    }
    where $g(\mX, \mx^*) = \nabla f(\mX+\mx^*) - \nabla f(\mx^*)$.

    Consider the following positive definite and radially unbounded candidate Lyapunov function:
    \eqe{
        V(\mx, \lambda) &= \frac{\varphi}{2\alpha}\mX\T \mX + \frac{\varphi+1}{2}\Lambda\T \Lambda \\
        &\quad + \frac{1}{2\alpha}(\mX+A\T \Lambda)\T (\mX+A\T \Lambda) \\
        &= \me\T \mE\me,
        \nonumber
    }
    where $\varphi > \max\lt\{0, \frac{l}{2}-1\rt\}$, $\me = [\mX\T, \Lambda\T]\T$, and
    $$\mE =\left[\begin{array}{cc}
                \frac{1}{2\alpha}(\varphi+1)I & \frac{1}{2\alpha}A\T                        \\
                \frac{1}{2\alpha}A            & \frac{\varphi +1}{2}I+\frac{1}{2\alpha}AA\T \\
            \end{array}\right].$$
    The positive definiteness of $V$ implies that $\mE \succ 0$, thus we have
    \eqe{ \label{eig2}
        \ue(\mE)\|\me\|^2 \leq V \leq \ove(\mE)\|\me\|^2.
    }
    Recall that $\lambda(0) = \0$, then
    \eqe{
        \Lambda(t) &= \int_0^t \lt((A\mx(s)-b)-\lambda_c^*\rt)ds \\
        &= \int_0^t \lt(A(\mx(s)-\mx^*)-\lambda_c^*\rt)ds \in \mC(A),
        \nonumber
    }
    hence
    \eqe{ \label{uL}
        \|A\T \Lambda(t)\|^2 \geq \us^2(A)\|\Lambda(t)\|^2.
    }
    The Lie derivative of $V$ along \cref{cena2} is
    \eqe{
        \dot{V} &= -(\varphi+1)\mX\T g(\mX, \mx^*) - \frac{\varphi}{\alpha}\|A\mX\|^2 - (A\T \Lambda)\T g(\mX, \mx^*) \\
        &\quad - \|A\T \Lambda\|^2 \\
        &= -(\varphi+1)\mX\T g(\mX, \mx^*) - \frac{\varphi}{\alpha}\|A\mX\|^2 - \frac{1}{2}\|A\T \Lambda\|^2 \\
        &\quad - \frac{1}{2}\|A\T \Lambda+g(\mX, \mx^*)\|^2 + \frac{1}{2}\|g(\mX, \mx^*)\|^2 \\
        &\leq -\lt(\varphi+1-\frac{l}{2}\rt)\mu\|\mX\|^2 - \frac{\us^2(A)}{2}\|\Lambda\|^2,
        \nonumber
    }
    where the inequality holds due to \cref{strPro}, \cref{lS}, and \cref{uL}. Define $\gamma = \min\lt\{\lt(\varphi+1-\frac{l}{2}\rt)\mu, \frac{\us^2(A)}{2}\rt\}$, obviously $\gamma>0$, hence we can obtain that
    \eqe{
        \dot{V} \leq -\gamma\|\me\|^2,
        \nonumber
    }
    applying \cref{eig2} gives that
    \eqe{
        \dot{V} \leq -\frac{\gamma}{\ove(\mE)}V,
        \nonumber
    }
    it follows that
    \eqe{
        V(t) \leq V(0)e^{-\frac{\gamma}{\ove(\mE)}t},
        \nonumber
    }
    we then have
    \eqe{
        \|\me(t)\| \leq \sqrt{\frac{V(0)}{\ue(\mE)}} e^{-\frac{\gamma}{2\ove(\mE)}t}.
        \nonumber
    }
    Therefore, the trajectory of \cref{cena2} converges exponentially to the origin, which completes the proof. $\hfill\blacksquare$

    \textit{Proof of \cref{optCon}.} If $(\mx^*, \ml^*, \mz^*)$ is an equilibrium point of IDEA, then
    \begin{subequations} \label{13}
        \begin{align}
            \0 & = -\nabla f(\mx^*) - \mA\T \ml^*, \label{13a} \\
            \0 & = -\mz^* + \mA\mx^*-\mb, \label{13b}          \\
            \0 & = \mL\ml^*. \label{13c}
        \end{align}
    \end{subequations}
    If $(\mx^*, \lambda^*)$ is a saddle point of $\cLa$, then
    \begin{subequations} \label{15}
        \begin{align}
            \0 & = \nabla f(\mx^*) + A\T \lambda^*, \label{15a} \\
            \0 & = A\mx^*-b. \label{15b}
        \end{align}
    \end{subequations}

    \textbf{Sufficiency: }Left multiplying \cref{13b} by $\1_n \otimes I_p$ gives that
    \eqe{ \label{113}
        A\mx^*-b = \sum_{i=1}^n(A_ix_i^*-b_i) = \sum_{i=1}^nz_i^* = \sum_{i=1}^nz_i(0) =\0,
        \nonumber
    }
    since $\mG$ is weight-balanced, we have
    \eqe{ \label{114}
        \sum_{i=1}^n\dot{z}_i = (\1_n \otimes I_p)\T \dot{\mz} = \alpha\beta (\1_n \otimes I_p)\T (L \otimes I_p)\ml = \0,
        \nonumber
    }
    it follows that
    \eqe{
        A\mx^*-b = \sum_{i=1}^nz_i^* = \sum_{i=1}^nz_i(0) =\0.
        \nonumber
    }
    Note that $\mG$ is strongly connected, hence \cref{13c} implies that $\ml^* = \1_n \otimes \lambda^*$, where $\lambda^* \in \mR^p$, we then have $\mA\T \ml^*=A\T \lambda^*$, as a result, \cref{15a} can be derived from \cref{13a}.

    \textbf{Necessity: }Let $\ml^* = \1_n \otimes \lambda^*$, obviously \cref{13c} holds and \cref{13a} can be derived from \cref{15a}. Let $\mz^* = \mA\mx^*-\mb$, \cref{13b} can be further obtained. $\hfill\blacksquare$

    \textit{Proof of \cref{the1}.} As stated in \cref{rem4}, if \cref{ug} holds, we have $L = \sqrt{L}\sqrt{L}$. With the aid of this decomposion, we can obtain the following equivalent form of IDEA:
    \eqe{ \label{sys4}
        \dot{\mx} &= -\alpha (\nabla f(\mx) + \mA\T \ml) -\mA\T (\mA\mx-\mb-\sL\my), \\
        \dot{\ml} &= \mA\mx-\mb-\sL\my- \beta \mL\ml, \\
        \dot{\my} &= \alpha\beta \sL\ml,
    }
    where $\sL = \sqrt{L} \otimes I_p$. It is obvious that IDEA can be recovered from \cref{sys4} by replacing $\sL\my$ with $\mz$. Let $(\mx^*, \ml^*, \my^*)$ be an equilibrium point of \cref{sys4}, which satisfies
    \begin{subequations} \label{es3}
        \begin{align}
            \0 & = -\nabla f(\mx^*) - \mA\T \ml^*, \label{s4a} \\
            \0 & = -\sL\my^* + \mA\mx^*-\mb, \label{s4b}       \\
            \0 & = \sL\ml^*, \label{s4c}
        \end{align}
    \end{subequations}
    obviously $(\mx^*, \ml^*, \sL\my^*)$ is an equilibrium point of IDEA, and according to \cref{optCon}, $\mx^*$ is an optimal solution of \cref{cenp}. Therefore, to obtain the desired conclusion, we only need to prove that $(\mx(t), \ml(t),\my(t))$ converges to $(\mx^*, \ml^*, \my^*)$.

    Let
    \eqe{ \label{tra2}
        \mX &= \mx - \mx^*, \\
        \mLa &= \ml - \ml^*, \\
        \mZ & = \my - \my^*,
    }
    then \cref{sys4} can be rewritten as
    \eqe{ \label{sys3}
        \dot{\mX} &= -\alpha (g(\mX, \mx^*) + \mA\T \mLa) -\mA\T (\mA\mX-\sL\mZ), \\
        \dot{\mLa} &= - \beta \mL\mLa + \mA\mX-\sL\mZ, \\
        \dot{\mZ} &= \alpha\beta \sL\mLa,
    }
    where $g(\mX, \mx^*) = \nabla f(\mX+\mx^*) - \nabla f(\mx^*)$.

    Consider the following positive definite and radially unbounded candidate Lyapunov function:
    \eqe{
        V(\mX, \mLa, \mZ) &= \frac{1}{2}\mX\T \mX + \frac{1}{2\beta}\mZ\T \mZ \\
        &\quad + \frac{1}{2\alpha}(\sL\mZ + \alpha\mLa)\T (\sL\mZ + \alpha\mLa),
        \nonumber
    }
    whose Lie derivative along \cref{sys3} is
    \eqe{
        \dot{V} &= -\alpha\mX\T g(\mX, \mx^*) - \|\mX\|^2 + 2(\mA\mX)\T \sL\mZ - \|\sL\mZ\|^2 \\
        &= -\alpha\mX\T g(\mX, \mx^*) - \|\mA\mX-\sL\mZ\|^2,
        \nonumber
    }
    the convexity of $f$ implies that
    $$\mX\T g(\mX, \mx^*)=(\mx - \mx^*)\T (\nabla f(\mX+\mx^*) - \nabla f(\mx^*))\geq0,$$
    hence $\dot{V}\leq 0$ and $(\mx(t), \ml(t),\my(t))$ is bounded.

    Define $\mathcal{M} = \left\{(\mX, \mLa, \mZ) \in \mR^d \times \mR^{np} \times \mR^{np} \ | \ \dot{V} = 0\right\}$, let the largest invariant set of $\mathcal{M}$ and the equilibrium point set of APGD be $\mathcal{I}$ and $\mathcal{E}$ respectively, obviously $\mathcal{E}\subset \mathcal{I}$, next we will prove $\mathcal{I}\subset \mathcal{E}$. For any $(\hat{\mX}, \hat{\mLa}, \hat{\mZ}) \in \mathcal{I} \in \mathcal{M}$, it holds that
    \eqe{ \label{ic2}
        \mA\hat{\mX}-\sL\hat{\mZ} &= \0, \\
        \hat{\mX}\T g(\hat{\mX}, \mx^*) &= 0,
    }
    we can easily obtain that
    \eqe{
        \0 &= (\1_n \otimes I_p)(\mA\hat{\mX}-\sL\hat{\mZ}) = A\hat{\mX} = A(\hat{\mx} - \mx^*), \\
        0 &= \hat{\mX}\T g(\hat{\mX}, \mx^*)=(\hat{\mx}-\mx^*)\T (\nabla f(\hat{\mx})-\nabla f(\mx^*)),
        \nonumber
    }
    where $\hat{\mx}=\hat{\mX}+\mx^*$. With the help of the technique used in the proof of \cref{cenCon}, it is easy to prove that $\hat{\mx}$ is an optimal solution of \cref{cenp}.

    Let $(\mX(t), \mLa(t), \mZ(t))$ be the trajectory of \cref{sys3} initiated from $(\hat{\mX}, \hat{\mLa}, \hat{\mZ})$, since $(\hat{\mX}, \hat{\mLa}, \hat{\mZ}) \in \mathcal{I}$, it holds that $(\mX(t), \mLa(t), \mZ(t)) \in \mathcal{I}$ for any $t \geq 0$. Consequently, $\mx(t) = \mX(t) + \mx^*$ is always an optimal solution of \cref{cenp} and $(\mX(t), \mLa(t), \mZ(t))$ always satisfies \cref{ic2}, which means that $(\mX(t), \mLa(t), \mZ(t))$ must be governed by the following system:
    \eqe{ \label{sys5}
        \dot{\mX}(t) &= -\alpha (g(\mX(t), \mx^*) + \mA\T \mLa(t)), \\
        \dot{\mLa}(t) &= - \beta \mL\mLa(t), \\
        \dot{\mZ}(t) &= \alpha\beta \sL\mLa(t).
    }
    There exist two cases: (i) $\hat{\mLa} = \1_n \otimes (\hat{\lambda}-\lambda^*)$; (ii) $\hat{\mLa} \neq \1_n \otimes (\hat{\lambda}-\lambda^*)$, where $\hat{\lambda} \in \mR^p$.

    For the first case, substituting $\hat{\mLa} = \1_n \otimes (\hat{\lambda}-\lambda^*)$ in \cref{sys5} gives that
    \eqe{ \label{sys7}
        \dot{\mX}(t) &= -\alpha (g(\mX(t), \mx^*) + A\T (\hat{\lambda}-\lambda^*)), \\
        \dot{\mLa}(t) &= \0, \\
        \dot{\mZ}(t) &= \0.
    }
    Since the strong duality holds and $\mx(t) $ is an optimal solution of \cref{cenp} for any $t \geq 0$, there must exist $\hat{\lambda}^*(t) \in \mR^p$ such that $(\mx(t), \hat{\lambda}^*(t))$ is a saddle point of $\cLa$ \cite[Theorem 3.34]{ruszczynski2011nonlinear}. As a saddle point of $\cLa$, $(\mx(t), \hat{\lambda}^*(t))$ must satisfy $\nabla f(\mx(t)) + A\T \hat{\lambda}^*(t) = \0$, and further applying $\nabla f(\mx^*) + A\T \lambda^* = \0$, we can obtain that
    \eqe{
        &g(\mX(t), \mx^*) + A\T (\hat{\lambda}-\lambda^*) \\
        &= \nabla f(\mx(t)) + A\T \hat{\lambda} \\
        &= A\T (\hat{\lambda}-\hat{\lambda}^*(t)),
        \nonumber
    }
    substituting it in \cref{sys5} gives that
    \eqe{
        \dot{\mX}(t) &= -\alpha A\T (\hat{\lambda}-\hat{\lambda}^*(t)).
        \nonumber
    }
    Recall that $(\mX(t), \mLa(t), \mZ(t)) \in \mathcal{I}$ for any $t\geq 0$, hence $A\mX(t) = \0$, which implies that $A\dot{\mX}(t) = \0$, it follows that
    $$0 = (\hat{\lambda}-\lambda^*(t))\T A\dot{\mX}(t) = -\alpha \|A\T (\hat{\lambda}-\lambda^*(t))\|^2,$$
    hence $\dot{\mX}(t) = \0$, which leads to the conclusion that $(\hat{\mX}, \hat{\mLa}, \hat{\mZ})$ is an equilibrium point of \cref{sys3}.

    For the second case, in light of \cref{sys5}, we can easily verify that $\lim_{t \rightarrow \infty}\mLa(t) = \1_n \otimes (\hat{\lambda}-\lambda^*)$. Therefore, based on the analysis of the first case, $(\mX(t), \mLa(t), \mZ(t))$ converges to an equilibrium point of \cref{sys3}. Notice that there is one-to-one correspondence between the equilibrium points of \cref{sys3,sys4}, without loss of generality, let $(\mx^*, \ml^*, \my^*)$ in \cref{tra2} be the equilibrium point of \cref{sys4} that corresponds to the equilibrium point of \cref{sys3} which $(\mX(t), \mLa(t), \mZ(t))$ converges to, we then have $\lim_{t \rightarrow \infty}V(t) = 0$. In consideration of the facts that $V(0) = V(\hat{\mX}, \hat{\mLa}, \hat{\mZ})$ and $\dot{V}(t) = 0$, if $(\hat{\mX}, \hat{\mLa}, \hat{\mZ})$ is not that equilibrium point of \cref{sys3} (which means that $V(0) > 0$), then a paradox would arise. Therefore, $(\hat{\mX}, \hat{\mLa}, \hat{\mZ})$ must be that equilibrium point.

    In conclusion, $(\hat{\mX}, \hat{\mLa}, \hat{\mZ})$ must be an equilibrium point of \cref{sys3}, and because $(\hat{\mX}, \hat{\mLa}, \hat{\mZ})$ can be any point of $\mathcal{I}$, we can conclude that $\mathcal{E}= \mathcal{I}$, i.e., the largest invariant set of $\mathcal{M}$ is the equilibrium point set of \cref{sys3}. Let $(\mX(t), \mLa(t), \mZ(t))$ be the trajectory of \cref{sys3} initiated from any point of $\mR^d \times \mR^{np} \times \mR^{np}$, according to LaSalle’s invariance principle, $(\mX(t), \mLa(t), \mZ(t))$ converges to the equilibrium point set of \cref{sys3}. Since $(\mx^*, \ml^*, \mz^*)$ in \cref{tra2} can be any equilibrium point of \cref{sys4}, we can easily verify that all equilibrium points of \cref{sys3} are Lyapunov stable, hence \cref{sys3} is semistable \cite[Theorem 4.20]{haddad2011nonlinear}, which implies that $(\mX(t), \mLa(t), \mZ(t))$ converges to an equilibrium point of \cref{sys3}, equivalently, $(\mx(t), \ml(t),\my(t))$ converges to an equilibrium point of \cref{sys4}, the proof is finished. $\hfill\blacksquare$

    \begin{lemma} \label{A}
        Suppose \cref{ug} holds and $A$ has full row rank, then $\mA\mA\T +c\mL \succ 0$, where $c>0$ is a constant.
    \end{lemma}

    \textit{Proof of \cref{A}.} Since $\mG$ is connected, $\mL \succeq 0$, it follows that $\mA\mA\T +c\mL \succeq 0$. Let $\mz=[z_1\T, \dots, z_n\T]\T\in \mR^{np}$, assume that $\mz\T (\mA\mA\T +c\mL)\mz = 0$, then we have $\mz\T \mA\mA\T \mz=0$ and $\mz\T \mL\mz = 0$. According to $\mz\T \mL\mz = 0$, we can obtain that $z_1 = z_2 = \cdots =z_n=z$, where $z \in \mR^p$, it follows that $\mz\T \mA\mA\T \mz = z\T \sum_{i=1}^nA_iA_i\T z=0$. $A$ has full row rank implies that $\sum_{i=1}^nA_iA_i\T  = AA\T\succ 0$, hence $z=\0$, which leads to the conclusion. $\hfill\blacksquare$

    \textit{Proof of \cref{the2}.} Since $f_i$ is strongly convex, the optimal solution of \cref{cenp} is unique. Combining with the condition that $A$ has full row rank, we can easily verify the uniqueness of the equilibrium point of IDEA based on the proof of \cref{optCon}.

    For the convenience of analysis, we still consider the equivalent form of IDEA--\cref{sys4}, and let $\my(0) = \0$. Let $(\mx^*, \ml^*, \my^*)$ be an equilibrium point of \cref{sys4}, in light of \cref{es3}, $\mx^*$ and $\ml^*$ are unique, but $\my^*$ is not since $\sL$ is singular. Even though, $\my_c^*$, the projection of $\my^*$ onto $\mC(\sL)$, is unique. Since $\my(0) \in \mC(\sL)$, we have $\my(t) \in \mC(\sL), \ \forall t\geq 0$. In this case, i.e., $(\mx(t), \ml(t),\my(t)) \in \mR^d \times \mR^{np} \times \mC(\sL)$, $(\mx^*, \ml^*, \my_c^*)$ is the unique equilibrium point of \cref{sys4}. Consequently, to obtain the desired conclusion, we only need to prove that $(\mx(t), \ml(t),\my(t))$ converges exponentially to $(\mx^*, \ml^*, \my_c^*)$, hence we still use the error system \cref{sys3}, but let $\mZ = \my - \my_c^*$.

    Consider the following positive definite and radially unbounded candidate Lyapunov function:
    \eqe{
        V(\mX, \mLa, \mZ) = V_1 + V_2 + V_3,
        \nonumber
    }
    where
    \eqe{
        V_1 = \frac{\varphi}{2}\mX\T \mX + \frac{\alpha}{2}\mLa\T \mLa + \frac{\varphi+1}{2\beta}\mZ\T \mZ,
        \nonumber
    }
    \eqe{
        V_2 = \frac{\varphi}{2\alpha}(\sL\mZ + \alpha\mLa)\T (\sL\mZ + \alpha\mLa),
        \nonumber
    }
    and
    \eqe{
        V_3 &= \frac{1}{2}\lt\|\mX+\frac{1}{\alpha}\mA\T (\sL\mZ + \alpha\mLa)\rt\|^2.
        \nonumber
    }
    The Lie derivatives of $V_1$, $V_2$, and $V_3$ along \cref{sys3} are
    \eqe{
        \dot{V}_1 &= -\varphi\alpha\mX\T g(\mX, \mx^*) - (\varphi-1)\alpha\mLa\T \mA\mX - \varphi\|\mA\mX\|^2 \\
        &\quad + \varphi(\mA\mX)\T \sL\mZ -\alpha\beta\mLa\T \mL\mLa +\varphi\alpha\mLa\T \sL\mZ,
        \nonumber
    }
    \eqe{
        \dot{V}_2 &= -\varphi\|\sL\mZ\|^2 + \varphi(\mA\mX)\T \sL\mZ - \varphi\alpha\mLa\T \sL\mZ \\
        &\quad + \varphi\alpha\mLa\T \mA\mX,
        \nonumber
    }
    and
    \eqe{
        \dot{V}_3 &= -\alpha\mX\T g(\mX, \mx^*)-\alpha\mLa\T \mA\mX - (\mA\T \sL\mZ)\T g(\mX, \mx^*) \\
        &\quad - (\mA\T \sL\mZ)\T \mA\T \mLa - \alpha(\mA\T \mLa)\T g(\mX, \mx^*) - \alpha\|\mA\T \mLa\|^2,
        \nonumber
    }
    hence the Lie derivative of $V$ along \cref{sys3} is
    \eqe{
        \dot{V} &= -(\varphi+1)\alpha\mX\T g(\mX, \mx^*) -\varphi\|\mA\mX\|^2 + 2\varphi(\mA\mX)\T \sL\mZ \\
        &\quad -\alpha\beta\mLa\T \mL\mLa -\varphi\|\sL\mZ\|^2 - (\mA\T \sL\mZ)\T g(\mX, \mx^*) \\
        &\quad - (\mA\T \sL\mZ)\T \mA\T \mLa - \alpha(\mA\T \mLa)\T g(\mX, \mx^*) -\alpha\|\mA\T \mLa\|^2 \\
        &= -(\varphi+1)\alpha\mX\T g(\mX, \mx^*) -\varphi\|\mA\mX\|^2 -\alpha\beta\mLa\T \mL\mLa \\
        &\quad  - \alpha\|\mA\T \mLa\|^2-\varphi\|\sL\mZ\|^2 +2\varphi\|\mA\mX\|^2 +\frac{\varphi}{2}\|\sL\mZ\|^2 \\
        &\quad -\varphi\lt\|\sqrt{2}\mA\mX-\frac{1}{\sqrt{2}}\sL\mZ\rt\|^2 +\frac{1}{2}\|\mA\T \sL\mZ\|^2 \\
        &\quad +\frac{1}{2}\|g(\mX, \mx^*)\|^2 -\frac{1}{2}\|\mA\T \sL\mZ+g(\mX, \mx^*)\|^2 \\
        &\quad +\frac{1}{2}\|\mA\T \sL\mZ\| +\frac{1}{2}\|\mA\T \mLa\|^2 -\frac{1}{2}\|\mA\T \sL\mZ+\mA\T \mLa\|^2 \\
        &\quad +\frac{\alpha}{2}\|\mA\T \mLa\|^2 +\frac{\alpha}{2}\|g(\mX, \mx^*)\|^2 -\frac{\alpha}{2}\|\mA\T \mLa+g(\mX, \mx^*)\|^2 \\
        &\leq -(\varphi+1)\alpha\mX\T g(\mX, \mx^*) +\varphi\|\mA\mX\|^2 -\alpha\beta\mLa\T \mL\mLa \\
        &\quad  - \frac{\alpha-1}{2}\|\mA\T \mLa\|^2 -\frac{\varphi}{2}\|\sL\mZ\|^2 +\|\mA\T \sL\mZ\|^2 \\
        &\quad +\frac{1+\alpha}{2}\|g(\mX, \mx^*)\|^2 \\
        &\leq -\lt[\lt(\varphi+1-\frac{l}{2}\rt)\alpha-\frac{l}{2}\rt]\mX\T g(\mX, \mx^*) +\varphi\os^2(\mA)\|\mX\|^2 \\
        &\quad -\frac{\alpha-1}{2}\mLa\T \left(\mA\mA\T +\frac{2\alpha\beta}{\alpha-1}\mL\right)\mLa \\
        &\quad -\lt(\frac{\varphi}{2}-\os^2(\mA)\rt)\|\sL\mZ\|^2 \\
        &\leq -\lt\{\lt[\lt(\varphi+1-\frac{l}{2}\rt)\alpha-\frac{l}{2}\rt]\mu-\varphi\os^2(\mA)\rt\}\|\mX\|^2 \\
        &\quad -\frac{\alpha-1}{2}\ue\lt(\mA\mA\T +\frac{2\alpha\beta}{\alpha-1}\mL\rt)\|\mLa\|^2 \\
        &\quad -\lt(\frac{\varphi}{2}-\os^2(\mA)\rt)\eta_2(L)\|\mZ\|^2,
        \nonumber
    }
    where the last two inequalities hold due to \cref{strPro}, \cref{A}, \cref{th2step}, and
    $$\mZ(t) = \my(t) - \my_c^* = \my(0)- \my_c^*+\int_{0}^t\alpha\beta \sL\ml(\tau)d\tau \in \mC(\sL).$$
    According to \cref{th2step}, we have
    \eqe{
        \lt[\lt(\varphi+1-\frac{l}{2}\rt)\alpha-\frac{l}{2}\rt]\mu-\varphi\os^2(\mA) > 0,
        \nonumber
    }
    define $\me=[\mX\T, \mLa\T, \mZ\T]\T$, it holds that
    \eqe{
        \dot{V} \leq -\gamma\|\me\|^2,
        \nonumber
    }
    where
    \eqe{
        \gamma &= \min\Bigg\{\lt[\lt(\varphi+1-\frac{l}{2}\rt)\alpha-\frac{l}{2}\rt]\mu-\varphi\os^2(\mA), \\
        &\quad \frac{\alpha-1}{2}\ue\lt(\mA\mA\T +\frac{2\alpha\beta}{\alpha-1}\mL\rt), \lt(\frac{\varphi}{2}-\os^2(\mA)\rt)\eta_2(L)\Bigg\} > 0.
        \nonumber
    }
    Notice that $V$ can be rewritten as
    \eqe{
        V(\mX, \mLa, \mZ) = \me\T \mE\me,
        \nonumber
    }
    where $\mE=\frac{1}{2}[\mE_1, \mE_2, \mE_3]$,
    $$\mE_1 = \left[\begin{array}{c}(\varphi+1)I \\ \mA \\ \frac{1}{\alpha}\sL\mA \\ \end{array}\right],$$
    $$\mE_2 = \left[\begin{array}{c}\mA\T\\ (\varphi+1)\alpha I + \mA\mA\T\\ \varphi\sL+\frac{1}{\alpha}\sL\mA\mA\T\\ \end{array}\right],$$
    and
    $$\mE_3 = \left[\begin{array}{c}\frac{1}{\alpha}\mA\T \sL \\ \varphi\sL+\frac{1}{\alpha}\mA\mA\T \sL \\ \frac{\varphi+1}{\beta}I + \frac{\varphi}{\alpha}\mL + \frac{1}{\alpha^2}\sL\mA\mA\T \sL \\ \end{array}\right].$$
    The positive definiteness of $V$ implies that $\mE \succ 0$, then we have
    \eqe{ \label{eig1}
        \ue(\mE)\|\me\|^2 \leq V \leq \ove(\mE)\|\me\|^2,
    }
    it follows that
    \eqe{
        \dot{V} \leq -\frac{\gamma}{\ove(\mE)}V,
        \nonumber
    }
    then
    \eqe{
        V(t) \leq V(0)e^{-\frac{\gamma}{\ove(\mE)}t},
        \nonumber
    }
    combining with \cref{eig1}, we can obtain that
    \eqe{
    \|\me(t)\| \leq \sqrt{\frac{V(0)}{\ue(\mE)}}e^{-\frac{\gamma}{2\ove(\mE)}t}.
    \nonumber
    }
    Therefore, $(\mx(t), \ml(t),\my(t))$ converges exponentially to $(\mx^*, \ml^*, \my_c^*)$, which completes the proof. $\hfill\blacksquare$

    \textit{Proof of \cref{the3}.} Let $(\mx^*, \ml^*, \mz^*)$ be an equilibrium point IDEA, according to \cref{optCon}, $\mx^*$ is an optimal solution of \cref{cenp}. The $\mu$-strong convexity of $f_i$ can guarantee the uniquenesses of $\mx^*$ and $\mz^*$, however, the uniqueness of $\ml^*$ cannot be guaranteed since we do not assume that $A$ has have full row rank. As a result, the equilibrium point of IDEA is not necessarily unique.

    Define $\mathbf{m} = m \otimes I_p$, $\mathbf{M} = M \otimes I_p$, and $T = [\mathbf{m}, \mathbf{M}]$, where $m$ and $M$ are defined in \cref{mM}. In light of the definitions of $m$ and $M$, we can easily verify that $T\T T = TT\T  = I_{np}$. For the convenience of analysis, we apply the following variables transformation to IDEA:
    \eqe{ \label{tra}
        \mX &= \mx - \mx^*, \\
        \mLa &= T\T (\ml - \bm{\lambda}^*), \\
        \mZ & = T\T (\mz - \mz^*),
    }
    where $(\mx^*, \ml^*, \mz^*)$ is an equilibrium point of IDEA. According to \cref{optCon}, $\mx^*$ is the unique optimal solution of \cref{cenp}. Let $\mLa = [\mLa_1\T, \mLa_2\T]\T$, where $\mLa_1 \in \mR^{p}$ and $\mLa_2 \in \mR^{(n-1)p}$, similarly, let $\mZ = [\mZ_1\T, \mZ_2\T]\T$. Notice that $\mathbf{m}\T\mL = \0$, $\mathbf{m}\T\mz = \0$, $\mathbf{m}\T (\mA\mx^*-\mb) = \0$, and $\nabla f(\mx^*) = -\mA\T \bm{\lambda}^*$, it follows that $\mZ_1(t) = \0$, then IDEA can be rewritten as
    \eqe{ \label{sys}
        \dot{\mX} \ &= - \alpha(g(\mX, \mx^*) + \mA\T T\mLa) -\mA\T (\mA\mX - \mathbf{M}\mZ_2), \\
        \dot{\mLa}_1 &= \mathbf{m}\T \mA\mX, \\
        \dot{\mLa}_2 &= \mathbf{M}\T \mA\mX - \mZ_2 - \beta \mathbf{M}\T \mL\mathbf{M}\mLa_2, \\
        \dot{\mZ}_1 &= \0, \\
        \dot{\mZ}_2 &= \alpha\beta \mathbf{M}\T \mL\mathbf{M}\mLa_2,
    }
    where $g(\mX, \mx^*) = \nabla f(\mX+\mx^*) - \nabla f(\mx^*)$.

    Consider the following positive definite and radially unbounded candidate Lyapunov function:
    \eqe{
        V(\mX, \mLa, \mZ_2) &= \frac{1}{2\alpha}(\alpha\mLa_2+\mZ_2)\T (\alpha\mLa_2+\mZ_2) + \frac{\varphi}{2}\alpha\mLa_2\T \mLa_2 \\
        &\quad + \frac{\varphi+1}{2}\alpha\mLa_1\T \mLa_1 + \frac{\varphi+1}{2}\mX\T\mX,
        \nonumber
    }
    let $\hat{\mL} =  \hat{L} \otimes I_p$, according to \cref{sec}, we can obtain that
    $$\mLa_2^{\top}\mathbf{M}^{\top} \mL\mathbf{M} \mLa_2 = \mLa_2^{\top}\mathbf{M}^{\top} \hat{\mL}\mathbf{M} \mLa_2 \geq \eta_2(\hat{L})\|\mLa_2\|^2,$$
    hence the Lie derivative of $V$ along \cref{sys} is
    \eqe{ \label{dotv}
        \dot{V} &= -(\varphi+1)\alpha\mLa_2\T \mZ_2 + (\varphi+2)\mZ_2\T \mathbf{M}\T \mA\mX - \mZ_2\T\mZ_2 \\
        &\quad - \varphi\alpha\beta \mLa_2\T \mathbf{M}\T\mL\mathbf{M} \mLa_2 - (\varphi+1)\alpha\mX\T g(\mX, \mx^*) \\
        &\quad - (\varphi+1)\mX\T \mA\T \mA\mX \\
        &\leq -\lt\|(\varphi+1)\alpha\mLa_2+\frac{1}{2}\mZ_2\rt\|^2 + (\varphi+1)^2\alpha^2\mLa_2^2 + \frac{1}{4}\|\mZ_2\|^2 \\
        &\quad - \lt\|(\varphi+2) \mathbf{M}\T \mA\mX - \frac{1}{2}\mZ_2\rt\|^2 + (\varphi+2)^2\|\mathbf{M}\T \mA\mX\|^2 \\
        &\quad + \frac{1}{4}\|\mZ_2\|^2 -\|\mZ_2\|^2 -\varphi\alpha\beta\eta_2(\hat{L})\|\mLa_2\|^2 \\
        &\quad - (\varphi+1)\alpha\mu\|\mX\|^2 - (\varphi+1)\|\mA\mX\|^2 \\
        &\leq -\frac{1}{2}\|\mZ_2\|^2 -(\varphi\alpha\beta\eta_2(\hat{L}) - (\varphi+1)^2\alpha^2)\|\mLa_2\|^2 \\
        &\quad -\mX\T ((\varphi+1)\alpha\mu I -(\varphi^2+3\varphi+3)\mA\T \mathbf{M}\mathbf{M}\T \mA)\mX,
    }
    where the first and the last inequalities hold because of \cref{strPro} and $I - \mathbf{M}\mathbf{M}\T  = \mathbf{m}\mathbf{m}\T\succeq 0$ respectively.

    According to \cref{115}, we have
    \eqe{
        (\varphi+1)\alpha\mu I -(\varphi^2+3\varphi+3)\mA\T \mathbf{M}\mathbf{M}\T \mA &\succ 0, \\
        \varphi\alpha\beta\eta_2(\hat{L}) - (\varphi+1)^2\alpha^2 &> 0,
        \nonumber
    }
    hence $\dot{V} \leq 0$, which implies the trajectory of \cref{sys} is bounded. Define $\mathcal{M} = \left\{(\mX, \mLa, \mZ) \in \mR^d \times \mR^{np} \times \mR^{np} \ | \ \dot{V} = 0\right\}$, obviously the largest invariant set of $\mathcal{M}$ is the equilibrium point set of \cref{sys}. According to LaSalle’s invariance principle, the trajectory of \cref{sys} converges to its equilibrium point set. Notice that the variables transformation \cref{tra} is affine and static, there is one-to-one correspondence between the equilibrium points of \cref{sys} and IDEA. Since $(\mx^*, \ml^*, \mz^*)$ in \cref{tra} can be any equilibrium point of IDEA, we can easily verify that all equilibrium points of \cref{sys} are Lyapunov stable, which implies that \cref{sys} is semistable \cite[Theorem 4.20]{haddad2011nonlinear}, hence the trajectory of \cref{sys} converges to one of its equilibrium points. Equivalently, the trajectory of IDEA converges to one of its equilibrium points, the proof is finished. $\hfill\blacksquare$

    \textit{Proof of \cref{the4}.} Let $(\mx^*, \ml^*, \mz^*)$ be an equilibrium point of IDEA, according to \cref{optCon}, $\mx^*$ is an optimal solution of \cref{cenp}. Notice that $f_i$ i strongly convex and $A_i = I_p$, hence $\mx^*$, $\ml^*$, and $\mz^*$ are all unique, which means that $(\mx^*, \ml^*, \mz^*)$ is the unique equilibrium point of IDEA. For the convenience of analysis, we also use the variables transformation \cref{tra}.

    Define $\me = [\mX\T, \mLa_1\T, \mLa_2\T, \mZ_2\T]\T$, consider the following positive definite and radially unbounded candidate Lyapunov function:
    \eqe{
        V(\mX, \mLa, \mZ_2) &= \frac{1}{2\alpha}(\alpha\mLa_2+\mZ_2)\T (\alpha\mLa_2+\mZ_2) + \frac{\varphi}{2}\alpha\mLa_2\T \mLa_2 \\
        &\quad + \frac{\varphi+1}{2}\alpha\mLa_1\T \mLa_1 + \frac{\varphi+1}{2}\mX\T\mX \\
        &\quad + \underbrace{\frac{1}{2\alpha}(\mX+\mathbf{m}\mLa_1)\T (\mX+\mathbf{m}\mLa_1)}_{V_1} \\
        &= \me\T \mE\me,
        \nonumber
    }
    where $\mE = $
    $$\frac{1}{2}\left[\begin{array}{cccc}
                (\varphi+1+\frac{1}{\alpha})I & \frac{1}{\alpha}\mathbf{m}            & \0                  & \0                \\
                \frac{1}{\alpha}\mathbf{m}\T  & ((\varphi+1)\alpha+\frac{1}{\alpha})I & \0                  & \0                \\
                \0                            & \0                                    & (\varphi+1)\alpha I & I                 \\
                \0                            & \0                                    & I                   & \frac{1}{\alpha}I
            \end{array}\right].$$
    The positive definiteness of $V$ implies that $\mE \succ 0$, then we have
    \eqe{ \label{111}
        \ue(\mE)\|\me\|^2 \leq V \leq \ove(\mE)\|\me\|^2.
    }
    The Lie derivative of $V_1$ along \cref{sys} is
    \eqe{
        \dot{V}_1 &= -\mX\T g(\mX, \mx^*) - \mX\T \mathbf{m}\mLa_1 - \mX\T \mathbf{M}\mLa_2 -\frac{1}{\alpha}\mX\T \mathbf{M}\mathbf{M}\T \mX \\
        &\quad +\frac{1}{\alpha}\mX\T \mathbf{M}\mZ_2 - \mLa_1\T \mathbf{m}\T g(\mX, \mx^*) - \mLa_1\T \mLa_1 \\
        &= -\|\mX+\frac{1}{2}\mathbf{m}\mLa_1\|^2 + \|\mX\|^2 + \frac{1}{4}\|\mLa_1\|^2 \\
        &\quad -\frac{1}{2}\|\mX+\mathbf{M}\mLa_2\|^2 + \frac{1}{2}\|\mX\|^2 + \frac{1}{2}\|\mLa_2\|^2 \\
        &\quad -\frac{1}{\alpha}\|\mathbf{M}\T \mX-\frac{1}{2}\mZ_2\|^2 + \frac{1}{\alpha}\|\mathbf{M}\T \mX\|^2 + \frac{1}{4\alpha}\|\mZ_2\|^2 \\
        &\quad -\|g(\mX, \mx^*) + \frac{1}{2}\mathbf{m}\mLa_1\|^2 + \|g(\mX, \mx^*)\|^2 + \frac{1}{4}\|\mLa_1\|^2 \\
        &\quad -\mX\T g(\mX, \mx^*) - \frac{1}{\alpha}\|\mathbf{M}\T \mX\|^2 - \|\mLa_1\|^2 \\
        &\leq -\lt(\mu - l^2 -\frac{3}{2}\rt)\|\mX\|^2 - \frac{1}{2}\|\mLa_1\|^2 + \frac{1}{2}\|\mLa_2\|^2 \\
        &\quad + \frac{1}{4\alpha}\|\mZ_2\|^2,
        \nonumber
    }
    where \cref{mM} is used in the above derivation and the inequality holds due to \cref{strPro}. Notice that $I - \mathbf{M}\mathbf{M}\T\succeq 0$ and recall \cref{dotv}, we can obtain that
    \eqe{
        \dot{V} &\leq -\lt[\lt((\varphi+1)\alpha+1\rt)\mu -(\varphi^2+3\varphi+3) -l^2-\frac{3}{2}\rt]\|\mX\|^2 \\
        &\quad -\lt[\varphi\alpha\beta\eta_2(\hat{L}) - (\varphi+1)^2\alpha^2 -\frac{1}{2}\rt]\|\mLa_2\|^2 \\
        &\quad - \frac{1}{2}\|\mLa_1\|^2 -\frac{1}{2}\lt(1-\frac{1}{2\alpha}\rt)\|\mZ_2\|^2 \\
        &\leq -\gamma\|\me\|^2,
        \nonumber
    }
    where
    \eqe{
        \gamma &= \min\Bigg\{\lt((\varphi+1)\alpha+1\rt)\mu -(\varphi^2+3\varphi+3) -l^2-\frac{3}{2}, \\
        &\quad \varphi\alpha\beta\eta_2(\hat{L}) - (\varphi+1)^2\alpha^2 -\frac{1}{2}, \ \frac{1}{2}\lt(1-\frac{1}{2\alpha}\rt)\Bigg\}.
        \nonumber
    }
    Note that \cref{116} guarantees that $\gamma>0$ and recall \cref{111}, we have
    \eqe{
        \dot{V} \leq -\frac{\gamma}{\ove(\mE)}V,
        \nonumber
    }
    it follows that
    \eqe{
        V(t) \leq V(0)e^{-\frac{\gamma}{\ove(\mE)}t},
        \nonumber
    }
    applying \cref{111} gives that
    \eqe{
        \|\me(t)\| \leq \sqrt{\frac{V(0)}{\ue(\mE)}} e^{-\frac{\gamma}{2\ove(\mE)}t},
        \nonumber
    }
    hence the trajectory of \cref{sys} converges exponentially to the origin. Since the variables transformation \cref{tra} is affine and static, the trajectory of IDEA converges to its unique equilibrium point, the proof is finished. $\hfill\blacksquare$

    \textit{Proof of \cref{optCon2}.} If $(\mw^*, \ml^*, \mz^*)$ is an equilibrium point of Proj-IDEA, then
    \begin{subequations} \label{pideaE}
        \begin{align}
            \0 & = \mw^*-\mx^*+\nabla f(\mx^*) + \mA\T \ml^*, \label{piEa} \\
            \0 & = \mA\mx^*-\mb^*-\mz^*, \label{piEb}                      \\
            \0 & = \mL\ml^*. \label{piEc}
        \end{align}
    \end{subequations}
    If $(\mx^*, \lambda^*)$ is a saddle point of $\cLa$, then
    \begin{subequations} \label{15p}
        \begin{align}
            \0 & \in \nabla f(\mx^*) + A\T \lambda^* + N_{\mathcal{X}}(\mx^*), \label{15pa} \\
            \0 & = A\mx^*-b, \label{15pb}
        \end{align}
    \end{subequations}
    where $N_{\mathcal{X}}(\mx^*)$ is the normal cone to $\mathcal{X}$ at $\mx^*$.
    Note that \cref{piEc} implies that $\ml^* = \1_n \otimes \lambda^*$, where $\lambda^* \in \mR^p$, hence \cref{piEa} can be rewritten as
    \eqe{
        \0 = \mw^*-\mx^*+\nabla f(\mx^*) + A\T \lambda^*,
        \nonumber
    }
    it follows that
    \eqe{ \label{mxx}
        \mx^* = \Px(\mw^*) = \Px(\mx^*-\nabla f(\mx^*) - A\T \lambda^*),
        \nonumber
    }
    which is equivalent to \cref{15pa}, according to \cite[Lemma 2.38]{ruszczynski2011nonlinear}. Conversely, we can obtain \cref{piEa} from \cref{15pa} by letting
    $$\mw^*= \mx^*-\nabla f(\mx^*) - A\T \lambda^*.$$
    The remainder of this proof is the same with \cref{optCon}, which is omitted thereby. $\hfill\blacksquare$

    \textit{Proof of \cref{the5}.} For the convenience of analysis, we study the following equivalent form of Proj-IDEA:
    \eqe{ \label{sys7}
        \dot{\mw} &= -\alpha (\mw-\mx+\nabla f(\mx) + \mA\T \ml) -\mA\T (\mA\mx-\mb-\sL\my), \\
        \dot{\ml} &= \mA\mx-\mb-\sL\my - \beta \mL\ml, \\
        \dot{\my} &= \alpha\beta \sL\ml, \\
        \mx &= \Px(\mw),
    }
    obviously Proj-IDEA can be obtained by letting $\mz = \sL\my$.
    Let $(\mw^*, \ml^*, \my^*)$ be an equilibrium point of \cref{sys7}, which satisfies
    \eqe{ \label{sys7E}
        \0 &= \mw^*-\mx^*+\nabla f(\mx^*) + \mA\T \ml^*, \\
        \0 &= \mA\mx^*-\mb-\sL\my^*, \\
        \0 &= \alpha\beta \sL\ml^*,
    }
    obviously $(\mw^*, \ml^*, \sL\my^*)$ is an equilibrium point of Proj-IDEA. In light of \cref{optCon2}, $\mx^* = \Px(\mw^*)$ is an optimal solution of \cref{cenp}. Therefore, to obtain the desired conclusion, we only need to prove that the trajectory of \cref{sys7} converges to an optimal solution of \cref{cenp}. Applying \cref{sys7E}, \cref{sys7} can be rewritten as
    \eqe{ \label{sys8}
        \dot{\mw} &= -\alpha \big(\mw-\mw^*-(\mx-\mx^*)+\nabla f(\mx)-\nabla f(\mx^*) \\
        &\quad + \mA\T (\ml-\ml^*)\big) -\mA\T \big(\mA(\mx-\mx^*)-\sL(\my-\my^*)\big), \\
        \dot{\ml} &= \mA(\mx-\mx^*)-\sL(\my-\my^*) - \beta \mL(\ml-\ml^*), \\
        \dot{\my} &= \alpha\beta \sL(\ml-\ml^*), \\
        \mx &= \Px(\mw).
    }

    Consider the following candidate Lyapunov function:
    \eqe{
        V(\mw, \ml, \my) =& \underbrace{\frac{1}{2}\lt(\|\mw-\mx^*\|^2-\|\mw-\Px(\mw)\|^2\rt)}_{V_1} \\
        & + \frac{1}{2\alpha}\|\sL(\my-\my^*) + \alpha(\ml-\ml^*)\|^2 \\
        & + \frac{1}{2\beta}\|\my-\my^*\|^2,
        \nonumber
    }
    \cref{prd} implies that $V_1 \geq \frac{1}{2}\|\mx-\mx^*\|^2$, hence $V$ is positive definite and radially unbounded. \cref{strPro,proj} gives that
    \eqe{
        (\mx-\mx^*)^{\top}(\nabla f(\mx)-\nabla f(\mx^*)) &\geq 0, \\
        (\mx-\mx^*)^{\top}(\mw-\mw^*-(\mx-\mx^*)) &\geq 0,
    }
    then the Lie derivative of $V$ along \cref{sys8} is
    \eqe{
        \dot{V}
        &= -\alpha(\mx-\mx^*)\T (\mw-\mw^*-(\mx-\mx^*)) \\
        &\quad -\alpha(\mx-\mx^*)\T (\nabla f(\mx)-\nabla f(\mx^*)) - \|\mA(\mx-\mx^*)\|^2 \\
        &\quad -\|\sL(\my-\my^*)\|^2 + 2(\mA(\mx-\mx^*))\T \sL(\my-\my^*) \\
        &= -\alpha(\mx-\mx^*)\T (\mw-\mw^*-(\mx-\mx^*)) \\
        &\quad -\alpha(\mx-\mx^*)\T (\nabla f(\mx)-\nabla f(\mx^*)) \\
        &\quad - \|\mA(\mx-\mx^*)-\sL(\my-\my^*)\|^2 \\
        &\leq 0,
        \nonumber
    }
    where the first equality holds due to \cref{prd}. Therefore, $\mx(t)$, $\ml(t)$, and $\my(t)$ are all bounded. Recall that $\nabla f$ is locally Lipschitz, thus $\nabla f(\mx(t))$ is also bounded. According to $\dot{\mw}(t)$, we have $\|\mw(t)\| \leq \|\mw(0)\|e^{-t} + (1-e^{-t})M$, where $M>0$ is an upper bound related to $\mx(t)$, $\ml(t)$, and $\my(t)$, hence $\mw(t)$ is also bounded, which implies that the trajectory of \cref{sys7} is bounded.

    Define $\mathcal{M} = \lt\{(\mw, \ml, \my) \in \mR^d \times \mR^{np} \times \mR^{np} \ | \ \dot{V} = 0\rt\}$ and let the largest invariant set in $\mathcal{M}$ be $\mathcal{I}$, for any $(\hat{\mw}, \hat{\ml}, \hat{\my}) \in \mathcal{M}$, we have
    \eqe{ \label{ic3}
        \mA(\hat{\mx}-\mx^*)-\sL(\hat{\my}-\my^*) &= \0, \\
        (\hat{\mx}-\mx^*)\T (\nabla f(\hat{\mx})-\nabla f(\mx^*)) &= 0, \\
        (\hat{\mx}-\mx^*)\T (\hat{\mw}-\mw^*-(\hat{\mx}-\mx^*)) &= 0,
    }
    where $\hat{\mx} = \Px(\hat{\mw})$.
    It follows that
    \eqe{
        \0 &= (\1_n \otimes I_p)\T (\mA(\hat{\mx}-\mx^*)-\sL(\hat{\my}-\my^*)) \\
        &= A(\hat{\mx}- \mx^*) = A\hat{\mx} - b,
        \nonumber
    }
    then we have
    \eqe{
        \ml^{*\top}\mA(\hat{\mx}-\mx^*) = \lambda^{*\top}\sum_{i=1}^nA_i(\hat{x}_i-x_i^*) = \lambda^{*\top}(A\hat{\mx}- b)= 0,
        \nonumber
    }
    combining it with \cref{ic3} gives that
    \eqe{
        0 &= (\hat{\mx}-\mx^*)\T (\hat{\mw}-\mw^*-(\hat{\mx}-\mx^*)+\nabla f(\hat{\mx})-\nabla f(\mx^*) \\
        &\quad +\mA\ml^*) \\
        &= (\hat{\mx}-\mx^*)\T (\hat{\mw}-\hat{\mx}+\nabla f(\hat{\mx})).
        \nonumber
    }
    According to \cref{proj}, we have
    $$(\hat{\mx}-\mx^*)\T (\hat{\mw}-\hat{\mx})\geq 0,$$
    hence
    $$(\hat{\mx}-\mx^*)\T \nabla f(\hat{\mx}) \leq 0,$$
    it follows that
    \eqe{
        f(\hat{\mx}) - f(\mx^*) \leq (\hat{\mx}-\mx^*)\T \nabla f(\hat{\mx}) \leq 0,
        \nonumber
    }
    which holds since $f$ is convex. Recall that $\hat{\mx} \in \mathcal{X}$, $A\hat{\mx} - b = \0$, and $\mx^*$ is an optimal solution of \cref{cenp}, obviously $f(\hat{\mx}) \geq f(\mx^*)$, which implies that $f(\hat{\mx}) = f(\mx^*)$, then we can conclude that $\hat{\mx}$ is also an optimal solution of \cref{cenp}. Therefore, for any $(\mw, \ml, \my) \in \mathcal{I}$, $\mx = \Px(\mw)$ is an optimal solution of \cref{cenp}.

    Let $\mX = (\mw, \ml, \my) \in \mR^d \times \mR^{np} \times \mR^{np}$, based on the above analysis, for any initial point $\mX(0)$, the trajectory of \cref{sys7}, i.e., $\mX(t)$, is bounded, then there must exist a strictly increasing sequence $\{t_n\}$, which satisfies $\lim_{k \rightarrow \infty}t_k = \infty$, such that $\lim_{k \rightarrow \infty}\mX(t_k) = \hat{\mX}$, where $\hat{\mX} = (\hat{\mw}, \hat{\ml}, \hat{\my})$. Obviously $\hat{\mX} \in \mathcal{I}$, hence $\hat{\mx} = \Px(\hat{\mw})$ is an optimal solution of \cref{cenp}.
    Define
    \eqe{
        V_2(\mX) =& \frac{1}{2}\lt(\|\mw-\hat{\mx}\|^2-\|\mw-\Px(\mw)\|^2\rt) \\
        & + \frac{1}{2\alpha}\lt\|\sL(\my-\hat{\my}) + \alpha(\ml-\hat{\ml})\rt\|^2 \\
        & + \frac{1}{2\beta}\|\my-\hat{\my}\|^2,
        \nonumber
    }
    we can easily obtain that $\dot{V}_2 \leq 0$, hence $V_2$ is nonincreasing with respect to $t$ and bounded below, then we have $\lim_{t \rightarrow \infty}V_2(\mX(t)) = c$, where $c \geq 0$ is a constant. Note that $\lim_{k \rightarrow \infty}V_2(\mX(t_k)) = 0$, hence $c = 0$, which implies that $\lim_{t \rightarrow \infty}\mX(t) = \hat{\mX}$, it follows that $\lim_{t \rightarrow \infty}\mx(t) = \hat{\mx}$, thereby finishing the proof. $\hfill\blacksquare$

    \textit{Proof of \cref{the6}.} Let $(\mw^*, \ml^*, \mz^*)$ be an equilibrium point of Proj-IDEA, according to \cref{optCon2} and the strong convexity of $f$, we can verify that $\mx^* = \Px(\mw^*)$ is the unique optimal solution of \cref{cenp}. Recall \cref{pideaE}, Proj-IDEA can be rewritten as
    \eqe{ \label{sys9}
        \dot{\mw} &= -\alpha \big(\mw-\mw^*-(\mx-\mx^*)+\nabla f(\mx)-\nabla f(\mx^*) \\
        &\quad + \mA\T (\ml-\ml^*)\big) -\mA\T \big(\mA(\mx-\mx^*)-(\mz-\mz^*)\big), \\
        \dot{\ml} &= \mA(\mx-\mx^*)-(\mz-\mz^*) - \beta \mL(\ml-\ml^*), \\
        \dot{\mz} &= \alpha\beta \mL(\ml-\ml^*), \\
        \mx &= \Px(\mw).
    }

    Consider the following positive definite and radially unbounded candidate Lyapunov function:
    \eqe{
        V(\mw, \ml, \mz) &= \frac{\varphi+1}{2}\lt(\|\mw-\mx^*\|^2-\|\mw-\Px(\my)\|^2\rt) \\
        &\quad +\frac{\varphi\alpha}{2}\|\ml-\ml^*\|^2 + \frac{1}{2\alpha}\|\alpha(\ml-\ml^*)+\mz-\mz^*\|^2,
        \nonumber
    }
    whose Lie derivative of $V$ along \cref{sys9} is
    \eqe{
        \dot{V} \leq &-(\varphi+1)\alpha(\mx-\mx^*)\T (\nabla f(\mx)-\nabla f(\mx^*)) \\
        &-(\varphi+1)\|\mA(\mx-\mx^*)\|^2 + (\varphi+2)(\mz-\mz^*)\T \mA(\mx-\mx^*) \\
        &-(\varphi+1)\alpha(\ml-\ml^*)\T (\mz-\mz^*) \\
        &- \varphi\alpha\beta(\ml-\ml^*)\T \mL(\ml-\ml^*)- \|\mz-\mz^*\|^2,
        \nonumber
    }
    which holds due to \cref{proj}.
    Let $\ml_C$ and $\ml_N$ be the projections of $\ml$ onto $\mC(\hat{\mL})$ and $\mN(\hat{\mL})$ respectively, where $\hat{\mL} =  \hat{L} \otimes I_p$, we then have $\ml = \ml_C + \ml_N$ since $\mC(\hat{\mL}) = \mC(\hat{\mL}\T ) \perp \mN(\hat{\mL})$. Furthermore, it follows that $\ml_N = \1_n \otimes \lambda$, where $\lambda \in \mR^p$, then we have
    \eqe{
        \dot{V} \leq &-(\varphi+1)\alpha\mu\|\mx-\mx^*\|^2 -(\varphi+1)\|\mA(\mx-\mx^*)\|^2 \\
        & -\|(\varphi+2)\mA(\mx-\mx^*)-\frac{1}{2}(\mz-\mz^*)\|^2 \\
        &+ (\varphi+2)^2\|\mA(\mx-\mx^*)\|^2 + \frac{1}{4}\|\mz-\mz^*\|^2 \\
        &-(\varphi+1)\alpha\ml_C\T (\mz-\mz^*)- \varphi\alpha\beta\eta_2(\hat{L})\|\ml_C\|^2 - \|\mz-\mz^*\|^2 \\
        \leq &-((\varphi+1)\alpha\mu-(\varphi^2+3\varphi+3)\os^2(\mA))\|\mx-\mx^*\|^2 \\
        &-\|(\varphi+1)\alpha\ml_C+\frac{1}{2}(\mz-\mz^*)\|^2 + (\varphi+1)^2\alpha^2\|\ml_C\|^2 \\
        &+ \frac{1}{4}\|\mz-\mz^*\|^2 - \varphi\alpha\beta\eta_2(\hat{L})\|\ml_C\|^2 - \frac{3}{4}\|\mz-\mz^*\|^2 \\
        \leq &-((\varphi+1)\alpha\mu-(\varphi^2+3\varphi+3)\os^2(\mA))\|\mx-\mx^*\|^2 \\
        &- (\varphi\alpha\beta\eta_2(\hat{L})-(\varphi+1)^2\alpha^2)\|\ml_C\|^2 - \frac{1}{2}\|\mz-\mz^*\|^2,
        \nonumber
    }
    where the first inequality holds due to the strong convexity of $f$, $\sum_{i=1}^nz_i(t) = \sum_{i=1}^nz_i(0) =\0$, $\mL(\ml_N-\ml^*) = \0$, and $\ml_C\T \mL\ml_C = \ml_C\T \hat{\mL}\ml_C \geq \eta_2(\hat{L})\|\ml_C\|^2$. Recall \cref{stepthe6}, it follows that $\dot{V} \leq 0$, similar to the proof of \cref{the5}, we can conclude that the trajectory of \cref{sys9} is bounded. Define $\mathcal{M} = \left\{(\mw, \ml, \my) \in \mR^d \times \mR^{np} \times \mR^{np} \ | \ \dot{V} = 0\right\}$ and let $\mathcal{I}$ be the largest invariant set in $\mathcal{M}$, obviously $\mx = \Px(\mw) = \mx^*$, for any $(\mw, \ml, \my) \in \mathcal{I}$. The remainder of this proof is similar to the proof of \cref{the5}, hence we omit it. $\hfill\blacksquare$

\end{appendix}

\end{document}